\documentclass[11pt, reqno]{amsart}

\usepackage{amssymb,mathrsfs,graphicx}
\usepackage{amsmath,amsfonts,amscd,amsthm}
\usepackage{mathtools}
\usepackage{epsfig}
\usepackage{comment}
\usepackage{caption}
\usepackage{subcaption}
\usepackage{colortbl}
\usepackage{bm}
\usepackage{tikz}
\usepackage{tikz-cd}
\usetikzlibrary{matrix,shapes,arrows,positioning, chains}
\usepackage[unicode,hypertexnames=false]{hyperref}
\hypersetup{colorlinks,breaklinks,linkcolor=blue,urlcolor=blue,anchorcolor=blue,citecolor=blue}

\usepackage[shortlabels]{enumitem}
\usepackage[autostyle]{csquotes}
   \MakeOuterQuote{"}

\title[Relaxation dynamics of Winfree oscillators on Riemannian manifolds]{Relaxation dynamics of Winfree oscillators\\on Riemannian manifolds}

\author[Ha]{Seung-Yeal Ha}
\address[Seung-Yeal Ha]{\newline Department of Mathematical Sciences and Research Institute of Mathematics \newline Seoul National University, Seoul 08826, Republic of Korea}
\email{syha@snu.ac.kr}

\author[Lee]{Chaejoo Lee}
\address[Chaejoo Lee]{\newline Department of Mathematical Sciences\newline Seoul National University, Seoul 08826, Republic of Korea}
\email{leechaeju@snu.ac.kr}

\author[Lee]{Eunjun Lee}
\address[Eunjun Lee]{\newline Department of Mathematical Sciences\newline Seoul National University, Seoul 08826, Republic of Korea}
\email{eunjun.lee@snu.ac.kr}

\author[Lee]{Jaemoon Lee}
\address[Jaemoon Lee]{\newline Department of Mathematical Sciences\newline Seoul National University, Seoul 08826, Republic of Korea}
\email{dlwoans0001@snu.ac.kr}

\author[Ryoo]{Seung-Yeon Ryoo}
\address[Seung-Yeon Ryoo]{Hausdorff Center for Mathematics,
\newline University of Bonn
Endenicher Allee 60
53115 Bonn, Germany}
\email{sryoo@uni-bonn.de}

\newtheorem{theorem}{Theorem}[section]
\newtheorem{lemma}{Lemma}[section]
\newtheorem{corollary}{Corollary}[section]
\newtheorem{proposition}{Proposition}[section]
\newtheorem{remark}{Remark}[section]

\newcommand{\bbc}{\mathbb C}
\newcommand{\bbs}{\mathbb S}

\newcommand*\di{\mathop{}\!\mathrm{d}}

\begin{document}
\tikzstyle{block} = [rectangle, draw, 
text width=15em, text centered, rounded corners, minimum height=3em]
\tikzstyle{line} = [draw, -latex']

\date{\today}

\subjclass{34C40, 34C15, 34D06} 
\keywords{Winfree model, synchronization, Riemannian manifolds, unitary group, phase-locking, exponential convergence}

\thanks{\textbf{Acknowledgment.} The work of S.-Y. Ha was supported by the NRF grant (RS-2025-00514472). The last author thanks Otte Hein\"avaara for the proof of Lemma \ref{L4.1}.}

\begin{abstract}	
We study exponential relaxation to equilibrium for Winfree-type models on Riemannian manifolds.
For the sphere model with entrainment point $p$ and canonical influence, we derive  convergence to equilibrium, whenever no oscillator starts at the antipode $-p$.
For the corresponding canonical model on $\mathrm U(d)$, we obtain the same convergence whenever no initial matrix has the eigenvalue $-1$. Both results hold above explicit coupling thresholds depending on the initial configuration and the frequency matrices.
For a broader class of unitary influence functions, the same spectral condition yields exponential convergence for sufficiently large coupling strengths.
Finite-time entrance into positively invariant stable regions and contraction estimates yield exponential convergence rates and equilibrium uniqueness within those regions.
A finite symmetry construction on Lie groups realizes the unitary model and a special unitary analogue as invariant Kuramoto--Lohe subsystems, yielding local relaxation results.
Finally, we formulate an abstract Winfree model on Riemannian manifolds and establish local exponential relaxation in a fixed coordinate neighborhood under quantitative attraction and regularity assumptions.
\end{abstract}

\maketitle


\section{Introduction}\label{sec:1}
\setcounter{equation}{0}
Synchronization can be roughly defined as a type of collective behavior in populations of interacting phase oscillators, in which phases or frequencies form a consensus. It appears across a remarkably wide range of physical and biological complex systems, from networks of neurons \cite{ermentrout1994introduction,breakspear2010generative} and cardiac pacemaker cells \cite{de1983synchronization,ypey1982synchronization} to Josephson junction arrays \cite{wiesenfeld1996new,wiesenfeld1998frequency} and flashing fireflies \cite{smith1935synchronous,buck1938synchronous,buck1988synchronous}, with recent applications in signal recovery and alignment problems in data science \cite{bandeira2015convex,bandeira2015ten} and machine learning \cite{zhang2020autosync}, etc. A central mathematical challenge is to develop rigorous qualitative and quantitative frameworks to understand when spontaneous synchronization occurs and how rapidly it sets in.

The mathematical study of synchronization has been dominated by two phase models: one introduced by Arthur Winfree in 1967 \cite{winfree1967biological} and stated formally below in \eqref{A-1}, and the other introduced in 1975 by Yoshiki Kuramoto \cite{kuramoto1975international}, where oscillators on the circle interact through a sinusoidal coupling based on their phase differences. A nonexhaustive selection of references is \cite{ermentrout1990oscillator,ariaratnam2001phase,oukil2017synchronization,oukil2019invariant} for the former, and \cite{strogatz2000kuramoto,acebron2005kuramoto,dorfler2014synchronization} for the latter. A central distinction is that the Kuramoto model is invariant under a common phase shift, whereas the Winfree model generally lacks this $U(1)$-symmetry.

Synchronization on higher-dimensional manifolds is a recent topic of interest. This generalization is motivated by two factors. First, many physical systems of interest, including models of coupled quantum oscillators \cite{ha2016collective}, multi-dimensional biological rhythms \cite{winfree1980geometry}, and transformer dynamics \cite{rigollet2026mean}, are naturally described by phases taking values in higher-dimensional spaces rather than the unit circle ${\mathbb S}^1$. Second, synchronization models often interact richly with the geometry and topology of the phase space; for example, the topology of the phase space can lead to multistability \cite{markdahl2020synchronization}. The Kuramoto model and its high-dimensional counterparts are being actively studied in the literature. Some of these studies are devoted to the analogues of the Kuramoto model on spheres or Hilbert spaces \cite{olfati2006swarms,lohe2010quantum,choi2014complete,zhang2018equilibria,chandra2019continuous,geshkovski2025mathematical}, Lie groups \cite{lohe2009non,ha2016emergencelohe,ha2017emergent,deville2019synchronization,ryoo2025asymptotic}, hyperbolic geometries \cite{chen2017hyperbolic,lipton2021kuramoto}, and matrix/tensor spaces \cite{lohe2019systems,ha2020lohe}. However, to the authors' knowledge, few consensus models on Riemannian manifolds incorporate Winfree-type interactions; we are aware only of \cite{ha2021collective,park2021generalization,ha2026emergent}. As such, in this paper, we suggest another viewpoint on some high-dimensional generalizations of the $\mathbb{S}^1$-Winfree model. To set up the stage, we begin with the Winfree model. 

Let $\theta_i = \theta_i(t)$ be the time-dependent phase of the $i$-th Winfree oscillator with a real-valued natural frequency $\nu_i$. For definiteness, we choose a real-valued influence function $I=I(\theta)$ and sensitivity function $S = S(\theta)$:
\begin{equation}\notag 
I(\theta) = 1 + \cos \theta, \quad S(\theta) = -\sin \theta. 
\end{equation}
In this setting, the temporal dynamics of $\theta_i$ is governed by the following Cauchy problem for the (classical) Winfree model:
\begin{equation}
	\begin{cases} \label{A-1}
		\displaystyle \dot{\theta_i} =\nu_i - \frac{\kappa}{N} \sum^{N}_{j=1}  (1 + \cos \theta_j) \sin \theta_i, \quad t > 0, \\
		\displaystyle \theta_i \Big|_{t = 0+} = \theta_i^0,  \quad  i \in [N]\coloneqq  \{1, \cdots, N \},
	\end{cases}
\end{equation}
where $\kappa\in\mathbb{R}$ 
is the nonnegative coupling strength, and $\{\theta_i^0\}_{i=1}^N$ is the initial phase configuration. The influence function $1+\cos\theta_j$ is sometimes replaced by a smooth approximation of the Dirac delta function. The global well-posedness of \eqref{A-1} follows from the standard Cauchy--Lipschitz theory.

Ha, Park, and Ryoo \cite{ha2015emergence} established sufficient conditions for the ensemble in \eqref{A-1} to converge to an equilibrium configuration in the large-coupling regime. Their result, recalled in Theorem~\ref{T2.1}, gives an explicit coupling threshold in terms of $\max_{i\in[N]}|\nu_i|$ and a quantitative exponential convergence rate. Related results concern locally coupled networks \cite{ha2016emergent} and partial phase-locking \cite{ha2017emergence,oukil2017synchronization,oukil2019invariant}; in the latter group of works, synchronization is formulated through the control of phase differences and should be distinguished from convergence to a stationary equilibrium. In this work, we address the following questions: 
	\begin{itemize}
\item "What are the natural high-dimensional counterparts of the classical $\mathbb {S}^1$-Winfree model \eqref{A-1}?"
\item "For such models, under what conditions on system parameters and initial data do they admit a unique equilibrium in a specified trapping region, and what does the relaxation dynamics look like?"
\end{itemize}
These questions were already addressed in the literature. A high-dimensional realization of the Winfree model on the unit Euclidean sphere ${\mathbb S}^{d-1}\subset{\mathbb R}^d$ was first studied by Park \cite{park2021generalization}, where a sufficient framework for relaxation dynamics was established. Inspired by this, the authors of this paper further extended Park's sphere Winfree model to the special orthogonal groups $\mathrm{SO}(d)$ and not only established similar convergence phenomena but also developed a leader-follower mechanism \cite{ha2026emergent}. 

The earlier sphere and special orthogonal analyses \cite{park2021generalization,ha2026emergent} use localized influence functions and include leader--follower conditions requiring an oscillator initially near the entrainment point. Their admissible initial configurations can already include followers outside a hemisphere or with rotation angles exceeding $\pi/2$. For the canonical influence functions considered here, we instead establish relaxation without requiring a nearby leader: every initial configuration avoiding the antipodal point, or the spectral value $-1$, respectively, converges for sufficiently large coupling depending on that configuration. There is also a need for a flexible theory of Winfree-type interactions on manifolds, allowing for many types of influence and sensitivity functions. Note that this was established in the one-dimensional, $\mathbb{S}^1$-setting by \cite{ha2015emergence}.

\subsection*{Main results of this paper}
The main results are threefold. First, for the Winfree sphere model with $I(x)=1+\langle x,p\rangle$, we show that every initial configuration satisfying $x_i^0\neq-p$ for all $i\in[N]$ converges to an equilibrium when the coupling is sufficiently large. Second, we introduce a Winfree matrix model on $\mathrm U(d)$ with canonical influence $I(U)=1+d^{-1}\operatorname{Re}\operatorname{tr}U$ and verify the analogous result whenever $-1\notin\operatorname{spec}(U_i^0)$ for every $i\in[N]$. In both cases, for each fixed initial configuration and frequency tuple, there exists a threshold $\kappa_0$ such that convergence holds for every $\kappa>\kappa_0$. The limiting equilibrium generally depends on $\kappa$ and is unique in the specified stable region; global uniqueness is not asserted. Third, we treat general influence functions on unitary groups and formulate an abstract Winfree model on Riemannian manifolds, for which we establish a local relaxation theorem in a fixed coordinate neighborhood. We also derive Winfree--Lohe subsystems by the finite symmetry reduction of Kuramoto--Lohe ensembles; this recovers the unitary model and supplies a special unitary example of the abstract framework.

In Sections~\ref{sec:3}--\ref{sec:5}, the proofs proceed through finite-time entrance into a positively invariant stable region, contraction between trajectories in that region, and comparison with time translates to obtain an equilibrium and exponential convergence. In Section~\ref{sec:7}, the initial data already lie in an invariant coordinate neighborhood, and contraction holds from time zero. In every case, hypotheses and convergence rates are quantitative.

\subsection*{Organization of the rest of the paper}

The rest of the paper is organized as follows. After establishing notation, in Section \ref{sec:2} we review the models and results of \cite{ha2015emergence,park2021generalization,ha2026emergent}.

In Section~\ref{sec:3}, we study the Winfree model on the unit sphere $\mathbb{S}^{d-1}$, where
 the state of each oscillator is a unit vector $x_i\in\mathbb{S}^{d-1}$ and the entrainment point is a
 fixed $p\in\mathbb{S}^{d-1}$. The relevant coupling threshold involves the restricted operator
 seminorm $\|\mathfrak A p\|_{\infty}$, defined in Section~\ref{subsec:2-2}, which measures
 the part of the frequency matrices $A_i$ that enters the angular dynamics relative to $p$. Our main
 result, Theorem~\ref{T3.1}, is a higher-dimensional counterpart of Theorem~\ref{T2.1} that
 establishes finite-time entrance and exponential convergence to an equilibrium that is unique in the specified stable region, with an explicit rate.

In Section~\ref{sec:4}, we study the Winfree model on the unitary group $\mathrm{U}(d)$,
 which can be viewed as the ``quantum'' generalization of \eqref{A-1}. Here the phases are
 unitary matrices $U_i\in\mathrm{U}(d)$, and the natural influence function involves the normalized real
 part of the trace. We prove Theorem~\ref{T4.1}, a $\mathrm{U}(d)$-analogue of
 Theorem~\ref{T2.1}, in which the role of the phase $\theta_i$ is played by the eigenvalues
 of $U_i$ and convergence is measured in Schatten--von Neumann norms. A key technical ingredient
 is Lemma~\ref{L4.1}, which controls the real part of the trace of a difference $U-V$
 in terms of the Schatten $1$-norm via a spectral separation argument.

In Section~\ref{sec:5}, we extend the unitary group analysis to allow a 
general influence function $I\colon \mathrm{U}(d)\to\mathbb{R}$ satisfying natural positivity, 
monotonicity, and regularity conditions. The canonical choice $I(U)=1+\frac{1}{d}\operatorname{Re}\operatorname{tr}U$ 
recovers the setting of Section~\ref{sec:4}. The general result, Theorem~\ref{T5.1}, 
gives explicit exponential rates in terms of the Lipschitz constant of $I$ on the relevant stable 
region.

In Section~\ref{sec:6}, we start from the Winfree-to-Kuramoto embedding of \cite{morenoearle2026relaxation} and develop a symmetry reduction in the Lie-group setting of \cite{ryoo2025asymptotic}. We identify an averaging condition that makes the reduced interaction separable and verify it for unitary and special unitary groups. The convergence theory of \cite{ryoo2025asymptotic} is again applicable, and yields exponential relaxation of the reduced systems to a unique local equilibrium.

In Section~\ref{sec:7}, we formulate a Winfree model on a general Riemannian 
manifold $\mathcal M$ and prove Theorem~\ref{T6.1}, which gives conditions on the drift vector 
fields $F_i$, the influence function $I$, and the entrainment vector field $S$ under which 
global-in-time well-posedness and exponential convergence to a unique equilibrium in the specified coordinate neighborhood hold. The sphere, unitary, and special unitary models correspond to special cases of the abstract equation considered in Section~\ref{sec:7}. The quantitative hypotheses of Theorem~\ref{T6.1} are local and
chart-dependent.

\vspace{0.5cm}

\noindent {\bf Gallery of Notation}:
For Cartesian vector components, superscripts denote component indices and subscripts denote particle indices; superscripts $0$ and $\infty$ denote initial and limiting states. Thus, $x_i^k$ denotes the $k$-th Cartesian component of the $i$-th particle:
\[ x_i=(x_i^1,\ldots,x_i^d)\in{\mathbb R}^d,\quad X\coloneqq(x_1,\ldots,x_N)\in{\mathbb R}^{dN}. \]
For $x=(x^1,\cdots,x^d)\in\mathbb{R}^d$ and $Y=(y_1,\cdots,y_N)\in(\mathbb{R}^d)^N$, we introduce the standard $\ell^p$-norm and the mixed norm:
\begin{equation}\label{A-2}
\hspace{-1cm}
\|x\|_p\coloneqq 
\begin{cases}
\displaystyle\left(\sum_{k=1}^d|x^k|^p\right)^{1/p},&1\le p<\infty,\\
\displaystyle\max_{1\le k\le d}|x^k|,&p=\infty,
\end{cases}
\quad
\|Y\|_{q,p}\coloneqq 
\begin{cases}
\displaystyle\left(\frac{1}{N}\sum_{i=1}^N\|y_i\|_p^q\right)^{1/q},&1\le q<\infty,\\
\displaystyle\max_{i\in[N]}\|y_i\|_p,&q=\infty.
\end{cases}
\end{equation}
For notational simplicity, we set $\|x\|\coloneqq\|x\|_2$. For a real-valued function $f$, we denote its upper right Dini derivative by
\[
D^+f(t)\coloneqq\limsup_{h\searrow0}\frac{f(t+h)-f(t)}{h}.
\]
Whenever the ordinary right derivative
\[
f'_+(t)\coloneqq\lim_{h\searrow0}\frac{f(t+h)-f(t)}{h}
\]
exists, it coincides with $D^+f(t)$.

\section{Previous models and relaxation results}\label{sec:2}
\setcounter{equation}{0}
In this section, we briefly review the relaxation dynamics of the Winfree model and its higher-dimensional generalizations on the unit sphere and the special orthogonal group, respectively.
\subsection{The Winfree circle model} \label{subsec:2-1}
In this subsection, we consider an ensemble of Winfree oscillators on the complete network, with each pair interacting through the same normalized coupling strength. More precisely, we consider the Cauchy problem \eqref{A-1}. 

For $\alpha\in \left(0,\pi\right)$, let $\alpha^\infty\in \left(0,\frac{\pi}{3}\right]$ be the unique solution to the following trigonometric equation:
\begin{equation} \label{B-2}
\sin\alpha^\infty(1+\cos\alpha^\infty)=\sin\alpha(1+\cos\alpha). 
\end{equation}
Note that for $\alpha \in (0, \frac{\pi}{3}]$, we have 
\[ \alpha^{\infty} = \alpha. \]
Throughout the paper, we set 
\[ 
\Theta\coloneqq(\theta_1,\ldots,\theta_N),\quad \mathcal{V}\coloneqq(\nu_1,\ldots,\nu_N), \quad \|\mathcal{V}\|_\infty\coloneqq\max_{i\in[N]}|\nu_i|.
\]
\begin{theorem}{\cite{ha2015emergence}} \label{T2.1}
Suppose system parameters and initial data satisfy 
\[
\alpha\in(0,\pi)\setminus\left\{\frac{\pi}{3}\right\},\quad
\kappa>\frac{\|\mathcal{V}\|_\infty}{\sin\alpha(1+\cos\alpha)},\quad
\theta_i^0\in[-\alpha,\alpha],\quad i\in[N],
\]
and let $\Theta$ be a global solution to \eqref{A-1}. Then, the following assertions hold. 
\begin{enumerate}
\item
\emph{(Existence and uniqueness of equilibrium)}:~there exists a unique equilibrium $\Theta^\infty\in (-\alpha^\infty,\alpha^\infty)^N$.
\vspace{0.1cm}
\item
\emph{(Finite-time entrance)}:~there exists $T\geq0$ such that
\[
T \le \frac{\pi}{\kappa\sin\alpha(1+\cos\alpha)-\|\mathcal{V}\|_\infty},
\quad
\Theta(t)\in(-\alpha^\infty,\alpha^\infty)^N,\quad t\geq T.
\]
\vspace{0.1cm}
\item
\emph{(Exponential relaxation)}:~set
\[
\mu_{\alpha^\infty}\coloneqq(2\cos\alpha^\infty-1)(1+\cos\alpha^\infty)>0.
\]
Then the solution $\Theta$ relaxes to $\Theta^\infty$ exponentially fast:
    \[
    \|\Theta(t)-\Theta^\infty\|_{1}\le e^{-\kappa\mu_{\alpha^\infty}(t-T)}\|\Theta(T)-\Theta^\infty\|_{1},\quad t\ge T.
    \]
\end{enumerate}
\end{theorem} 
\begin{remark}\label{R2.1} We make two observations about the preceding theorem.
   \begin{itemize}
    \item [(1)] If $\Theta^0\in (-\alpha^\infty,\alpha^\infty)^N$, we can take $T=0$.
    \vspace{0.1cm}
    \item [(2)] If $\kappa>\frac{4}{3\sqrt{3}}{\|\mathcal{V}\|_\infty},$ we may find $\alpha\in (\frac{\pi}3,\pi)$ such that 
    \[\kappa>\frac{{\|\mathcal{V}\|_\infty}}{\sin\alpha(1+\cos\alpha)},\]
    and Theorem \ref{T2.1} tells us that system \eqref{A-1} has a stable equilibrium in $(-\frac{\pi}{3},\frac{\pi}{3})^N$.
    \end{itemize} 
\end{remark}
\subsection{The Winfree sphere model} \label{subsec:2-2}
In this subsection, we consider the Winfree model on the unit sphere \cite{park2021generalization}. Let $d\geq2$, and let $x_i = x_i(t) \in {\mathbb S}^{d-1} \subset {\mathbb R}^d$ be the state of the $i$-th Winfree oscillator at time $t$. Then, the dynamics of $x_i$ is governed by the Cauchy problem for the Winfree sphere model on ${\mathbb S}^{d-1}$:
\begin{equation} \label{B-3}
\begin{cases}
\displaystyle \dot{x}_i = A_i x_i+\frac \kappa {N} \sum_{j=1}^N I(x_j) (p-\langle x_i,p\rangle x_i),\quad t>0, \\
\displaystyle x_i \Big|_{t = 0+} = x_i^0, \quad i \in [N],
\end{cases}
\end{equation}
where each $A_i$ is a $d\times d$ skew-symmetric matrix, $p\in\mathbb{S}^{d-1}$ is a fixed unit vector, and $I:\mathbb{S}^{d-1}\to\mathbb{R}$ is the influence function. We write
\[
\mathfrak A\coloneqq(A_1,\ldots,A_N),\quad X\coloneqq(x_1,\ldots,x_N).
\]
We define the restricted operator seminorm relative to the fixed vector $p$ by
\begin{equation*}
	 \|\mathfrak Ap\|_{\infty} \coloneqq
	\max_{i\in[N]}\|A_ip\|,
\end{equation*}
and remark that for a linear map $L:\mathbb{R}^d\to\mathbb{R}^d$,
\[
\|Lp\|\leq\|L\|_{\mathrm{op}}.
\]
Moreover, if $L$ is skew-symmetric, we have the identity
\begin{equation}\label{B-4-1}
\|Lp\|=
\sup_{\substack{V\leq\mathbb{R}^d\\ \dim V=2,\ p\in V}}
\ \sup_{\substack{x,y\in V\\ \|x\|,\|y\|\leq1}}
|\langle Lx,y\rangle|.
\end{equation}
Indeed, the compression\footnote{The \emph{compression} of a linear operator $L$ to a subspace $V$ is the operator $P_V\circ L|_V:V\to V$, where $P_V$ denotes the orthogonal projection onto $V$.} of $L$ to $\operatorname{span}\{p,v\}$, where $v\perp p$ and $\|v\|=1$, is a two-by-two skew-symmetric operator with operator norm $|\langle Lp,v\rangle|$. Taking the supremum over $v$ gives the identity, since $Lp\perp p$. Next, we introduce the order parameter measuring the mean influence of the ensemble:
\begin{equation}\notag
R\coloneqq \frac 1N \sum_{j=1}^N I(x_j).
\end{equation}
With this notation, system \eqref{B-3} can be written in the mean-field form:
\begin{equation}\notag
\dot{x}_i= A_i x_i+\kappa R(p-\langle x_i,p\rangle x_i), \quad i \in [N].
\end{equation}
Note that relaxation results for the sphere model \eqref{B-3} were established in \cite[Section 4]{park2021generalization} under hypotheses involving a localized influence function and a nearby leader. In Section~\ref{sec:3}, we study the canonical influence function $I(x)=1+\langle x,p\rangle$ and prove a counterpart of Theorem~\ref{T2.1} for every initial configuration avoiding $-p$, provided the coupling is sufficiently large.
\subsection{\texorpdfstring{The Winfree matrix model on $\mathrm{SO}(n)$}{The Winfree matrix model on SO(n)}}\label{subsec:2-3}
In this subsection, we recall the Winfree matrix model on the special orthogonal group and its relaxation dynamics established in \cite{ha2026emergent}. 

Let $R_i=R_i(t)\in\mathrm{SO}(n)$ be the state of the $i$-th oscillator at time $t$, and let $\Omega_i\in\mathfrak{so}(n)$ be its natural frequency matrix. Then, the corresponding Cauchy problem reads as follows.
\begin{equation}\label{B-5}
\begin{cases}
\displaystyle \dot R_i=\Omega_iR_i+\frac{\kappa}{2N}\sum_{j=1}^N I(R_j)(I_n-R_i^2),\quad t>0,\\
\displaystyle R_i\big|_{t=0}=R_i^0\in\mathrm{SO}(n),\quad i\in[N],
\end{cases}
\end{equation}
where $I_n$ is the $n\times n$ identity matrix and $I:\mathrm{SO}(n)\to[0,\infty)$ is the influence function. For $A,B\in\mathbb R^{n\times n}$, set
\[
\langle A,B\rangle\coloneqq\frac12\operatorname{tr}(A^\top B),\quad \|A\|\coloneqq\sqrt{\langle A,A\rangle},
\]
and let $d$ be the induced geodesic distance on $\mathrm{SO}(n)$. We write
\[
B_r\coloneqq\{R\in\mathrm{SO}(n):d(I_n,R)<r\}.
\]
The influence function is assumed to be radial, namely
\[
I(R)=\widetilde I\bigl(d(I_n,R)\bigr),
\]
where $\widetilde I$ is nonnegative, decreasing, globally Lipschitz continuous, satisfies 
\[ \lim_{r\to0+}\widetilde I(r)=1, \]
and has support $[0,\beta]$ for some $0<\beta<\frac{\pi}{2}$. For $\mathcal R=(R_1,\ldots,R_N)$, we define
\[
\|\mathcal R\|_{\ell^1}\coloneqq\sum_{i=1}^N\|R_i\|,\quad
\operatorname{Lip}^{*}(I)\coloneqq\sup_{\substack{R,S\in\mathrm{SO}(n)\\R\neq S}}\frac{|I(R)-I(S)|}{\|R-S\|}.
\]
Let $\gamma_0$ and $\gamma$ satisfy
\[
0<\gamma_0<2\sin\frac{\beta}{2},\quad \gamma\coloneqq2\arcsin\frac{\gamma_0}{2}<\beta,
\]
and we set
\begin{align*}
\begin{aligned}
& \lambda_1 \coloneqq\kappa\bigl(\cos\gamma\,\widetilde I(\gamma)-\gamma\operatorname{Lip}^{*}(I)\bigr),\quad
\Gamma_i(\gamma)
\coloneqq\left(2+2\sqrt{1-\frac{N^2\|\Omega_i\|^2}{\kappa^2\widetilde I(\gamma)^2}}\right)^{1/2},\\
& \kappa_c(\gamma,1) \coloneqq\frac{N\max_{j\in[N]}\|\Omega_j\|}{\sin\bigl(2\sin(\gamma/2)\bigr)\widetilde I(\gamma)}.
\end{aligned}
\end{align*}
Then, the emergent dynamics of \eqref{B-5} can be summarized as follows.
\begin{theorem}\cite[Corollary 3.1, Theorem 3.1, and Theorem 4.1]{ha2026emergent}\label{T2.2}
Suppose that there exists $i_*\in[N]$ such that
\[
R_{i_*}^0\in B_{\gamma_0},\quad R_i^0\in B_{\Gamma_i(\gamma)},\quad i\in[N],\quad \kappa>\kappa_c(\gamma,1),\quad \lambda_1>0,
\]
and let $\mathcal R$ be a global solution to \eqref{B-5}. Then, the following assertions hold.
\begin{enumerate}
\item \emph{(Finite-time entrance)}:~there exists $t_0\geq0$ such that $\mathcal R(t)\in(B_\gamma)^N$ for all $t\geq t_0$.
\item \emph{(Existence and uniqueness of equilibrium)}:~there exists a unique equilibrium $\mathcal R^\infty\in(B_\gamma)^N$.
\item \emph{(Exponential relaxation)}:~the solution relaxes to $\mathcal R^\infty$ exponentially fast:
\[
\|\mathcal R(t)-\mathcal R^\infty\|_{\ell^1}\leq e^{-\lambda_1(t-t_0)}\|\mathcal R(t_0)-\mathcal R^\infty\|_{\ell^1},\quad t\geq t_0.
\]
\end{enumerate}
\end{theorem}
\begin{remark}
For the identical-frequency ensemble $\Omega_i=\Omega$, the same work further verifies the exponential complete state synchronization and characterizes the equilibrium configurations through a fixed-point equation for the mean influence.
\end{remark}

\section{Winfree oscillators on the unit spheres}\label{sec:3}
\setcounter{equation}{0}
In this section, we establish quantitative relaxation for the Winfree sphere model \eqref{B-3} with the explicit influence function $I(x)=1+\langle x,p\rangle$. For this, we first establish a trapping region, then prove exponential stability in that region, and finally compare a solution with its time translates to obtain the existence and uniqueness of an equilibrium and exponential relaxation toward it.
 
\subsection{Trapping and finite-time entrance} \label{subsec:3-1}
Consider the Cauchy problem for the Winfree sphere model:
\begin{equation}\label{C-1}
\begin{cases}
\displaystyle \dot{x}_i = A_i x_i+\frac{\kappa}{N} \sum_{j=1}^N \Big (1+\langle x_j, p\rangle \Big)(p-\langle x_i, p \rangle x_i),\quad t>0,\\
\displaystyle x_i \Big|_{t = 0} =x_i^0, \quad \| x_i^0 \| = 1, \quad  i \in [N].
\end{cases}
\end{equation}
The well-posedness of \eqref{C-1} follows from local existence results together with the norm-preserving property established in the following lemma. This ensures that trajectories remain on a compact manifold and hence exist globally in time. 
\begin{lemma} \label{L3.1}
Let $X=(x_1,\ldots,x_N)$ be the solution to \eqref{C-1} on its maximal interval of existence $[0,T_{\max})$. Then the following identities hold. For $i \in [N]$, 
\begin{equation} \label{C-2}
\langle A_i x_i, x_i \rangle = 0, \quad  \| x_i(t) \| = 1, \quad \| p - \langle p, x_i \rangle x_i \|^2 = 1 - \langle p, x_i \rangle^2, \quad 0\le t<T_{\max}. 
\end{equation}
\end{lemma}
\begin{proof} 
In the sequel, we will provide the proofs for each estimate. \newline

\noindent (i)~We use $A^{T}_i = -A_i$ to find 
\[ 
\langle A_i x_i, x_i \rangle = \langle x_i, A_i^T x_i \rangle = -\langle x_i, A_i x_i \rangle = -\langle A_i x_i, x_i \rangle,
\]
which implies
\begin{equation}  \label{C-3}
 \langle A_i x_i, x_i \rangle = 0.
\end{equation} 
\noindent (ii)~ We set 
\begin{equation} \label{C-4}
R = \frac{1}{N} \sum_{j=1}^N \Big (1+\langle x_j,p\rangle \Big),
\end{equation}
and rewrite $\eqref{C-1}_1$ as 
\begin{equation} \label{C-5}
\dot{x}_i = A_i x_i+ \kappa R (p-\langle x_i, p \rangle x_i).
\end{equation}
Then, we use \eqref{C-3} and \eqref{C-5} to find
\begin{align*}
\begin{aligned}
\frac{d}{dt} \|x_i \|^2 &=  2\langle {\dot x}_i, x_i \rangle = 2 \langle A_i x_i+ \kappa R (p-\langle x_i, p \rangle x_i), x_i \rangle   \\
&= 2 \langle A_i x_i, x_i \rangle +  2 \kappa R \langle p, x_i \rangle  \Big( 1 - \langle x_i, x_i \rangle  \Big) = 2 \kappa R \langle p, x_i \rangle  \Big( 1 - \| x_i \|^2  \Big).
\end{aligned}
\end{align*}
This yields
\[ 
\frac{d}{dt}  \Big( 1 - \| x_i \|^2  \Big) =  -2 \kappa R \langle p, x_i \rangle  \Big( 1 - \| x_i \|^2  \Big).
\]
Using Gronwall's lemma and the initial condition $1 - \|x_i^0 \|^2 = 0$, we obtain the desired identity:
\[
1 - \| x_i(t) \|^2 = (1 - \| x^0_i \|^2) \exp \Big(  -2 \kappa \int_0^t R(s) \langle p, x_i(s) \rangle \di s     \Big) = 0, \quad 0\le t<T_{\max}.
\]
\noindent (iii)~We use $\|p\| = \| x_i \|^2 = 1$ to see
\begin{align*}
\begin{aligned} 
 \| p - \langle p, x_i \rangle x_i \|^2 &= \langle p - \langle p, x_i \rangle x_i, p - \langle p, x_i \rangle x_i \rangle       \\
&= \langle p, p \rangle - \langle p, x_i \rangle \langle p,  x_i \rangle -  \langle p, x_i \rangle \langle x_i, p \rangle + \langle p, x_i \rangle^2 \langle x_i, x_i \rangle  =   1-  \langle p, x_i \rangle^2.
\end{aligned}
\end{align*}
Thus $X(t)$ remains in the compact set $(\mathbb S^{d-1})^N$ throughout its maximal interval of existence. Since the vector field in \eqref{C-1} is smooth on $(\mathbb R^d)^N$, the continuation theorem gives $T_{\max}=\infty$.
\end{proof}
\begin{lemma} \label{L3.2}
Let $x$ be a value in $[-1, 1]$. Then, the following trigonometric identity holds.
\[ \sin \Big(\cos^{-1} x \Big) =  \sqrt{1-x^2}. \]
\end{lemma}
\begin{proof}
We set
\[ \theta\coloneqq\cos^{-1}x\in[0,\pi],\quad X\coloneqq\sin(\cos^{-1}x)=\sin\theta\in[0,1]. \]
Then we have
\[ X^2=\sin^2\theta=1-\cos^2\theta=1-x^2, \quad \mbox{and hence} \quad X=\sqrt{1-x^2}. \]
\end{proof}
\noindent For $\alpha\in(0,\pi)$, let $\alpha^\infty\in\left(0,\frac{\pi}{3}\right]$ be the constant appearing in \eqref{B-2}, and we set 
\begin{equation}\notag 
O(\mathbb{S}^{d-1})_{\alpha}\coloneqq \left\{x\in \mathbb{S}^{d-1}:\langle x,p\rangle >\cos\alpha \right\}.
\end{equation}
In the next proposition, we show that the set $O(\mathbb{S}^{d-1})_{\alpha}^N$ is positively invariant  under the dynamics \eqref{C-1}.
 \begin{proposition} \label{P3.1}
Suppose that system parameters and initial data satisfy 
\[
\alpha \in (0, \pi), \quad \kappa> \frac{\|\mathfrak A p\|_{\infty}}{\sin\alpha(1+\cos\alpha)}, \quad x_i^0 \in O(\mathbb{S}^{d-1})_{\alpha}, \quad  \forall~i \in [N], \]
and let $X= \{ x_i  \}_{i=1}^{N}$ be a global solution to \eqref{C-1}. Then, the flow $X(t)$ is confined to $O(\bbs^{d-1})_\alpha^N$.
\[ X(t)\in O(\mathbb{S}^{d-1})_{\alpha}^N \quad t > 0. \]
\end{proposition}
\begin{proof}
Below, we provide the proof in two steps. First, we derive a differential inequality for $\langle x_i,p\rangle$, and then we use a first-exit argument to prove positive invariance. \newline

\noindent\textbf{$\bullet$ Step A} (Differential inequality for $\langle x_i,p\rangle$): We use \eqref{C-2} and \eqref{C-3} to obtain
\begin{align}\label{C-6}
\begin{aligned}
\frac{d}{dt}\langle x_i,p\rangle &= \Big \langle  A_i x_i+ \kappa R (p-\langle x_i, p \rangle x_i), p \Big \rangle \\
&= \langle {A_i}x_i,p\rangle+\kappa R(1-\langle x_i,p\rangle^2)\\
&= \left\langle {A_i}x_i,p-\langle p,x_i\rangle x_i\right\rangle+\kappa R(1-\langle x_i,p\rangle^2)  \\
&\overset{(*)}{\ge} -\| A_i p \| \sqrt{1-\langle x_i,p\rangle^2}+\kappa R(1-\langle x_i,p\rangle^2) \\
&\geq -\|\mathfrak Ap\|_{\infty} \sqrt{1-\langle x_i,p\rangle^2}+\kappa R(1-\langle x_i,p\rangle^2).
\end{aligned}
\end{align}
Here in $(*)$, we use \eqref{B-4-1} and $\eqref{C-2}$ to see
\[
\Big | \left\langle {A_i}x_i,p-\langle p,x_i\rangle x_i\right\rangle \Big| \leq \| A_i p \| \| x_i \| \cdot \| p-\langle p,x_i\rangle x_i \| \leq  \| A_i p \| \sqrt{1 - \langle x_i, p \rangle^2}.
\]
\noindent\textbf{$\bullet$ Step B} (Positive invariance of $O(\mathbb{S}^{d-1})_{\alpha}^N$): We define a set ${\mathcal T}$ and its supremum $\tau_1^*$ as follows:
\begin{equation}\notag 
{\mathcal T} \coloneqq  \Big \{ \tau>0:\langle x_i(t),p\rangle>  \cos \alpha, \quad t\in [0, \tau),~~\forall~i  \in [N] \Big \}, \quad  \tau_1^* \coloneqq  \sup  {\mathcal T}.    
\end{equation}
By the continuity of the map $t~\mapsto~\langle x_i(t),p\rangle$ and the initial condition, the set ${\mathcal T}$ is nonempty. Hence 
\[ \tau_1^* \in (0, \infty]. \]
We claim that 
\begin{equation} \label{C-7}
 \tau_1^* = \infty. 
 \end{equation}
Suppose the contrary holds, i.e., 
\[ \tau_1^* < \infty \quad \mbox{and} \quad \langle x_i(t),p\rangle> \cos \alpha, \quad t\in [0, \tau_1^*),~~\forall~i  \in [N]. \]
Therefore, there exists $i_0 \in [N]$ such that 
\begin{align}
\begin{aligned} \label{C-8}
&  \langle x_{i_0}(t), p\rangle > \cos \alpha, \quad t \in (0, \tau_1^*), \quad \langle x_{i_0}(\tau_1^*), p\rangle = \cos \alpha, \\
& \langle x_{j}(t), p\rangle \geq \cos \alpha, \quad t \in [0, \tau_1^*],~~\forall~j \in [N] - \{i_0 \}.
\end{aligned}
\end{align}
This yields
\begin{equation} \label{C-9}
R(t)=\frac{1}{N}\sum_{j=1}^N\bigl(1+\langle x_j(t),p\rangle\bigr)
\geq1+\cos\alpha,\quad t\in[0,\tau_1^*].
\end{equation}
It follows from \eqref{C-6}, \eqref{C-8}, \eqref{C-9} and the assumption on $\kappa$ that 
\begin{align*}
\begin{aligned}
\frac{d}{dt}\langle x_{i_0}, p\rangle \Big|_{t = \tau_1^*} &\geq -\|\mathfrak Ap\|_{\infty} \sqrt{1-\langle x_{i_0}(\tau_1^*),p\rangle^2}+\kappa R(\tau_1^*)(1-\langle x_{i_0} (\tau_1^*),p\rangle^2) \\
&\ge  \sin \alpha \Big(  -\|\mathfrak Ap\|_{\infty}  + \kappa (1 + \cos \alpha) \sin \alpha \Big) > 0.
\end{aligned}
\end{align*}
Therefore, $\langle x_{i_0}, p\rangle$ is strictly increasing in a neighborhood of $t = \tau_1^*$. Moreover, since $\langle x_{i_0}(\tau_1^*), p\rangle = \cos \alpha$, for some $\varepsilon \ll 1$, 
\[ 
\langle x_{i_0}(t), p\rangle < \cos \alpha, \quad  t \in (\tau_1^* - \varepsilon, \tau_1^*),
\]
which is contradictory to $\eqref{C-8}_1$. Therefore, we verified the claim \eqref{C-7} and
\[
X(t)\in O(\mathbb{S}^{d-1})_{\alpha}^N, \quad t > 0.
\]
\end{proof}
Suppose that $y:[0,\infty)\to\mathbb R$ is continuous and satisfies
\begin{equation} \label{C-9-1}
\begin{cases}
\displaystyle D^+y(t)\leq\|\mathfrak Ap\|_{\infty}
-\kappa(1+\cos y(t))\sin y(t),\quad t\geq0, \\
\\
y(0)\in[\alpha^\infty,\alpha], \quad \alpha \in (\pi/3, \pi).
\end{cases}
\end{equation}
where we define $\alpha^\infty \in (0, \pi/3)$ as in \eqref{B-2}. \newline

Next, we introduce an auxiliary lemma to prove finite-time entrance into $O(\bbs^{d-1})_{\alpha^\infty}^N$. 
\begin{lemma} \label{L3.3}
Suppose that $\kappa$ satisfies the same condition as in Proposition~\ref{P3.1}:
\[\kappa > \frac{\|\mathfrak Ap\|_{\infty}}{\sin\alpha(1+\cos\alpha)},\] and let $y:[0,\infty)\to\mathbb R$ be a continuous function satisfying \eqref{C-9-1}.
Then, there exists a constant $T>0$ such that 
\begin{equation}\label{C-10}
y(T) < \alpha^\infty \quad \mbox{and} \quad  T<T_1^*\coloneqq
\frac{\pi}{\kappa(1+\cos\alpha)\sin\alpha-\|\mathfrak Ap\|_{\infty}}.
\end{equation}
\end{lemma}
\begin{proof}
Suppose the contrary holds: $y(t)$ is bounded below by $\alpha^\infty$ for all $t\in[0,T_1^*]$. 
\begin{equation}\label{C-11}
    y(t) \ge \alpha^\infty, \quad \forall\, t\in [0, T_1^*].
\end{equation}
First, the condition on $\kappa$ implies
\[D^+ y\Big|_{y = \alpha} \le \|\mathfrak Ap\|_{\infty}-\kappa (1+\cos \alpha) \sin \alpha < 0,\]
and the right-hand side of the differential inequality $\eqref{C-9-1}_1$ remains strictly negative in a neighborhood of $\alpha$ by continuity. Scalar comparison for upper right Dini derivatives therefore prevents $y$ from reaching $\alpha$ from below. Consequently,
\begin{equation}\label{C-12}
    y(t) < \alpha,\quad \forall\, t>0.
\end{equation}
We combine \eqref{C-11} and \eqref{C-12} to obtain
\begin{equation}\label{C-13}
    \alpha^\infty \le y(t) \le \alpha, \quad \forall\, t\in[0, T_1^*].
\end{equation}

Next, since the map $x\mapsto \sin x(1+\cos x)$ is increasing on $(0, \frac{\pi}{3}]$ and decreasing on $[\frac{\pi}{3}, \pi)$ and takes the same value at $\alpha^\infty$ and $\alpha$, we have
\begin{align}
\begin{aligned}  \label{C-14}
& \sin y(1+\cos y) \ge \sin\alpha(1+\cos\alpha)
\quad
\text{whenever }
y\in[\alpha^\infty,\alpha] \\
& \hspace{1cm} \Longrightarrow \|\mathfrak Ap\|_{\infty}-\kappa (1+\cos y) \sin y \leq \|\mathfrak Ap\|_{\infty}-\kappa (1+\cos \alpha) \sin \alpha < 0.
\end{aligned}
\end{align}
Let $z = z(t)$ be a solution to the following Cauchy problem:
\begin{equation}\notag
\begin{cases} 
\displaystyle {\dot z} =  \|\mathfrak Ap\|_{\infty}-\kappa (1+\cos \alpha) \sin \alpha, \quad t > 0,\\
\displaystyle z \Big|_{t = 0+} = \alpha.
\end{cases}
\end{equation}
From \eqref{C-13} and \eqref{C-14}, we have
\[y(0) \le z(0),\quad D^+ y(t)\le \dot{z}(t), \quad \forall\, t\in(0, T_1^*],\]
so  we have
\[ D^+(y-z)=D^+y-\dot z\leq0. \]
The scalar comparison principle for upper right Dini derivatives shows that $y-z$ is non-increasing on $(0,T_1^*)$; continuity at the endpoints therefore gives
\begin{equation}\notag 
y(t) \leq z(t) =  \alpha + \Big(   \|\mathfrak Ap\|_{\infty}-\kappa (1+\cos \alpha) \sin \alpha \Big)t. 
\end{equation}
In particular, it follows from \eqref{C-10} that 
\[y(T_1^*) \le z(T_1^*) = \alpha + \Big(   \|\mathfrak Ap\|_{\infty}-\kappa (1+\cos \alpha) \sin \alpha \Big) T_1^* = \alpha - \pi < \alpha^\infty. \]
This contradicts \eqref{C-11}. Hence there exists $T<T_1^*$ such that
\[y(T) < \alpha^\infty.\]
\end{proof} 
Combining Proposition~\ref{P3.1} with Lemma~\ref{L3.3}, we show the finite-time entrance of the flow into $O(\mathbb{S}^{d-1})_{\alpha^\infty}^N$. This gives the following corollary.
\begin{corollary} \label{C3.1}
Suppose that system parameters and initial data satisfy 
\[
\alpha \in (0, \pi), \quad \kappa> \frac{\|\mathfrak Ap\|_{\infty}}{\sin\alpha(1+\cos\alpha)}, \quad x_i^0 \in O(\mathbb{S}^{d-1})_{\alpha}, \quad  \forall~i \in [N], \]
and let $X= \{ x_i  \}_{i=1}^{N}$ be a global solution to \eqref{C-1}. Then, there exists a time $T\le T_1^*$ such that 
\[ X(t)\in O(\mathbb{S}^{d-1})_{\alpha^\infty}^N \quad \mbox{for $t\ge T$}, \]
where $T_1^*$ is defined in \eqref{C-10}.
\end{corollary}
\begin{proof} 
We claim that there exists a finite time $T$ such that 
\begin{equation} \label{C-15}
X(T)\in O(\mathbb{S}^{d-1})_{\alpha^\infty}^N.
\end{equation}
Then, since the set $O(\mathbb{S}^{d-1})_{\alpha^\infty}^N$ is positively invariant by Proposition~\ref{P3.1}, we have the desired result:
\[  X(t)\in O(\mathbb{S}^{d-1})_{\alpha^\infty}^N \quad \mbox{for $t\ge T$}. \]
\noindent {\it Proof of \eqref{C-15}}: For $\alpha \in (0, \pi/3]$, we have $\alpha^{\infty} = \alpha$, and we can choose 
\[ T= 0, \quad  X^0 \in O(\mathbb{S}^{d-1})_{\alpha^\infty}^N. \]
Thus, we may assume that $\alpha \in (\pi/3, \pi).$ We use the computations in \eqref{C-6} to see the differential inequality for $\langle x_i, p\rangle$:
\[\frac{d}{dt}\langle x_i,p\rangle\geq -\|\mathfrak Ap\|_{\infty} \sqrt{1-\langle x_i,p\rangle^2}+\kappa R(1-\langle x_i,p\rangle^2).
\]
From this, we derive the differential inequality for $\cos^{-1}\langle x_i,p\rangle$. We split into two cases:
\[\text{Case 1: } \langle x_i, p \rangle \ne 1 \quad \text{or}\quad \text{Case 2: } \langle x_i, p\rangle = 1.\]

\vspace{0.1cm}

\noindent\textbf{$\bullet$ Case 1} ($\langle x_i, p\rangle \ne 1$): In this case, the inequality \eqref{C-6} implies the following relation for $t\in [0,T_1^*)$:
\begin{align}
\begin{aligned} \label{C-16}
\frac{d}{dt}\cos^{-1}\langle x_i,p\rangle &= -\frac{1}{\sqrt{1 - \langle x_i, p \rangle^2}} \frac{d}{dt} \langle x_i, p \rangle \leq \| A_i p \| - \kappa R \sqrt{1-\langle x_i,p\rangle^2} \\
&\le \|\mathfrak Ap\|_{\infty} - \frac{\kappa}{N} \sum_{j=1}^N \Big (1+\langle x_j,p\rangle \Big) \sin\left(\cos^{-1}\langle x_i,p\rangle\right), \\
&\le \|\mathfrak Ap\|_{\infty} - \frac{\kappa}{N} \sum_{j=1}^N \Big (1+ \cos (\cos^{-1} \langle x_j,p\rangle) \Big) \sin\left(\cos^{-1}\langle x_i,p\rangle\right),
\end{aligned}
\end{align}
where we use \eqref{C-4} to write
\[
R=\frac{1}{N}\sum_{j=1}^N\Big(1+\langle x_j,p\rangle\Big)=\frac{1}{N}\sum_{j=1}^N\Big(1+\cos\big(\cos^{-1}\langle x_j,p\rangle\big)\Big).
\]

\vspace{0.5em}
\noindent\textbf{$\bullet$ Case 2} ($\langle x_i, p \rangle = 1$): At such a time, $x_i=p$ and the coupling term vanishes, so $\dot{x}_i=A_ip$ in general. The local expansion of a differentiable curve on the sphere at $p$ gives
\[
D^+\cos^{-1}\langle x_i,p\rangle
= \|A_ip\|
\le \|\mathfrak Ap\|_{\infty}.
\]
Since $\sin(\cos^{-1}\langle x_i,p\rangle)=0$, we have
\begin{equation}\label{C-17}
\hspace{1cm} D^+\cos^{-1}\langle x_i,p\rangle
\leq \|\mathfrak Ap\|_{\infty} - \frac{\kappa}{N} \sum_{j=1}^N \Big (1+ \cos (\cos^{-1} \langle x_j,p\rangle) \Big) \sin\left(\cos^{-1}\langle x_i,p\rangle\right).
\end{equation}
We set 
\[ \vartheta_i(t)\coloneqq\cos^{-1}\langle x_i(t),p\rangle. \]
From \eqref{C-16} and \eqref{C-17}, we obtain, for every $i\in[N]$,
\begin{equation}\label{C-18}
D^+\vartheta_i(t)
\leq \|\mathfrak Ap\|_{\infty} - \frac{\kappa}{N} \sum_{j=1}^N \big(1+\cos\vartheta_j(t)\big)\sin\vartheta_i(t).
\end{equation}
Now let $\mathcal{I}(t)$ be the set of maximal indices:
\[\mathcal{I}(t) \coloneqq \big\{j\in [N]: \vartheta_j(t) = \max_{i\in[N]} \vartheta_i(t) \big\}.\]
For each $i_0\in\mathcal{I}(t)$, \eqref{C-18} gives
\begin{equation}\notag
D^+\vartheta_{i_0}(t)
\leq \|\mathfrak Ap\|_{\infty}
-\kappa\left(1+\cos\Big(\max_{i\in[N]}\vartheta_i(t)\Big)\right)
\sin\Big(\max_{i\in[N]}\vartheta_i(t)\Big).
\end{equation}
Using the upper-Dini derivative rule for a finite maximum, we obtain
\begin{align*}
\begin{aligned}
D^+\max_{i \in [N]}\vartheta_i(t)
&\leq \max_{i_0\in \mathcal{I}(t)}D^+\vartheta_{i_0}(t) \\
&\leq \|\mathfrak Ap\|_{\infty}-\kappa \left(1+\cos \Big( \max_{i \in [N]}\vartheta_i(t) \Big)\right) \sin \Big( \max_{i \in [N]}\vartheta_i(t) \Big),\quad t\in [0, T_1^*).
\end{aligned}
\end{align*}
If $\max_{i \in [N]}\vartheta_i(0) < \alpha^{\infty}$, then 
\[
X^0\in O(\mathbb S^{d-1})_{\alpha^\infty}^N.
\]
By Proposition~\ref{P3.1}, this set is positively invariant.
Hence we may take $T=0$.
Otherwise, Lemma~\ref{L3.3} yields a time $T<T_1^*$ such that
\[
\max_{i \in [N]}\vartheta_i(T) < \alpha^{\infty}, \quad \mbox{i.e.,} \quad \cos^{-1} \langle x_i(T),p\rangle < \alpha^{\infty}, \quad \forall~i \in [N].
\]
This implies
\[
\langle x_i(T),p\rangle > \cos \alpha^{\infty}, \quad \forall~i \in [N],
\]
or equivalently, 
\[ X(T)\in O(\mathbb{S}^{d-1})_{\alpha^\infty}^N. \]
This verifies the claim \eqref{C-15}.
\end{proof}
\subsection{Exponential stability}\label{subsec:3-2}
In this subsection, we study the exponential stability for two solutions in $O(\mathbb{S}^{d-1})_{\alpha^\infty}^N$. Throughout this subsection, the parameters $\alpha\in(\frac{\pi}{3},\pi)$ and $\alpha^\infty\in(0,\frac{\pi}{3})$ are related by \eqref{B-2}, so the coefficient $\mu_{\alpha^\infty}$ defined by
\[\mu_{\alpha^\infty}\coloneqq(2\cos\alpha^\infty-1)(1+\cos\alpha^\infty)\] is positive. 

\vspace{0.1cm}

Recall that the mixed norm $\|\cdot\|_{q,2}$, $1\le q\le\infty$, was defined in \eqref{A-2}. Let $X=(x_1,\cdots,x_N)$ and $\widetilde X=(\widetilde x_1,\cdots,\widetilde x_N)$ be two global solutions to \eqref{C-1} with the same tuple of natural frequency matrices $\mathfrak A$. Before rewriting their equations, we define the corresponding order parameters by
\begin{equation} \label{C-19}
R(X)\coloneqq \frac{1}{N}\sum_{j=1}^N(1+\langle x_j,p\rangle),\quad R(\widetilde X)\coloneqq \frac{1}{N}\sum_{j=1}^N(1+\langle\widetilde x_j,p\rangle).
\end{equation}
Then, $X$ and $\widetilde X$ satisfy
\begin{align}
\begin{aligned} \label{C-20}
\dot{x}_i&=A_i x_i+\kappa R(X)(p-\langle x_i,p\rangle x_i),\quad i\in[N].\\
\dot{\widetilde x}_i&=A_i\widetilde x_i+\kappa R(\widetilde X)(p-\langle\widetilde x_i,p\rangle\widetilde x_i).
\end{aligned}
\end{align}
We subtract the two equations in \eqref{C-20} to obtain
\begin{align}
\begin{aligned} \label{C-21}
\frac{d}{dt}(x_i-\widetilde x_i)&=A_i(x_i-\widetilde x_i)+\kappa R(X)(\langle\widetilde x_i,p\rangle\widetilde x_i-\langle x_i,p\rangle x_i)\\
&\quad+\kappa(R(X)-R(\widetilde X))(p-\langle\widetilde x_i,p\rangle\widetilde x_i).
\end{aligned}
\end{align}
Taking the inner product of \eqref{C-21} with $2(x_i-\widetilde x_i)$ gives
\begin{align}
\begin{aligned} \label{C-22}
\frac{d}{dt}\|x_i-\widetilde x_i\|^2&=2\langle A_i(x_i-\widetilde x_i),x_i-\widetilde x_i\rangle\\
&\quad+2\kappa R(X)\langle\langle\widetilde x_i,p\rangle\widetilde x_i-\langle x_i,p\rangle x_i,x_i-\widetilde x_i\rangle\\
&\quad+2\kappa(R(X)-R(\widetilde X))\langle p-\langle\widetilde x_i,p\rangle\widetilde x_i,x_i-\widetilde x_i\rangle\\
&=: \mathcal I_{11}+\mathcal I_{12}+\mathcal I_{13}.
\end{aligned}
\end{align}
In the following lemma, we estimate three terms ${\mathcal I}_{1i},~i=1,2,3$ one by one. 
\begin{lemma} \label{L3.4}
Suppose that $X,\widetilde X\in O(\mathbb{S}^{d-1})_{\alpha^\infty}^N$. Then, for every $i\in[N]$,
\begin{align*}
\begin{aligned}
& (i)~\mathcal I_{11}=0,\quad \mathcal I_{12}\le-2\kappa\cos\alpha^\infty(1+\cos\alpha^\infty)\|x_i-\widetilde x_i\|^2, \\
& (ii)~|\mathcal I_{13}|\le2\kappa\sin^2\alpha^\infty\|X-\widetilde X\|_{1,2}\|x_i-\widetilde x_i\|.
\end{aligned}
\end{align*}
\end{lemma}
\begin{proof}
\noindent (i) Since $A_i^T=-A_i$, we have
\[
\mathcal I_{11}=2\langle A_i(x_i-\widetilde x_i),x_i-\widetilde x_i\rangle=0.
\]
Since $\|x_i\|=\|\widetilde x_i\|=1$, direct calculation yields
\begin{align*}
\begin{aligned}
&\langle\langle\widetilde x_i,p\rangle\widetilde x_i-\langle x_i,p\rangle x_i,x_i-\widetilde x_i\rangle =\langle\widetilde x_i,p\rangle(\langle\widetilde x_i,x_i\rangle-1)-\langle x_i,p\rangle(1-\langle x_i,\widetilde x_i\rangle)\\
&\quad=-(\langle x_i,p\rangle+\langle\widetilde x_i,p\rangle)(1-\langle x_i,\widetilde x_i\rangle)=-\frac{1}{2}(\langle x_i,p\rangle+\langle\widetilde x_i,p\rangle)\|x_i-\widetilde x_i\|^2.
\end{aligned}
\end{align*}
Hence, we have
\[
\mathcal I_{12}=-\kappa R(X)(\langle x_i,p\rangle+\langle\widetilde x_i,p\rangle)\|x_i-\widetilde x_i\|^2.
\]
The conditions $X,\widetilde X\in O(\mathbb{S}^{d-1})_{\alpha^\infty}^N$ imply
\[
R(X)\ge1+\cos\alpha^\infty,\quad \langle x_i,p\rangle+\langle\widetilde x_i,p\rangle\ge2\cos\alpha^\infty,
\]
and therefore we have
\[
\mathcal I_{12}\le-2\kappa\cos\alpha^\infty(1+\cos\alpha^\infty)\|x_i-\widetilde x_i\|^2.
\]
\noindent (ii) By Lemma~\ref{L3.1} and \eqref{C-2},
\begin{align} \label{C-23}
\begin{aligned}
|\mathcal I_{13}|&\le2\kappa|R(X)-R(\widetilde X)|\|p-\langle\widetilde x_i,p\rangle\widetilde x_i\|\|x_i-\widetilde x_i\|\\
&= 2\kappa|R(X)-R(\widetilde X)|\sqrt{1-\langle \widetilde x_i, p\rangle^2}\|x_i-\widetilde x_i\| \\
&\le2\kappa\sin\alpha^\infty|R(X)-R(\widetilde X)|\|x_i-\widetilde x_i\|.
\end{aligned}
\end{align}
We next estimate the difference of the order parameters. For $x_j,\widetilde x_j\in O(\mathbb{S}^{d-1})_{\alpha^\infty}$, we use the orthogonal decomposition
\[x_j = \langle x_j, p \rangle p + (x_j - \langle x_j, p \rangle p)\]
and the Cauchy--Schwarz inequality to obtain
\begin{align*}
\begin{aligned}
\langle x_j,\widetilde x_j\rangle
&= \langle x_j,p\rangle\langle\widetilde x_j,p\rangle
 +\langle x_j-\langle x_j,p\rangle p,\widetilde x_j-\langle\widetilde x_j,p\rangle p\rangle\\
&\leq \langle x_j,p\rangle\langle\widetilde x_j,p\rangle
 +\sqrt{1-\langle x_j,p\rangle^2}\sqrt{1-\langle\widetilde x_j,p\rangle^2}\\
&=\cos\left(\cos^{-1}\langle x_j,p\rangle-\cos^{-1}\langle\widetilde x_j,p\rangle\right).
\end{aligned}
\end{align*}
Consequently, we have
\[
\|x_j-\widetilde x_j\|^2 =2(1-\langle x_j,\widetilde x_j\rangle) \ge4\sin^2\left(\frac{\cos^{-1}\langle x_j,p\rangle-\cos^{-1}\langle\widetilde x_j,p\rangle}{2}\right).
\]
Since $0\le\cos^{-1}\langle x_j,p\rangle,\cos^{-1}\langle\widetilde x_j,p\rangle<\alpha^\infty<\frac{\pi}{2}$, we also have
\begin{align*}
|\langle x_j-\widetilde x_j,p\rangle|&=\left|\cos\left(\cos^{-1}\langle x_j,p\rangle\right)-\cos\left(\cos^{-1}\langle\widetilde x_j,p\rangle\right)\right|\\
&=2\sin\left(\frac{\cos^{-1}\langle x_j,p\rangle+\cos^{-1}\langle\widetilde x_j,p\rangle}{2}\right)\left|\sin\left(\frac{\cos^{-1}\langle x_j,p\rangle-\cos^{-1}\langle\widetilde x_j,p\rangle}{2}\right)\right|\\
&\le\sin\alpha^\infty\|x_j-\widetilde x_j\|.
\end{align*}
Using \eqref{C-19} and the definition of $\|\cdot\|_{1,2}$, we obtain
\begin{equation} \label{C-24}
|R(X)-R(\widetilde X)|\le\frac{1}{N}\sum_{j=1}^N|\langle x_j-\widetilde x_j,p\rangle|\le\sin\alpha^\infty\|X-\widetilde X\|_{1,2}.
\end{equation}
Finally, we combine \eqref{C-23} and \eqref{C-24} to verify the estimate for $\mathcal I_{13}$.
\end{proof}
Next, we study the exponential stability of the Cauchy problem \eqref{C-1} with respect to initial data. 
\begin{proposition} \label{P3.2}
Suppose that system parameters and initial data satisfy
\[
\alpha\in\left(\frac{\pi}{3},\pi\right),\quad \kappa>\frac{\|\mathfrak Ap\|_{\infty}}{\sin\alpha(1+\cos\alpha)},\quad X^0,\widetilde X^0\in O(\mathbb{S}^{d-1})_{\alpha^\infty}^N, \quad 1\le q\le\infty, 
\]
and let $X$ and $\widetilde X$ be the global solutions to \eqref{C-1} corresponding to initial data $X^0$ and $\widetilde{X}^0$, respectively. Then, we have exponential stability: \begin{equation} \label{C-25}
\|X(t)-\widetilde X(t)\|_{q,2}\le e^{-\kappa\mu_{\alpha^\infty}t}\|X^0-\widetilde X^0\|_{q,2},\quad t\ge0.
\end{equation}
\end{proposition}
\begin{proof}
It follows from \eqref{B-2} that 
\[
\sin\alpha(1+\cos\alpha)=\sin\alpha^\infty(1+\cos\alpha^\infty).
\]
Thus, Proposition \ref{P3.1}, applied with $\alpha^\infty$ in place of $\alpha$, yields
\begin{equation} \label{C-26}
X(t),\widetilde X(t)\in O(\mathbb{S}^{d-1})_{\alpha^\infty}^N,
\quad t\ge0.
\end{equation}
By \eqref{C-22}, Lemma \ref{L3.4}, and \eqref{C-26},
\[
\frac{d}{dt}\|x_i-\widetilde x_i\|^2\le-2\kappa\cos\alpha^\infty(1+\cos\alpha^\infty)\|x_i-\widetilde x_i\|^2 +2\kappa\sin^2\alpha^\infty\|X-\widetilde X\|_{1,2}\|x_i-\widetilde x_i\|.
\]
\newline
\noindent $\bullet$~{\bf Case 1}:~If $x_i\ne\widetilde x_i$, we divide the above relation by $2\|x_i-\widetilde x_i\|$ to find 
\[
\frac{d}{dt}\|x_i-\widetilde x_i\|\le-\kappa\cos\alpha^\infty(1+\cos\alpha^\infty)\|x_i-\widetilde x_i\|+\kappa\sin^2\alpha^\infty\|X-\widetilde X\|_{1,2}.
\]
\newline
\noindent $\bullet$~{\bf Case 2}:~If $x_i=\widetilde x_i$, we use \eqref{C-21}, Lemma \ref{L3.1} and \eqref{C-24} to derive the same estimate for the upper Dini derivative. Hence, for every $i\in[N]$,
\begin{equation} \label{C-27}
D^+\|x_i-\widetilde x_i\|\le-\kappa\cos\alpha^\infty(1+\cos\alpha^\infty)\|x_i-\widetilde x_i\|+\kappa\sin^2\alpha^\infty\|X-\widetilde X\|_{1,2}.
\end{equation}
Now, we divide the analysis into three cases:
\[\text{Case A: }q=1, \quad \text{Case B: }1<q<\infty, \quad \text{Case C: }q=\infty.\]

\noindent\textbf{$\diamond$ Case 2.1} ($q=1$): We average \eqref{C-27} over $i\in[N]$ to find
\begin{align*}
D^+\|X-\widetilde X\|_{1,2}&\le-\kappa\cos\alpha^\infty(1+\cos\alpha^\infty)\|X-\widetilde X\|_{1,2}+\kappa\sin^2\alpha^\infty\|X-\widetilde X\|_{1,2}\\
&=-\kappa\mu_{\alpha^\infty}\|X-\widetilde X\|_{1,2}.
\end{align*}

\noindent\textbf{$\diamond$ Case 2.2} ($1<q<\infty$): If $\|X-\widetilde X\|_{q,2}=0$, uniqueness of solutions to \eqref{C-1} implies that $X$ and $\widetilde X$ coincide from that time onward. Otherwise, the chain rule and \eqref{C-27} imply
\begin{align*}
\begin{aligned}
D^+\|X-\widetilde X\|_{q,2}&\le\|X-\widetilde X\|_{q,2}^{1-q}\frac{1}{N}\sum_{i=1}^N\|x_i-\widetilde x_i\|^{q-1}D^+\|x_i-\widetilde x_i\|\\
&\le-\kappa\cos\alpha^\infty(1+\cos\alpha^\infty)\|X-\widetilde X\|_{q,2}\\
&~+\kappa\sin^2\alpha^\infty\|X-\widetilde X\|_{q,2}^{1-q}\left(\frac{1}{N}\sum_{i=1}^N\|x_i-\widetilde x_i\|^{q-1}\right)\|X-\widetilde X\|_{1,2}.
\end{aligned}
\end{align*}
By H\"older's inequality, we have
\[
\frac{1}{N}\sum_{i=1}^N\|x_i-\widetilde x_i\|^{q-1}\le\|X-\widetilde X\|_{q,2}^{q-1},\quad \|X-\widetilde X\|_{1,2}\le\|X-\widetilde X\|_{q,2}.
\]
Therefore, we have
\[
D^+\|X-\widetilde X\|_{q,2}\le-\kappa\mu_{\alpha^\infty}\|X-\widetilde X\|_{q,2}.
\]
\newline
\noindent\textbf{$\diamond$ Case 2.3} ($q=\infty$):~We use the upper right Dini derivative rule for a finite maximum and \eqref{C-27} to get 
\begin{align*}
D^+\|X-\widetilde X\|_{\infty,2}&\le\max_{\substack{i\in[N]\\\|x_i-\widetilde x_i\|=\|X-\widetilde X\|_{\infty,2}}}D^+\|x_i-\widetilde x_i\|\\
&\le-\kappa\cos\alpha^\infty(1+\cos\alpha^\infty)\|X-\widetilde X\|_{\infty,2}+\kappa\sin^2\alpha^\infty\|X-\widetilde X\|_{1,2}\\
&\le-\kappa\mu_{\alpha^\infty}\|X-\widetilde X\|_{\infty,2}.
\end{align*}
Now, we combine Case A, Case B and Case C to see that for every $1\le q\le\infty$,
\begin{equation} \label{C-28}
D^+\|X-\widetilde X\|_{q,2}\le-\kappa\mu_{\alpha^\infty}\|X-\widetilde X\|_{q,2}.
\end{equation}
Finally, we apply Gronwall's lemma to \eqref{C-28} to get \eqref{C-25}.
\end{proof}

\subsection{Relaxation dynamics} \label{subsec:3-3}
In this subsection, we use the exponential-stability estimate in Proposition~\ref{P3.2} and compare a solution with its time translates to prove the existence and uniqueness of an equilibrium in $O(\mathbb{S}^{d-1})_{\alpha^\infty}^N$ and exponential relaxation toward it.

\begin{proposition} \label{P3.3}
Suppose system parameters satisfy 
\[
\alpha\in\left(\frac{\pi}{3},\pi\right),\quad \kappa>\frac{\|\mathfrak Ap\|_{\infty}}{\sin\alpha(1+\cos\alpha)}.
\]
Then, system \eqref{C-1} admits a unique equilibrium $X^\infty\in O(\mathbb{S}^{d-1})_{\alpha^\infty}^N$.
\end{proposition}
\begin{proof}
For $X^0\in O(\mathbb{S}^{d-1})_{\alpha^\infty}^N$, let $X=X(t)$ be the corresponding global solution to \eqref{C-1}. By \eqref{B-2}, the coupling hypothesis also holds with $\alpha^\infty$ in place of $\alpha$. Hence Proposition~\ref{P3.1} yields
\begin{equation} \label{C-29}
X(t)\in O(\mathbb{S}^{d-1})_{\alpha^\infty}^N,\quad t\geq0.
\end{equation}
For $s\geq0$, we define 
\[ \widetilde X(t)\coloneqq X(t+s). \] 
Since \eqref{C-1} is autonomous, $\widetilde X$ is another solution in the same stable region. Proposition~\ref{P3.2} with $q=2$ gives
\begin{equation}\label{C-30}
\|X(t+s)-X(t)\|_{2,2}
\leq e^{-\kappa\mu_{\alpha^\infty}t}\|X(s)-X(0)\|_{2,2},
\quad t,s\geq0.
\end{equation}
Because $(\mathbb{S}^{d-1})^N$ is compact, there exists $C>0$ such that
\[
\|X(s)-X(0)\|_{2,2}\leq C,\quad s\geq0.
\]
For $t_2\geq t_1\geq0$, we take $t=t_1$ and $s=t_2-t_1$ in \eqref{C-30} to see
\[
\|X(t_2)-X(t_1)\|_{2,2} \leq e^{-\kappa\mu_{\alpha^\infty}t_1}\|X(t_2-t_1) - X(0)\|_{2,2}
\leq Ce^{-\kappa\mu_{\alpha^\infty}t_1}.
\]
Thus $X(t)$ is Cauchy as $t\to\infty$. Since $(\mathbb{S}^{d-1})^N$ is complete, there exists $X^\infty\in(\mathbb{S}^{d-1})^N$ such that
\begin{equation}\label{C-31}
\lim_{t\to\infty}\|X(t)-X^\infty\|_{2,2}=0.
\end{equation}
By \eqref{C-29}, the limit satisfies $X^\infty\in\overline{O(\mathbb{S}^{d-1})_{\alpha^\infty}}^N$. Let $\Phi_s$ denote the flow of \eqref{C-1}. For every $s\geq0$, continuity of the flow and \eqref{C-31} give
\[
\Phi_s(X^\infty)
=\lim_{t\to\infty}\Phi_s(X(t))
=\lim_{t\to\infty}X(t+s)
=X^\infty,
\]
so $X^\infty$ is an equilibrium. If some component $x_{i_0}^\infty$ lies on the boundary, then 
\[ R(X^\infty)\geq1+\cos\alpha^\infty \]
and the relation \eqref{C-6} yields
\[
0=\frac{d}{dt}\langle x_{i_0}^\infty,p\rangle\geq \sin\alpha^\infty\Big(-\|\mathfrak Ap\|_{\infty}+\kappa(1+\cos\alpha^\infty)\sin\alpha^\infty\Big)>0,
\]
where the last inequality follows from \eqref{B-2} and the hypothesis on $\kappa$. This contradiction shows that 
\[ X^\infty\in O(\mathbb{S}^{d-1})_{\alpha^\infty}^N. \]
To prove uniqueness, let $X^\infty$ and $\widetilde X^\infty$ be two equilibria in $O(\mathbb{S}^{d-1})_{\alpha^\infty}^N$. Applying Proposition~\ref{P3.2} to the corresponding stationary solutions with $q=\infty$, we obtain
\[
\|X^\infty-\widetilde X^\infty\|_{\infty,2}
\leq e^{-\kappa\mu_{\alpha^\infty}t}
\|X^\infty-\widetilde X^\infty\|_{\infty,2},
\quad t\geq0.
\]
Letting $t\to\infty$ yields $X^\infty=\widetilde X^\infty$.
\end{proof}
Next, we combine Corollary~\ref{C3.1}, Proposition~\ref{P3.2}, and Proposition~\ref{P3.3} to obtain the exponential relaxation result. 
\begin{theorem} \label{T3.1}
Suppose that system parameters and initial data satisfy 
\[
\alpha\in\left(\frac{\pi}{3},\pi\right),\quad \kappa>\frac{\|\mathfrak Ap\|_{\infty}}{\sin\alpha(1+\cos\alpha)},\quad X^0\in O(\mathbb{S}^{d-1})_\alpha^N,
\]
and let $X=X(t)$ and $X^\infty\in O(\mathbb{S}^{d-1})_{\alpha^\infty}^N$  be the global solution to \eqref{C-1} with initial data $X^0$, and the unique equilibrium given by Proposition \ref{P3.3}, respectively. Then, there exists $T\le T_1^*$ such that, for every $1\le q\le\infty$,
\begin{equation} \label{C-32}
\|X(t)-X^\infty\|_{q,2}\le e^{-\kappa\mu_{\alpha^\infty}(t-T)}\|X(T)-X^\infty\|_{q,2},\quad t\ge T,
\end{equation}
where $T_1^*$ is defined in \eqref{C-10}.
\end{theorem}
\begin{proof}
By Corollary \ref{C3.1}, there exists $T\le T_1^*$ such that
\[
X(t)\in O(\mathbb{S}^{d-1})_{\alpha^\infty}^N,\quad t\ge T.
\]
The time-translated curve $t\mapsto X(t+T)$ and the constant equilibrium curve $t\mapsto X^\infty$ are two solutions of \eqref{C-1} whose initial data lie in $O(\mathbb{S}^{d-1})_{\alpha^\infty}^N$. Proposition \ref{P3.2} therefore gives
\[
\|X(t+T)-X^\infty\|_{q,2}\le e^{-\kappa\mu_{\alpha^\infty}t}\|X(T)-X^\infty\|_{q,2},\quad t\ge0.
\]
Replacing $t$ by $t-T$ verifies \eqref{C-32}. If $X^0\in O(\mathbb{S}^{d-1})_{\alpha^\infty}^N$, Proposition \ref{P3.1} shows that the trajectory remains in this set for all $t\ge0$, so $T=0$ is admissible.
\end{proof}

\begin{remark}\label{R3.1}
We make the following remarks on Theorem \ref{T3.1}:
\begin{enumerate}
\item
If $X^0\in O(\mathbb{S}^{d-1})_{\alpha^\infty}^N$, then one may take $T=0$.
\vspace{0.1cm}
\item 
As mentioned in the introduction, Theorem \ref{T3.1} shows that every initial configuration avoiding the antipodal point $-p$ converges exponentially to an equilibrium, provided the coupling strength is sufficiently large, with a threshold depending on the initial configuration and natural frequency matrices.
\vspace{0.1cm}
\item
One could also consider the sphere model \eqref{B-3} with a general influence function $I:\mathbb{S}^{d-1}\to [0,\infty)$ in place of $1+\langle x_j,p\rangle$. Such a model is locally compatible with the abstract framework of 
Section~\ref{sec:7} (Theorem~\ref{T6.1}), after choosing a coordinate chart around $p$. To apply Theorem~\ref{T6.1}, however, the coordinate representatives of the vector field
\[
S(x)=p-\langle x,p\rangle x
\]
and of $I$ must satisfy the local assumptions $({\mathcal F}_{B}1)$--$({\mathcal F}_{B}4)$, in particular the strong dissipativity condition \eqref{F-4}. 
\end{enumerate}
\end{remark}
\section{Winfree oscillators on unitary groups I}\label{sec:4}
\setcounter{equation}{0}
In this section, we propose the canonical Winfree model on $\mathrm{U}(d)$, and study finite-time entrance into a stable spectral region, exponential stability there, the existence of a unique equilibrium, and exponential relaxation, following the presentation in the previous section.  We first begin with the reformulation of the Winfree model on ${\mathrm U}(1)={\mathbb S}^1$. For $z_i=e^{{\mathrm i}\theta_i}$, we have
\begin{equation} \label{D-1}
{\dot \theta}_i  =  - {\mathrm i} {\dot z}_i \overline{z_i}, \quad  1 + \cos \theta_i  = 1 + \mbox{Re} z_i, \quad  \sin \theta_i = \frac{z_i - \overline{z_i}}{2 {\mathrm i}}, 
\end{equation}
where $\overline{z}$ is the complex conjugate of  $z$.  Then, we have
\begin{align*}
\begin{aligned}
&  \dot{\theta}_i=\nu_i- \kappa \Big( \frac{1}{N} \sum_{j=1}^N  (1+\cos\theta_j) \Big) \sin\theta_i, \\
& \hspace{0.5cm} \Longleftrightarrow \quad - {\mathrm i} {\dot z}_i \overline{z_i} = \nu_i  + \frac{\kappa {\mathrm i}}{N}  \sum_{j=1}^{N} \frac{(1 + \mbox{Re} z_j)}{2} (z_i - \overline{z_i}) \quad \mbox{by \eqref{D-1}} \\
&  \hspace{0.5cm} \Longleftrightarrow \quad {\dot z}_i \overline{z_i} = {\mathrm i} \nu_i + \frac{\kappa}{N}  \sum_{j=1}^{N} \frac{(1 + \mbox{Re} z_j)}{2} (\overline{z_i} - z_i) \quad \mbox{by multiplying by ${\mathrm i}$},
\end{aligned}
\end{align*}
i.e., the complex counterpart of the Winfree model becomes 
\begin{equation}\notag 
{\dot z}_i \overline{z_i} = {\mathrm i} \nu_i + \frac{\kappa}{N}  \sum_{j=1}^{N} \frac{(1 + \mbox{Re} z_j)}{2} (\overline{z_i} - z_i), \quad i \in [N].
\end{equation}
For the unitary group $\mathrm{U}(d)$ with $d \geq 2$, we use the following correspondence:
\[
(z_i, {\bar z}_i,   \mathrm{i} \nu_i, \mathrm{Re}z_j) \quad \Longleftrightarrow \quad  \Big (U_i, U_i^{\dagger}, \Omega_i, \frac{\operatorname{Re}\operatorname{tr}{U_j}}{d} \Big ),
\]
where $U^{\dagger}$ is the Hermitian adjoint of $U$. Under the above correspondence, the  ${\mathrm U}(d)$-counterpart of \eqref{A-1} reads as follows.
\begin{equation}\label{D-2}
\dot{U}_i{U_i}^{\dagger}={\Omega_i}+\frac{\kappa}{N} \sum_{j=1}^N \frac{d+\operatorname{Re}\operatorname{tr}{U_j}(t)}{2d}\left({U_i}^{\dagger}(t)-{U_i}(t)\right),\quad i \in [N],
\end{equation}
where $\Omega_1,\cdots, \Omega_N\in \mathbb{C}^{d\times d}$ are skew-Hermitian matrices:
\[ \Omega_i^{\dagger} = -\Omega_i, \quad i \in [N]. \]
In what follows, we study the relaxation dynamics of the Winfree matrix model on unitary groups, which can be regarded as the ``quantum'' counterpart of \eqref{A-1}:
\begin{equation}\label{D-3}
\begin{cases}
\displaystyle \dot{U}_i{U_i}^{\dagger}={\Omega_i}+\frac{\kappa}{2N} \sum_{j=1}^N \frac{d+\operatorname{Re}\operatorname{tr}{U_j}(t)}{d}\left({U_i}^{\dagger}(t)-{U_i}(t)\right),\quad t>0,\\
\displaystyle U_i \Big|_{t = 0} = U_i^0\in\mathrm{U}(d), \quad i \in [N],
\end{cases}
\end{equation}
where $\Omega_1,\ldots,\Omega_N\in\mathbb{C}^{d\times d}$ are skew-Hermitian matrices, and $\dagger$ denotes the Hermitian adjoint. We write
\[
\mathfrak V\coloneqq(\Omega_1,\ldots,\Omega_N),\quad
\mathcal U\coloneqq(U_1,\ldots,U_N),
\]
and equip the frequency tuple with the maximum operator norm
\[
\|\mathfrak V\|_\infty\coloneqq\max_{i\in[N]}\|\Omega_i\|_{\mathrm{op}}.
\]
Here $\|\cdot\|_{\mathrm{op}}$ denotes the operator norm on $\mathbb{C}^{d\times d}$ induced by the Euclidean norm on $\mathbb{C}^d$. We will use the Schatten--von Neumann norm to measure distances on $\mathrm{U}(d)$. For a complex $d\times d$ matrix $U$, let $s_1\ge\cdots\ge s_d\ge0$ be its singular values, equivalently the nonnegative square roots of the eigenvalues of $UU^\dagger$. Then for $p\in [1,\infty]$, the Schatten--von Neumann $p$-norm is defined as
\[
\|U\|_{S_p^d}=
\begin{cases}
\left(\sum_{k=1}^d s_{k}^p\right)^{1/p},&1\le p<\infty,\\
\max_{k=1,\cdots,d}s_{k},&p=\infty.
\end{cases}
\]
If ${U_i}^0$ and ${\Omega_i}$ are simultaneously diagonalizable for each fixed $i$, for instance if the ${\Omega_i}$'s are constant multiples of $I_d$, we may without loss of generality take 
\[ {U_i}(t)=\operatorname{diag}(e^{\mathrm{i}\theta_{i1}},\cdots,e^{\mathrm{i}\theta_{id}}), \quad  \Omega_i=\operatorname{diag}(\mathrm{i}\omega_{i1},\cdots,\mathrm{i}\omega_{id}),  \]
and \eqref{D-2} becomes
\begin{equation}\notag
\dot\theta_{ij}=\omega_{ij}-\frac \kappa {Nd}\sum_{k=1}^N\sum_{l=1}^d(1+\cos\theta_{kl})\sin\theta_{ij},\quad i \in [N], \quad j \in [d],
\end{equation}
namely, it becomes the Winfree model \eqref{A-1} over the eigenvalues of $U_1(t),\cdots, U_N(t)$. Conversely, we may as well view \eqref{A-1} as the evolution of the eigenvalues of the single matrix $U_1$ in \eqref{D-2} with $N=1$. We will see that in general, the eigenvalues of $U_1(t),\cdots,U_N(t)$ satisfy the same differential \emph{inequality} as that of \eqref{A-1}. Throughout this section, we use $\succeq$ as a comparison operator between matrices:
\[A\succeq B \quad \Longleftrightarrow \quad A-B \text{ is positive semidefinite}.\]
Similarly, we denote
\[A \succ B \quad \Longleftrightarrow \quad A-B \text{ is positive definite}.\]
\subsection{Finite-time entrance into the stable region}\label{subsec:4-1}
For the $\mathrm{U}(d)$-model, one can introduce the order parameter:
\begin{equation}\label{D-4}
R\coloneqq \frac 1N \sum_{j=1}^N \frac{d+\operatorname{Re}\operatorname{tr}{U_j}}{d},
\end{equation}
which measures how close the matrices $U_i$ are to $I_d$ and how far they are from $-I_d$. In this case, \eqref{D-2} becomes
\begin{equation}\label{D-5}
\dot{U}_i{U_i}^{\dagger}={\Omega_i}+\frac {\kappa R}{2}({U_i}^{\dagger}-{U_i}).
\end{equation}
Again, global well-posedness follows by the same argument as in Section \ref{sec:3}. For $\alpha\in\left(\frac{\pi}{3},\pi\right)$, let $\alpha^\infty\in\left(0,\frac{\pi}{3}\right)$ be the unique solution of
\[
\sin\alpha^\infty(1+\cos\alpha^\infty)=\sin\alpha(1+\cos\alpha).
\]
We set
\begin{align}
\begin{aligned} \label{D-5-1}
&T_2^*\coloneqq\frac{\pi}{\kappa\sin\alpha(1+\cos\alpha)-\|{\mathfrak V}\|_\infty}, \quad \beta\in(0,\pi], \\
&O(\mathrm{U}(d))_{\beta}\coloneqq
\left\{U\in\mathrm{U}(d):\operatorname{Re}\lambda>\cos\beta\ \text{for every eigenvalue $\lambda$ of $U$}\right\}.
\end{aligned}
\end{align}
In what follows, we provide several elementary lemmas to be used later.
\begin{lemma}\label{L4.1}
For $U,V\in O(\mathrm{U}(d))_{\alpha^\infty}$, we have
\[
\operatorname{Re}\operatorname{tr}(U-V)\leq\sin\alpha^\infty\|U-V\|_{S_1^d}.
\]
\end{lemma}
\begin{proof}
We provide its proof in Appendix~\ref{app:1}.
\end{proof}
Next, we derive a differential inequality for the real parts of the eigenvalues of $U_j(t)$. For this, we recall an instance of the Pythagorean theorem.
\begin{lemma}\label{L4.2}
For a unitary matrix $U$, let $H$ and $K$ be its Hermitian and skew-Hermitian parts, respectively:
\[H \coloneqq \frac{U + U^\dagger}2, \quad K\coloneqq\frac{U - U^\dagger}2.\]
Then, for any vector $x\in\mathbb C^d$, we have
\[\|x\|^2 = \|H x\|^2 + \|K x\|^2.\]
\end{lemma}
\begin{proof}
While the proof can be easily found in basic linear algebra textbooks, we provide the full proof in Appendix \ref{app:2}.
\end{proof}

\begin{lemma}\label{L4.3}
Let $\{ U_i(t) \}$ be a global solution to \eqref{D-3}, and denote the eigenvalues of $U_i(t)$ by $\lambda_{ik}(t)$, $k\in[d]$, ordered by decreasing real part:
\[\mathrm{Re}\lambda_{i1}(t)\ge \cdots \ge \mathrm{Re}\lambda_{ik}(t)\ge\cdots\ge \mathrm{Re}\lambda_{id}(t).\]
Then, the following assertions hold.
\begin{enumerate}
\item
For each $k$, the real-valued function $t\mapsto\operatorname{Re}\lambda_{ik}(t)$ is right-differentiable.
\vspace{0.1cm}
\item
The real part of the eigenvalue $\lambda_{ik}$ satisfies the following differential inequality:
\begin{equation*}\label{D-6}
D^+ (\mathrm{Re}\lambda_{ik}(t)) \ge -\|\Omega_i\|_{\mathrm{op}}\sqrt{1-(\mathrm{Re}\lambda_{ik}(t))^2}+\kappa R(1-(\mathrm{Re}\lambda_{ik}(t))^2).
\end{equation*}
\end{enumerate}
\end{lemma}
\begin{proof}
(1)
First, we define 
\[H_i(t)\coloneqq \frac{1}{2}(U_i(t) + U_i(t)^\dagger).\]
Since $U_i$ is unitary, the $k$-th largest eigenvalue of $H_i$ is $\operatorname{Re}\lambda_{ik}$. By \cite[Theorem 3.3]{hiriart1999clarke}, the ordered eigenvalue functions have one-sided directional derivatives; the same compression formula for Hermitian matrices follows by realification, which repeats each eigenvalue twice. Since $H_i(t+h)=H_i(t)+h\dot H_i(t)+o(h)$ and ordered Hermitian eigenvalues are Lipschitz continuous in the matrix argument, the ordinary right derivative exists at every time, including eigenvalue crossings. Therefore, it coincides with the upper right Dini derivative, and we obtain
\[D^+{\mathrm{Re}\lambda_{ik}(t)} = (k-\hat k+1)\text{-th largest eigenvalue of } V_{i,k}^\dagger \dot{H}_i V_{i, k},\]
where
\[
\hat{k}=\min\{l:\operatorname{Re}\lambda_{il}=\operatorname{Re}\lambda_{ik}\},\quad\bar{k}=\max\{l:\operatorname{Re}\lambda_{il}=\operatorname{Re}\lambda_{ik}\},
\]
and $V_i$ is a unitary matrix that diagonalizes $H_i$:
\[
V_i^\dagger H_i V_i=\operatorname{diag}(\operatorname{Re}\lambda_{i1},\cdots,\operatorname{Re}\lambda_{id}),
\]
and $V_{i,k}$ is the submatrix of $V_i$ consisting of the columns indexed by $\hat{k},\hat{k}+1,\cdots,\bar{k}$. \newline

\noindent (2)~Next, we have to analyze the eigenvalues of $V_{i, k}^\dagger \dot H_i V_{i, k}$. Let $(\eta,y)$ be any eigenpair of $V_{i,k}^\dagger\dot H_iV_{i,k}$ with $\|y\|=1$, so that
\[
(V_{i,k}^\dagger\dot H_iV_{i,k})y=\eta y,
\quad \|y\|=1.
\]
Then, we have
\[
(V_{i,k}y)^\dagger\dot H_i(V_{i,k}y)=\eta.
\]
Hence, it is enough to establish an inequality on $x^\dagger \dot{H}_i x$, where $x\in \bbc^d$ is a linear combination of the columns of $V_{i,k}$. Moreover, since $V_{i, k}$ satisfies 
\[V_{i,k}^\dagger V_{i,k} = I_{\bar k - \hat k +1},\]
we may assume that $x$ is a unit vector. Accordingly, let $x_k$ be any unit vector in the column space of $V_{i,k}$. We estimate $x_k^\dagger\dot H_i x_k$ below. Note that, from the definition of $V_{i, k}$, $\hat k$, and $\bar k$, we also have
\begin{equation}\label{D-7}
{H}_i x_k = (\mathrm{Re}\lambda_{ik})x_k.
\end{equation}

Next, we derive the dynamics of $H_i$. From  \eqref{D-3}, we have
\[\dot U_i = \Omega_i U_i  + \frac {\kappa R}2 (I_d - U_i^2).\]
Consequently, we have the following differential equation for $H_i$:
\begin{equation}\label{D-8}
\begin{aligned}
\dot{H}_i &= \frac{1}{2}(\dot U_i + \dot U_i ^\dagger) \\
&= \frac12\bigl(\Omega_iU_i+U_i^\dagger\Omega_i^\dagger\bigr)
+ \frac{\kappa R}{2}\left(I_d-\frac{U_i^2 + (U_i^\dagger)^2}{2} \right) \\
&=\frac12\bigl(\Omega_iU_i+U_i^\dagger\Omega_i^\dagger\bigr)
+\kappa R\bigl(I_d- H_i^2\bigr).
\end{aligned}
\end{equation}
We use the following identity in the last equality:
\[H_i^2 = \frac 1 4 (U_i^2 + (U_i^\dagger)^2 + 2I_d) \quad \Longleftrightarrow \quad I_d - H_i^2 = \frac{1}{2}\left(I_d - \frac{U_i^2 + (U_i^\dagger)^2}{2}\right). \]
We combine \eqref{D-7} and \eqref{D-8} to obtain
\begin{equation}\label{D-9}
\begin{aligned}
    x_k^\dagger \dot H_i x_k &= \frac 1 2 x_k^\dagger (\Omega_i U_i + U_i^\dagger \Omega_i^\dagger)x_k + \kappa R (1 - (\mathrm{Re}\lambda_{ik})^2) \\
    &= \mathrm{Re}(x_k^\dagger \Omega_i U_i x_k) + \kappa R (1 - (\mathrm{Re}\lambda_{ik})^2).
\end{aligned}
\end{equation}
Finally, we estimate the first term of \eqref{D-9} as follows. For this, we define
\[K_i \coloneqq \frac 1 2 (U_i - U_i^\dagger).\]
Then, since $\Omega_i$ is skew-Hermitian, we have
\begin{equation}\label{D-10}
\mathrm{Re}(x_k^\dagger \Omega_i H_i x_k) = \mathrm{Re}\lambda_{ik}\mathrm{Re}(x_k^\dagger \Omega_i x_k) = 0.
\end{equation}
It follows from \eqref{D-10}, Lemma~\ref{L4.2}, the Cauchy--Schwarz inequality, and the definition of $x_k$ that 
\begin{equation}\label{D-11}
    \begin{aligned}
        \mathrm{Re}(x_k^\dagger \Omega_i U_i x_k) &= \mathrm{Re}(x_k^\dagger \Omega_i H_i x_k) + \mathrm{Re}(x_k^\dagger \Omega_i K_i x_k) \\
        &= \mathrm{Re}\big((x_k^\dagger \Omega_i)\cdot (K_i x_k)\big) \\
        &\ge -\|x_k^\dagger \Omega_i\|\cdot \|K_i x_k\| \\
        &\ge -\|\Omega_i\|_{\mathrm{op}} \sqrt{1 - (\mathrm{Re}\lambda_{ik})^2 }.
    \end{aligned}
\end{equation}
Finally, we combine \eqref{D-9} and \eqref{D-11} to obtain the desired result:
\[D^+ (\mathrm{Re} \lambda_{ik}) \ge -\|\Omega_i\|_{\mathrm{op}} \sqrt{1 - (\mathrm{Re}\lambda_{ik})^2 } + \kappa R (1 - (\mathrm{Re}\lambda_{ik})^2).\]
\end{proof}
Next, we establish finite-time entrance into the invariant set $O(\mathrm{U}(d))_{\alpha^\infty}^N$. 
\begin{lemma} \label{L4.4}
Suppose that parameters and initial data satisfy
\[
\alpha\in\left(\frac{\pi}{3},\pi\right),\quad
\kappa>\frac{\|{\mathfrak V}\|_\infty}{\sin\alpha(1+\cos\alpha)},\quad
U_i^0\in\overline{O(\mathrm{U}(d))_{\alpha}},\quad i\in[N],
\]
and let ${\mathcal U}(t)= \{ U_i(t) \}_{i=1}^{N}$ be the global solution to \eqref{D-3}. Then there exists $T<T_2^*$ such that
\[
{\mathcal U}(t)\in O(\mathrm{U}(d))_{\alpha^\infty}^N,\quad t\geq T,
\]
where $T_2^*$ is a positive constant defined in \eqref{D-5-1}.
\end{lemma}
\begin{proof} 

The proof consists of two steps. First, we derive a differential inequality for the maximal spectral angle. Second, we use the unimodality of $(1+\cos\theta)\sin\theta$ to obtain finite-time entrance.

\vspace{0.5em}

\noindent\textbf{$\bullet$ Step A} (Differential inequality for the maximal spectral angle):
For each $i\in[N]$, denote the eigenvalues of $U_i(t)$ by $\lambda_{ik}(t)$, $k\in[d]$, ordered by decreasing real part:
\[\mathrm{Re}\lambda_{i1}(t)\ge \cdots \ge \mathrm{Re}\lambda_{ik}(t)\ge\cdots\ge \mathrm{Re}\lambda_{id}(t).\]
Now we define 
\[\rho_i(t) = \cos^{-1} \big(\mathrm{Re}\lambda_{id}(t)\big), \quad \rho(t)\coloneqq \max_{i\in [N]}\rho_i(t).\]
This definition yields, for each $\beta \in (0, \pi)$,
\[U_i(t) \in O(\mathrm{U}(d))_\beta \quad \Longleftrightarrow\quad \rho_i(t) < \beta, \]
and
\[\mathcal{U}(t)\in O(\mathrm{U}(d))_\beta^N \quad \Longleftrightarrow \quad \rho(t) < \beta. \]
We apply Lemma~\ref{L4.3} with $k=d$ to obtain
\begin{equation}\label{D-12}
D^+ \mathrm{Re}\lambda_{id}(t)\ge -\|\Omega_i\|_{\mathrm{op}} \sqrt{1 - (\mathrm{Re}\lambda_{id})^2 } + \kappa R (1 - (\mathrm{Re}\lambda_{id})^2).
\end{equation}
For $\rho_i(t)\in(0,\pi)$, Lemma~\ref{L4.3} and the chain rule for ordinary right derivatives give
\begin{equation}\label{D-13}
D^+\rho_i(t)
\leq\|\Omega_i\|_{\mathrm{op}}-\kappa R(t)\sin \rho_i(t).
\end{equation}
If $\rho_i(t)=0$, then $U_i(t)=I_d$ and \eqref{D-5} gives $\dot U_i(t)=\Omega_i$.  Now, we use
\[
\|U_i(t+h)-I_d\|_{\mathrm{op}}=2\sin\frac{\rho_i(t+h)}{2} \quad \mbox{and} \quad U_i(t+h)=I_d+h\Omega_i+o(h) 
\]
to obtain 
\[ D^+\rho_i(t)=\|\Omega_i\|_{\mathrm{op}}. \]
So the relation \eqref{D-13} also holds when $\rho_i(t)=0$. On the other hand, at $\rho_i(t)=\pi$, the bound $\rho_i\leq\pi$ yields 
\[ D^+\rho_i(t)\leq0\leq\|\Omega_i\|_{\mathrm{op}}. \]
Thus, the relation \eqref{D-13} holds at both endpoints as well as through eigenvalue multiplicities and crossings. The definition of $\rho(t)$ implies
\[\mathrm{Re}\lambda_{ik}(t) \ge \cos\rho(t), \quad \forall\, i\in [N],\ k\in [d]. \]
Hence, we have
\[
R(t)=1+\frac{1}{Nd}\sum_{j=1}^N\operatorname{Re}\operatorname{tr}U_j(t)
= 1 + \frac1{Nd} \sum_{j=1}^{N} \sum_{k=1}^d \mathrm{Re}\lambda_{jk}(t)
\geq1+\cos \rho(t).
\]
Using \eqref{D-13} and the upper right Dini derivative rule for a finite maximum, we obtain
\begin{equation}\label{D-14}
D^+ \rho(t)
\leq\|{\mathfrak V}\|_\infty-\kappa(1+\cos \rho(t))\sin \rho(t).
\end{equation}

\vspace{0.5em}

\noindent\textbf{$\bullet$ Step B} (Entrance into $O(\mathrm{U}(d))_{\alpha^\infty}^N$): We set 
\[
g(\theta)\coloneqq(1+\cos\theta)\sin\theta.
\]
The function $g$ is increasing on $(0,\pi/3]$ and decreasing on $[\pi/3,\pi)$, and $g(\alpha^\infty)=g(\alpha)$. Thus we have 
\[ g(\theta)\geq g(\alpha) \quad \mbox{for $\theta\in[\alpha^\infty,\alpha]$}. \]
It follows from \eqref{D-14} that, whenever $\rho(t)\in[\alpha^\infty,\alpha]$,
\begin{equation}\label{D-15}
D^+\rho(t)
\leq-\left(\kappa\sin\alpha(1+\cos\alpha)-\|\mathfrak V\|_\infty\right)<0.
\end{equation}
Because $U_i(0)\in\overline{O(\mathrm{U}(d))_\alpha}$, we have
\[ \rho(0)\leq\alpha. \]
The strict inequality in \eqref{D-15} at $\rho=\alpha$ prevents an outward crossing. If $\rho(0)<\alpha^\infty$, we take $T=0$; otherwise, we set
\[
\tau\coloneqq\inf\{t\geq0:\rho(t)\leq\alpha^\infty\}.
\]
While $\rho(t)\in[\alpha^\infty,\alpha]$, the relation \eqref{D-15} gives
\[
\tau
\leq\frac{\rho(0)-\alpha^\infty}{\kappa\sin\alpha(1+\cos\alpha)-\|{\mathfrak V}\|_\infty}
\leq\frac{\alpha-\alpha^\infty}{\kappa\sin\alpha(1+\cos\alpha)-\|{\mathfrak V}\|_\infty}
<T_2^*.
\]
At $\rho=\alpha^\infty$, the strict inequality \eqref{D-15} still holds because $g(\alpha^\infty)=g(\alpha)$. Hence continuity and the definition of $D^+$ provide $\varepsilon>0$ such that
\[
\rho(t)<\alpha^\infty,\quad t\in(\tau,\tau+\varepsilon].
\]
Since $\tau<T_2^*$, choose
\[
T\in\bigl(\tau,\min\{\tau+\varepsilon,T_2^*\}\bigr).
\]
Then, we have
\[ \rho(T)<\alpha^\infty. \]
The same strict inequality in \eqref{D-15} at the boundary rules out the first exit, so
\[
\rho(t)<\alpha^\infty,\quad t\geq T.
\]
Equivalently, we have
\[ {\mathcal U}(t)\in O(\mathrm{U}(d))_{\alpha^\infty}^N \quad \mbox{for every $t\geq T$}. \]
\end{proof}
\begin{remark} \label{R4.1}
If ${\mathcal U}^0\in O(\mathrm{U}(d))_{\alpha^\infty}^N$, then one may take $T=0$.
\end{remark}
\subsection{Exponential stability}\label{subsec:4-2}
In this subsection, we establish exponential stability for \eqref{D-3} with respect to initial data. For this, we define a mixed norm for an $N$-tuple ${\mathcal C}=(C_1,\ldots,C_N)$ of $d\times d$ complex matrices:
\[
\|{\mathcal C}\|_{\ell_q^N(S_p^d)}\coloneqq
\begin{cases}
\displaystyle\left(\frac{1}{N}\sum_{i=1}^N\|C_i\|_{S_p^d}^q\right)^{1/q},&1\leq q<\infty,\\
\displaystyle\max_{i\in[N]}\|C_i\|_{S_p^d},&q=\infty.
\end{cases}
\]
Recall the positive coefficient $\mu_{\alpha^\infty}$ used in Section~\ref{subsec:3-2}:
\[\mu_{\alpha^\infty} = (2\cos\alpha^\infty - 1)(1+\cos\alpha^\infty).\]
\begin{lemma} \label{L4.5}
Suppose that parameters and initial data satisfy 
\[
1\leq p\leq\infty, \quad 1\leq q\leq\infty,\quad
\alpha\in\left(\frac{\pi}{3},\pi\right),\quad
\kappa>\frac{\|{\mathfrak V}\|_\infty}{\sin\alpha(1+\cos\alpha)},
\]
and let ${\mathcal U}$ and $\widetilde{\mathcal U}$ be the global solutions to \eqref{D-3} with initial data satisfying
\[
U_i^0,~~\widetilde U_i^0\in\overline{O(\mathrm{U}(d))_\alpha},\quad \forall~i\in[N].
\] Then there exists $T<T_2^*$ such that
\[
\|{\mathcal U}(t)-\widetilde{\mathcal U}(t)\|_{\ell_q^N(S_p^d)}
\leq e^{-\kappa\mu_{\alpha^\infty}(t-T)}
\|{\mathcal U}(T)-\widetilde{\mathcal U}(T)\|_{\ell_q^N(S_p^d)},
\quad t\geq T.
\]
\end{lemma}
\begin{proof}
By Lemma~\ref{L4.4}, each solution enters $O(\mathrm{U}(d))_{\alpha^\infty}^N$ before the instant $T_2^*$. Let $T$ be the larger one of the two entrance times, and for $t\geq T$ we set
\[
W_i(t)\coloneqq U_i(t)-\widetilde U_i(t),\quad
\mathcal W(t)\coloneqq(W_1(t),\ldots,W_N(t)).
\]
The proof consists of four steps in what follows: we derive the difference equation, estimate the Schatten norm, control the order-parameter difference, and pass to the mixed norm. It suffices to prove
\begin{equation}\label{D-16}
D^+\|\mathcal W(t)\|_{\ell_q^N(S_p^d)}
\leq-\kappa\mu_{\alpha^\infty}\|\mathcal W(t)\|_{\ell_q^N(S_p^d)},\quad t\geq T.
\end{equation}

\noindent\textbf{$\bullet$ Step A} (Difference equation and stable-region bounds):
The crucial noncommutative identity is
\[
U_i^2-\widetilde U_i^2=U_iW_i+W_i\widetilde U_i,
\]
where $W_i$ was defined above. Let $R$ and $\widetilde R$ be the order parameters associated with $\mathcal U$ and $\widetilde{\mathcal U}$, respectively, as in \eqref{D-4}. Subtracting the equations for $U_i$ and $\widetilde U_i$ in \eqref{D-3} yields
\begin{equation}\label{D-17}
\dot W_i
=\Omega_iW_i
-\frac{\kappa R}{2}\bigl(U_iW_i+W_i\widetilde U_i\bigr)
+\frac{\kappa(R-\widetilde R)}{2}\bigl(I_d-\widetilde U_i^2\bigr).
\end{equation}
For $t\geq T$, we have
\begin{equation}\label{D-18}
\frac{U_i+U_i^\dagger}{2}\succeq \cos\alpha^\infty I_d,
\quad
\frac{\widetilde U_i+\widetilde U_i^\dagger}{2}\succeq \cos\alpha^\infty I_d,
\quad
R\geq1+\cos\alpha^\infty.
\end{equation}
Moreover, we have
\[
\|I_d-\widetilde U_i^2\|_{\mathrm{op}}\leq2\sin\alpha^\infty.
\]
\noindent\textbf{$\bullet$ Step B} (Schatten norm estimate):
Fix $t\geq T$, we set $a=\alpha^\infty$, and evaluate all coefficients at $t$, and rewrite \eqref{D-17} as $\dot W_i=A_iW_i+W_iB_i+F_i$, where
\[
A_i=\Omega_i-\tfrac{\kappa R}{2}U_i,\qquad
B_i=-\tfrac{\kappa R}{2}\widetilde U_i,\qquad
F_i=\tfrac{\kappa(R-\widetilde R)}2(I_d-\widetilde U_i^2).
\]
Set $\gamma=\kappa R\cos(a)/2$. By \eqref{D-18} and the
skew-Hermiticity of $\Omega_i$, the Hermitian parts of $A_i$ and $B_i$
are bounded above by $-\gamma I_d$. Differentiating $\|e^{sA_i}v\|^2$
and $\|e^{sB_i}v\|^2$ gives
\begin{equation}\label{D-19}
	\|e^{hA_i}\|_{\mathrm{op}},\ \|e^{hB_i}\|_{\mathrm{op}}
	\leq e^{-\gamma h},\qquad h\geq0.
\end{equation}
To the first order as $h\downarrow0$,
\[
W_i(t+h)=e^{hA_i}W_i(t)e^{hB_i}+hF_i+o(h).
\]
For $1\leq p\leq\infty$, the Schatten ideal inequality yields
\[ \|P W Q\|_{S_p^d}\leq\|P\|_{\mathrm{op}}\|W\|_{S_p^d}\|Q\|_{\mathrm{op}}. \]
Since $\|I_d-\widetilde U_i^2\|_{S_p^d}\leq2d^{1/p}\sin a$,
taking the upper right difference quotient gives
\begin{equation}\label{D-20}
	D^+\|W_i\|_{S_p^d}\leq
	-\kappa R\cos(a)\|W_i\|_{S_p^d}
	+\kappa d^{1/p}\sin(a)|R-\widetilde R|,
\end{equation}
with the convention $1/\infty=0$. This also holds when $W_i(t)=0$.
\newline

\noindent\textbf{$\bullet$ Step C} (Difference of order parameters):
By \eqref{D-4} and Lemma~\ref{L4.1} applied in both directions, we have
\begin{equation} \label{D-21}
\hspace{-0.5cm} |R-\widetilde R| \leq\frac{\sin\alpha^\infty}{Nd}\sum_{j=1}^N\|W_j\|_{S_1^d} \leq\sin\alpha^\infty d^{-1/p}\frac1N\sum_{j=1}^N\|W_j\|_{S_p^d} =\sin\alpha^\infty d^{-1/p}\|\mathcal W\|_{\ell_1^N(S_p^d)}.
\end{equation}
Combining \eqref{D-18}, \eqref{D-20}, and \eqref{D-21}, we obtain
\begin{equation}\label{D-22}
D^+\|W_i\|_{S_p^d}
\leq-\kappa\cos\alpha^\infty(1+\cos\alpha^\infty)\|W_i\|_{S_p^d}
+\frac{\kappa\sin^2\alpha^\infty}{N}\sum_{j=1}^N\|W_j\|_{S_p^d}.
\end{equation}

\vspace{0.5em}

\noindent\textbf{$\bullet$ Step D} (Mixed-norm estimate):~
For $q=1$, averaging \eqref{D-22} over $i$ and using the
subadditivity of the upper right Dini derivative gives
\[
D^+\left(\frac1N\sum_{j=1}^N\|W_j\|_{S_p^d}\right)
\leq-\kappa\mu_{\alpha^\infty}\frac1N\sum_{j=1}^N\|W_j\|_{S_p^d}.
\]
For $1\leq q<\infty$, we set
\[
Y_q(t)\coloneqq\left(\frac1N\sum_{j=1}^N\|W_j(t)\|_{S_p^d}^q\right)^{1/q}.
\]
So the above inequality reads
\[
D^+Y_1 \le -\kappa \mu_{\alpha^\infty} Y_1.
\]
For $1<q<\infty$, at times with $Y_q(t)>0$, \eqref{D-22}, H\"older's inequality, and $Y_1\leq Y_q$ imply
\begin{align*}
	D^+Y_q
	&\leq-\kappa\cos\alpha^\infty(1+\cos\alpha^\infty)Y_q
	+\kappa\sin^2\alpha^\infty Y_1Y_q^{1-q}\frac1N\sum_{i=1}^N\|W_i\|_{S_p^d}^{q-1}\\
	&\leq-\kappa\bigl(\cos\alpha^\infty(1+\cos\alpha^\infty)-\sin^2\alpha^\infty\bigr)Y_q
	=-\kappa\mu_{\alpha^\infty}Y_q.
\end{align*}
When $Y_q(t)=0$, the two full states agree at time $t$. Uniqueness of solutions then gives 
\[ \mathcal W(s)=0 \quad \mbox{for all $s\geq t$}, \]
so we have
\[ D^+Y_q(t)=0. \]

For $q=\infty$, we set $Y_\infty(t)=\max_i\|W_i(t)\|_{S_p^d}$.
At an index attaining the maximum, \eqref{D-22} and
$Y_1\leq Y_\infty$ bound the component's upper right Dini derivative
by $-\kappa\mu_{\alpha^\infty}Y_\infty$.
Let $\mathcal I(t)=\{i:\|W_i(t)\|_{S_p^d}=Y_\infty(t)\}$.
Since the component norms are continuous and $N$ is finite,
\[
D^+Y_\infty(t)
=\max_{i\in\mathcal I(t)}D^+\|W_i(t)\|_{S_p^d}
\leq -\kappa\mu_{\alpha^\infty}Y_\infty(t),
\]
where the inequality follows from \eqref{D-22} and $Y_1\leq Y_\infty$. Gronwall's lemma now yields the asserted contraction for all
$1\leq p,q\leq\infty$.
\end{proof}
\begin{remark} \label{R4.2}
If both initial states belong to $O(\mathrm{U}(d))_{\alpha^\infty}^N$, then one may take $T=0$.
\end{remark}

\subsection{Relaxation dynamics}\label{subsec:4-3}
In this subsection, we use Lemma~\ref{L4.5} and compare a solution with its time translates to show the existence and uniqueness of an equilibrium in $O(\mathrm{U}(d))_{\alpha^\infty}^N$ and exponential relaxation toward it.
\begin{proposition} \label{P4.1}
Suppose that the system parameters satisfy
\[
\alpha\in\left(\frac{\pi}{3},\pi\right),\quad
\kappa>\frac{\|{\mathfrak V}\|_\infty}{\sin\alpha(1+\cos\alpha)}.
\]
Then system \eqref{D-3} admits a unique equilibrium ${\mathcal U}^\infty\in O(\mathrm{U}(d))_{\alpha^\infty}^N$.
\end{proposition}
\begin{proof} 
We choose ${\mathcal U}^0\in O(\mathrm{U}(d))_{\alpha^\infty}^N$, and let ${\mathcal U}(t)$ be the corresponding global solution. Now, we use Lemma~\ref{L4.4} with $T=0$ to see that this region is positively invariant. For $s\geq0$, we define
\[
\widetilde{\mathcal U}(t)\coloneqq{\mathcal U}(t+s).
\]
The system is autonomous, so $\widetilde{\mathcal U}$ is another solution in the same stable region. Lemma~\ref{L4.5} with $p=q=2$ and $T=0$ yields
\begin{equation}\label{D-23}
\|{\mathcal U}(t+s)-{\mathcal U}(t)\|_{\ell_2^N(S_2^d)}
\leq e^{-\kappa\mu_{\alpha^\infty}t}
\|{\mathcal U}(s)-{\mathcal U}(0)\|_{\ell_2^N(S_2^d)}.
\end{equation}
Because $\mathrm{U}(d)^N$ is compact, the last factor in \eqref{D-23} is bounded uniformly in $s$. For $t_2\geq t_1\geq0$, we take 
\[ t=t_1 \quad \mbox{and} \quad s=t_2-t_1 \quad \mbox{in \eqref{D-23}} \]
to obtain
\[
\|{\mathcal U}(t_2)-{\mathcal U}(t_1)\|_{\ell_2^N(S_2^d)}
\leq Ce^{-\kappa\mu_{\alpha^\infty}t_1},
\]
where $C$ is independent of $t_1$ and $t_2$. Thus ${\mathcal U}(t)$ is Cauchy and converges to some ${\mathcal U}^\infty\in\mathrm{U}(d)^N$. If $\Phi_h$ denotes the flow of \eqref{D-3}, then for every $h\geq0$,
\[
\Phi_h({\mathcal U}^\infty)
=\lim_{t\to\infty}\Phi_h({\mathcal U}(t))
=\lim_{t\to\infty}{\mathcal U}(t+h)
={\mathcal U}^\infty.
\]
Hence, ${\mathcal U}^\infty$ is an equilibrium.
It remains to show that the limiting equilibrium belongs to the open stable region. Continuity first gives
\[
{\mathcal U}^\infty\in\overline{O(\mathrm{U}(d))_{\alpha^\infty}}^N.
\]
Suppose that for some $i \in [N]$, the Hermitian part
\[
H_i^\infty:=\frac{U_i^\infty+U_i^{\infty\dagger}}{2}
\]
has the smallest eigenvalue $\cos\alpha^\infty$. We choose a unit vector $x$ in that eigenspace.
Then $H_i^\infty x=\cos\alpha^\infty x$, and Lemma~\ref{L4.2} gives the exact identity
\[
\Big \|\frac{U_i^\infty-U_i^{\infty\dagger}}{2} x \Big \|^2=1-\|H_i^\infty x\|^2
=\sin^2\alpha^\infty.
\]
Also, skew-Hermiticity of $\Omega_i$ implies
\[
\operatorname{Re}(x^\dagger\Omega_iH_i^\infty x)
=\cos\alpha^\infty\operatorname{Re}(x^\dagger\Omega_ix)=0.
\]
Evaluating \eqref{D-8} at the equilibrium and using these identities, we obtain
\begin{align*}
0=x^\dagger\dot H_ix
&=\operatorname{Re}(x^\dagger\Omega_i\frac{U_i^\infty-U_i^{\infty\dagger}}{2} x)
+\kappa R(\mathcal U^\infty)\sin^2\alpha^\infty\\
&\geq-\|\Omega_i\|_{\mathrm{op}}\sin\alpha^\infty
+\kappa(1+\cos\alpha^\infty)\sin^2\alpha^\infty\\
&\geq\sin\alpha^\infty
\left[\kappa\sin\alpha(1+\cos\alpha)-\|\mathfrak V\|_\infty\right]>0.
\end{align*}
The first inequality follows from the Cauchy--Schwarz inequality and $R(\mathcal U^\infty)\geq1+\cos\alpha^\infty$; the second uses \eqref{B-2}. This contradiction shows that $\mathcal U^\infty\in O(\mathrm{U}(d))_{\alpha^\infty}^N$.

Finally, let ${\mathcal U}^\infty$ and ${\mathcal V}^\infty$ be two equilibria in $O(\mathrm{U}(d))_{\alpha^\infty}^N$. Lemma~\ref{L4.5} with $p=q=2$ and $T=0$ applied to two constant solutions gives
\[
\|{\mathcal U}^\infty-{\mathcal V}^\infty\|_{\ell_2^N(S_2^d)}
\leq e^{-\kappa\mu_{\alpha^\infty}t}
\|{\mathcal U}^\infty-{\mathcal V}^\infty\|_{\ell_2^N(S_2^d)},
\quad t\geq0.
\]
Letting $t\to\infty$ yields 
\[ {\mathcal U}^\infty={\mathcal V}^\infty. \]
\end{proof}
\begin{remark}\label{R4.3}
Note that the uniqueness assertion is restricted to $O(\mathrm U(d))_{\alpha^\infty}^N$. Global uniqueness fails in general: when all $\Omega_i=0$, any configuration of unitary involutions $U_i^2=I_d$ is an equilibrium. The contraction argument gives precisely the uniqueness needed for the relaxation theorem.
\end{remark}
Next, we show the exponential convergence toward the equilibrium whose existence is guaranteed by Proposition \ref{P4.1}. 

\begin{theorem}  \label{T4.1}
Suppose that parameters and initial data satisfy 
\[
\begin{gathered}
1\leq p\leq\infty, \quad 1\leq q\leq\infty,\quad \alpha\in\left(\frac{\pi}{3},\pi\right),\\
\kappa>\frac{\|{\mathfrak V}\|_\infty}{\sin\alpha(1+\cos\alpha)},\quad U_i^0\in\overline{O(\mathrm{U}(d))_\alpha},\quad i\in[N],
\end{gathered}
\]
and let ${\mathcal U}(t)$ and~${\mathcal U}^\infty\in O(\mathrm{U}(d))_{\alpha^\infty}^N$ be the global solution to \eqref{D-3} and the unique equilibrium from Proposition~\ref{P4.1}, respectively. Then there exists $T<T_2^*$ such that
\begin{equation}\label{D-24}
\|{\mathcal U}(t)-{\mathcal U}^\infty\|_{\ell_q^N(S_p^d)}
\leq e^{-\kappa\mu_{\alpha^\infty}(t-T)}
\|{\mathcal U}(T)-{\mathcal U}^\infty\|_{\ell_q^N(S_p^d)},
\quad t\geq T.
\end{equation}
\end{theorem}
\begin{proof}
Note that Lemma~\ref{L4.4} provides $T<T_2^*$ such that 
\[ {\mathcal U}(t)\in O(\mathrm{U}(d))_{\alpha^\infty}^N \quad \mbox{for all $t\geq T$.}  \]
The equilibrium ${\mathcal U}^\infty$ generates a constant solution in the same region. Applying Lemma~\ref{L4.5} to the shifted solution ${\mathcal U}(T+\cdot)$ and the constant equilibrium gives \eqref{D-24}. If the initial state is already in the stable region, Lemma~\ref{L4.4} allows $T=0$.
\end{proof}

\begin{remark}\label{R4.4}
We make the following remarks on Theorem \ref{T4.1}:
\begin{enumerate}
\item
If $\mathcal U^0\in O(\mathrm{U}(d))_{\alpha^\infty}^N$,
then one may take $T=0$.
\vspace{0.1cm}
\item 
As mentioned in the introduction, Theorem \ref{T4.1} shows that every initial configuration whose constituent matrices have no eigenvalue equal to $-1$ converges exponentially to an equilibrium, provided the coupling strength is sufficiently large, with a threshold depending on the initial configuration and natural frequency matrices.
\vspace{0.1cm}
\item
One could also consider the unitary model \eqref{D-3} with a general influence function $I:\mathrm{U}(d)\to[0,\infty)$ in place of $1+(\mathrm{Re}\,\mathrm{tr}\,U_j)/d$.
Such a generalization is studied in Section~\ref{sec:5}, where exponential relaxation is established under the framework $({\mathcal F}_{A}1)$--$({\mathcal F}_{A}3)$ and suitable conditions on the coupling strength and the influence function (Theorem~\ref{T5.1}).

Moreover, this model is locally compatible with
the abstract framework of Section~\ref{sec:7}
(Theorem~\ref{T6.1}),
after choosing a coordinate chart around $I_d$.
To apply Theorem~\ref{T6.1}, however, the coordinate
representatives of the vector field
\[
S(U)=\frac{1}{2}(I_d-U^2)
\]
and of $I$ must satisfy the local assumptions
$({\mathcal F}_{B}1)$--$({\mathcal F}_{B}4)$,
in particular the strong dissipativity condition
\eqref{F-4}.
\end{enumerate}
\end{remark}
\begin{remark}[Critical coupling for the existence of unitary equilibria]
	\label{R:unitary-critical-coupling}
	Let $\mathrm{i}\omega_{i1},\ldots,\mathrm{i}\omega_{id}$ be the
	eigenvalues of $\Omega_i$, counted with multiplicity. At an equilibrium
	of the canonical unitary model, the relation \eqref{D-5} gives
	\[
	\Omega_i=\frac{\kappa R}{2}(U_i-U_i^\dagger).
	\]
	Hence $U_i$ commutes with $\Omega_i$. Since both matrices are normal,
	for each $i$ we may choose a common orthonormal eigenbasis and write
	$e^{\mathrm{i}\theta_{ik}}$ for the eigenvalue of $U_i$ corresponding
	to $\mathrm{i}\omega_{ik}$. In these bases, the equilibrium equations
	and the order parameter \eqref{D-4} become
	\[
	\omega_{ik}=\kappa R\sin\theta_{ik},
	\qquad
	R=1+\frac{1}{Nd}\sum_{j=1}^N\sum_{\ell=1}^d\cos\theta_{j\ell}.
	\]
	These are precisely the equilibrium equations for the scalar Winfree
	model \eqref{A-1} with $Nd$ oscillators and frequency tuple
	$(\omega_{ik})_{i,k}$. Conversely, any solution of these scalar equations
	gives a unitary equilibrium by taking $U_i$ to be diagonal with entries
	$e^{\mathrm{i}\theta_{ik}}$ in an orthonormal eigenbasis of $\Omega_i$.
	
	Consequently, for $\mathfrak V\ne0$ and $\kappa>0$, the unitary model
	admits an equilibrium if and only if $\kappa\geq\kappa_{\mathrm c}$,
	where \cite[Proposition~63]{ryoo2026oscillator} gives
	\[
	\kappa_{\mathrm c}
	=\min_{u\geq\|\mathfrak V\|_\infty}
	\frac{u}{
		1+\dfrac{1}{Nd}\displaystyle\sum_{i=1}^N\sum_{k=1}^d
		\sqrt{1-\omega_{ik}^2/u^2}}.
	\]
	If $\mathfrak V=0$, the identity configuration is an equilibrium for
	every $\kappa\geq0$, so $\kappa_{\mathrm c}=0$.
\end{remark}

\section{Winfree oscillators on unitary groups II}\label{sec:5}
\setcounter{equation}{0}
In this section, we extend the analysis of Section~\ref{sec:4} to a Winfree matrix model with a general influence function. Consider the Cauchy problem:
\begin{equation}\label{E-1}
\begin{cases}
\displaystyle \dot{U}_iU_i^{\dagger}=\Omega_i+\frac{\kappa}{2N}\sum_{j=1}^N I(U_j)\left(U_i^{\dagger}-U_i\right),\quad t>0,\\
\displaystyle U_i\big|_{t=0}=U_i^0\in\mathrm{U}(d),\quad i\in[N],
\end{cases}
\end{equation}
where \(\Omega_1,\ldots,\Omega_N\in\mathbb{C}^{d\times d}\) are skew-Hermitian matrices and \(I:\mathrm{U}(d)\to\mathbb{R}\) is the influence function. For \(\beta\in(0,\pi]\), we recall the stable region
\begin{equation}\notag
O(\mathrm{U}(d))_{\beta}\coloneqq \left\{U\in\mathrm{U}(d):\operatorname{Re}\lambda>\cos\beta\ \text{for every eigenvalue }\lambda\text{ of }U\right\}.
\end{equation}
Next, we impose the following structural conditions on the influence function \(I\).
\vspace{0.1cm}
\begin{itemize}
\item
$({\mathcal F}_{A}1)$: For every \(\beta\in(0,\pi)\), the function \(I\) is Lipschitz continuous on \(O(\mathrm{U}(d))_{\beta}\) with respect to the Schatten--von Neumann \(1\)-norm:
\[
|I(U)-I(\widetilde{U})|\leq \operatorname{Lip}_{O(\mathrm{U}(d))_{\beta}}(I)\|U-\widetilde{U}\|_{S_1^d}.
\]

\item
$({\mathcal F}_{A}2)$: Let \(\lambda_1,\ldots,\lambda_d\in\mathbb{S}^1\) be the eigenvalues of \(U\in\mathrm{U}(d)\), ordered so that
\[
\operatorname{Re}\lambda_1\geq\cdots\geq\operatorname{Re}\lambda_d.
\]
For \(\theta_k\in(-\pi,\pi]\) such that \(\lambda_k=e^{\mathrm{i}\theta_k}\), we assume that there exists a non-increasing continuous function \(J:[0,\pi]\to[0,\infty)\) such that
\[
I(U)\geq J(|\theta_d|).
\]
Since \(|\theta_d|=\max_{k\in[d]}|\theta_k|\), this gives a lower bound in terms of the largest spectral angle.
\vspace{0.1cm}
\item
$({\mathcal F}_{A}3)$: The function
\[
F(\theta)\coloneqq J(\theta)\sin\theta \quad \mbox{is unimodal on \((0,\pi]\)}.
\]
More precisely, there exists a unique \(\theta_*\in(0,\pi)\) such that \(F\) is strictly increasing on \((0,\theta_*]\) and strictly decreasing on \([\theta_*,\pi]\). Thus, for every \(\alpha\in(\theta_*,\pi)\), there exists a unique \(\alpha^\infty\in(0,\theta_*)\) satisfying
\begin{equation}\label{E-2}
J(\alpha)\sin\alpha=J(\alpha^\infty)\sin\alpha^\infty, \quad J(\theta)\sin\theta\geq J(\alpha)\sin\alpha,\quad \theta\in[\alpha^\infty,\alpha].
\end{equation}
\end{itemize}
\begin{remark}\label{R5.1}
The canonical influence function from Section~\ref{sec:4},
\[
I(U)=1+\frac{1}{d}\operatorname{Re}\operatorname{tr}U,
\]
satisfies $({\mathcal F}_{A}1)$--$({\mathcal F}_{A}3)$ with \(J(\theta)=1+\cos\theta\) and \(\theta_*=\pi/3\). Moreover, for $0<a<\pi/3$, Lemma~\ref{L4.1}, applied after interchanging $U$ and $V$ as well, gives
\[
\operatorname{Lip}_{O(\mathrm U(d))_a}(I)=\frac{\sin a}{d}.
\]
The reverse inequality follows by taking diagonal matrices that differ in one entry, with both corresponding angles approaching $a$. Consequently, the coefficient in \eqref{E-5} satisfies, for $a=\alpha^\infty$,
\[
\mu_I=\cos a(1+\cos a)-\sin^2a
=\cos a+\cos(2a)=\mu_a>0,
\]
so the general theorem recovers the contraction rate of Section~\ref{sec:4}.
\end{remark}
\vspace{0.2cm}

As in Section~\ref{sec:4}, we introduce the order parameter:
\begin{equation}\label{E-3}
R\coloneqq\frac{1}{N}\sum_{j=1}^N I(U_j).
\end{equation}
Then, the equation $\eqref{E-1}_1$ can be written as
\begin{equation}\label{E-4}
\dot{U}_iU_i^{\dagger}=\Omega_i+\frac{\kappa}{2}R\left(U_i^{\dagger}-U_i\right),\quad i\in[N].
\end{equation}
For \(\alpha\in(\theta_*,\pi)\), let \(\alpha^\infty\) be determined by \eqref{E-2}, and set
\begin{equation}\label{E-5}
T_3^*\coloneqq \frac{\pi}{\kappa J(\alpha)\sin\alpha-\|{\mathfrak V}\|_\infty},\quad
\mu_I\coloneqq \cos\alpha^\infty J(\alpha^\infty)-d\sin\alpha^\infty\operatorname{Lip}_{O(\mathrm{U}(d))_{\alpha^\infty}}(I).
\end{equation}

\subsection{Finite-time entrance into the stable region}\label{subsec:5-1}
In this subsection, we prove a result parallel to that of Section~\ref{subsec:4-1}.

\begin{lemma}\label{L5.1}
Suppose that $({\mathcal F}_{A}1)$--$({\mathcal F}_{A}3)$ hold and
\[
\alpha\in(\theta_*,\pi),\quad
\kappa>\frac{\|{\mathfrak V}\|_\infty}{J(\alpha)\sin\alpha},\quad
U_i^0\in\overline{O(\mathrm{U}(d))_{\alpha}},\quad i\in[N],
\]
and let \({\mathcal U}(t) \) be the global solution to \eqref{E-1}. Then there exists $T<T_3^*$ such that
\[
{\mathcal U}(t)\in O(\mathrm{U}(d))_{\alpha^\infty}^N,\quad t\geq T.
\]
\end{lemma}
\begin{proof}
By $({\mathcal F}_{A}1)$, the vector field is locally Lipschitz on $O(\mathrm{U}(d))_\pi^N$, so the initial data admit a unique maximal solution. We first derive a scalar differential inequality for the maximal spectral angle and then use $({\mathcal F}_{A}3)$ to prove confinement in $\overline{O(\mathrm{U}(d))_\alpha}^{\,N}$ and the entrance-time estimate. Since this closed region is compactly contained in $O(\mathrm{U}(d))_\beta^N$ for every $\beta\in(\alpha,\pi)$, the confinement also gives continuation for all $t\geq0$. \newline

\noindent\textbf{$\bullet$ Step A} (Differential inequality for the maximal spectral angle):
For each $i\in[N]$, we define
\begin{equation}\notag
\rho_i(t)\coloneqq
\cos^{-1}\!\left(\lambda_{\min}\!\left(\frac{U_i(t)+U_i(t)^\dagger}{2}\right)\right),
\quad
\rho(t)\coloneqq\max_{i\in[N]}\rho_i(t).
\end{equation}
The spectral calculation in \eqref{D-12}--\eqref{D-14} applies to \eqref{E-4}, with the lower bound for $R$ modified as follows. By $({\mathcal F}_{A}2)$ and the monotonicity of $J$, we have
\[
R(t)\geq\frac{1}{N}\sum_{j=1}^N J(\rho_j(t))\geq J(\rho(t)).
\]
Consequently, the maximal-angle estimate gives
\begin{equation}\label{E-6}
D^+\rho(t)
\leq \|{\mathfrak V}\|_\infty-\kappa R(t)\sin\rho(t)
\leq \|{\mathfrak V}\|_\infty-\kappa J(\rho(t))\sin\rho(t).
\end{equation}
\vspace{0.2cm}

\noindent\textbf{$\bullet$ Step B} (Entrance into \(O(\mathrm{U}(d))_{\alpha^\infty}^N\)):
Whenever \(\rho(t)\in[\alpha^\infty,\alpha]\), $({\mathcal F}_{A}3)$ and \eqref{E-6} imply
\[
D^+\rho(t)\leq \|{\mathfrak V}\|_\infty-\kappa J(\alpha)\sin\alpha<0.
\]
Hence \(\rho\) cannot cross \(\alpha\) outward and reaches \(\alpha^\infty\) no later than
\[
\frac{\alpha-\alpha^\infty}{\kappa J(\alpha)\sin\alpha-\|{\mathfrak V}\|_\infty}
< \frac{\pi}{\kappa J(\alpha)\sin\alpha-\|{\mathfrak V}\|_\infty}=T_3^*.
\]

\vspace{0.1cm}

\noindent If $\rho(0)<\alpha^\infty$, the same strict boundary inequality shows positive invariance, and we may take $T=0$. Otherwise, let $\tau$ be the first time at which $\rho\leq\alpha^\infty$; the preceding bound gives $\tau<T_3^*$. 

\noindent If $\rho(\tau)=\alpha^\infty$, \eqref{E-2} and \eqref{E-6} give $D^+\rho(\tau)<0$. Hence there is $\varepsilon>0$ such that 
\[ \rho(t)<\alpha^\infty \quad \mbox{for $t\in(\tau,\tau+\varepsilon]$}. \]
Choose $T\in(\tau,\min\{\tau+\varepsilon,T_3^*\})$. The strict boundary inequality prevents a subsequent exit, so $\rho(t)<\alpha^\infty$ for all $t\geq T$, as required.
\end{proof}
\begin{remark}
If \({\mathcal U}^0\in O(\mathrm{U}(d))_{\alpha^\infty}^N\), then one may take \(T=0\).
\end{remark}

\subsection{Exponential stability}\label{subsec:5-2}
Next, we establish exponential stability in the stable region. Throughout this subsection, we assume that the coefficient $\mu_I$ defined in \eqref{E-5} is positive.

\begin{lemma}\label{L5.2}
Suppose that $({\mathcal F}_{A}1)$--$({\mathcal F}_{A}3)$ hold,  and parameters satisfy 
\[ \kappa>0, \quad \mu_I>0, \quad 1\leq p\leq\infty \quad \mbox{and} \quad 1\leq q\leq\infty, \]
and let ${\mathcal U}$ and $\widetilde{\mathcal U}$ be two global solutions to \eqref{E-1} satisfying 
\[
{\mathcal U}(t),~~\widetilde{\mathcal U}(t)\in O(\mathrm{U}(d))_{\alpha^\infty}^N,\quad t\geq0.
\]
Then, we have
\begin{equation}\label{E-7}
\|{\mathcal U}(t)-\widetilde{\mathcal U}(t)\|_{\ell_q^N(S_p^d)}
\leq e^{-\kappa\mu_I t}
\|{\mathcal U}(0)-\widetilde{\mathcal U}(0)\|_{\ell_q^N(S_p^d)},\quad t\geq0.
\end{equation}
\end{lemma}

\begin{proof}
The proof consists of four steps. First, we derive the difference equation and record the stable-region bounds. Second, we estimate the Schatten norm. Third, we control the difference of the order parameters. Finally, we pass to the mixed norm. \newline

\noindent\textbf{$\bullet$ Step A} (Difference equation and stable-region bounds): We set 
\[
W_i\coloneqq U_i-\widetilde U_i,\quad
\mathcal W\coloneqq(W_1,\ldots,W_N).
\]
The condition $\mu_I>0$ implies $\cos\alpha^\infty J(\alpha^\infty)>0$, and hence $\alpha^\infty<\pi/2$. Since both solutions lie in $O(\mathrm{U}(d))_{\alpha^\infty}^N$, we obtain
\[
R,\widetilde R\geq J(\alpha^\infty),\quad
\frac{U_i+U_i^\dagger}{2},\frac{\widetilde U_i+\widetilde U_i^\dagger}{2}
\succeq\cos\alpha^\infty I_d, \quad 
\|I_d-\widetilde U_i^2\|_{\mathrm{op}}\leq2\sin\alpha^\infty.
\]
Using $U_i^2-\widetilde U_i^2=U_iW_i+W_i\widetilde U_i$, we obtain
\[
\dot W_i
=\Omega_iW_i-\frac{\kappa R}{2}\bigl(U_iW_i+W_i\widetilde U_i\bigr)
+\frac{\kappa(R-\widetilde R)}{2}\bigl(I_d-\widetilde U_i^2\bigr).
\]

\vspace{0.5em}

\noindent\textbf{$\bullet$ Step B} (Schatten norm estimate):
Applying the argument for \eqref{D-20} to the difference equation in
Step A and using $R\geq J(\alpha^\infty)$ and $\cos\alpha^\infty>0$,
we obtain, for every $1\leq p\leq\infty$,
\begin{equation}\label{E-8}
	D^+\|W_i\|_{S_p^d}
	\leq-\kappa\cos\alpha^\infty J(\alpha^\infty)\|W_i\|_{S_p^d}
	+\kappa d^{1/p}\sin\alpha^\infty|R-\widetilde R|.
\end{equation}
\vspace{0.5em}

\noindent\textbf{$\bullet$ Step C} (Difference of the order parameters):
By \eqref{E-3} and $({\mathcal F}_{A}1)$,
\begin{equation} \label{E-9}
|R-\widetilde R| \leq\operatorname{Lip}_{O(\mathrm{U}(d))_{\alpha^\infty}}(I)
\|\mathcal W\|_{\ell_1^N(S_1^d)} \leq d^{1-1/p}\operatorname{Lip}_{O(\mathrm{U}(d))_{\alpha^\infty}}(I)
\|\mathcal W\|_{\ell_1^N(S_p^d)}.
\end{equation}
We combine \eqref{E-8} and \eqref{E-9} to obtain
\[
D^+\|W_i\|_{S_p^d}
\leq-\kappa\cos\alpha^\infty J(\alpha^\infty)\|W_i\|_{S_p^d}
+\kappa d\sin\alpha^\infty
\operatorname{Lip}_{O(\mathrm{U}(d))_{\alpha^\infty}}(I)
\|\mathcal W\|_{\ell_1^N(S_p^d)}.
\]
\vspace{0.5em}
\noindent\textbf{$\bullet$ Step D} (Mixed-norm estimate):
For $1\leq q\leq\infty$, we set
\[
Y_q(t)\coloneqq\|\mathcal W(t)\|_{\ell_q^N(S_p^d)}.
\]
Fix $1\leq p\leq\infty$. For $q=1$, averaging the component
estimate in Step C and using \eqref{E-5} gives
\[
D^+Y_1(t)
\leq \frac1N\sum_{i=1}^N D^+\|W_i(t)\|_{S_p^d}
\leq -\kappa\mu_I Y_1(t).
\]
For $1<q<\infty$ and $Y_q(t)>0$, the chain rule, the component
estimate in Step C, H\"older's inequality, and $Y_1(t)\leq Y_q(t)$ give
\begin{align*}
	D^+Y_q
	&\leq-\kappa\cos\alpha^\infty J(\alpha^\infty)Y_q +\kappa d\sin\alpha^\infty
	\operatorname{Lip}_{O(\mathrm{U}(d))_{\alpha^\infty}}(I)Y_q =-\kappa\mu_IY_q.
\end{align*}
If $Y_q(t)=0$, uniqueness gives $\mathcal W(s)=0$ for all $s\geq t$, and hence $D^+Y_q(t)=0$. Thus, for $q<\infty$,
\[
D^+Y_q(t)\leq-\kappa\mu_IY_q(t),\quad t\geq0.
\]
For $q=\infty$, apply the maximum formula from Step D of the proof
of Lemma~\ref{L4.5} to the component estimate in Step C above.
Since $Y_1(t)\leq Y_\infty(t)$, we obtain
\[
D^+Y_\infty(t)\leq-\kappa\mu_IY_\infty(t).
\]
Gronwall's lemma now gives \eqref{E-7} for all $1\leq p,q\leq\infty$.
\end{proof}

\subsection{Relaxation dynamics}\label{subsec:5-3} 
In this subsection, we show the existence and uniqueness of an equilibrium, and establish convergence of every solution satisfying the stated hypotheses.
\begin{proposition}\label{P5.1}
Suppose that $({\mathcal F}_{A}1)$--$({\mathcal F}_{A}3)$ hold and system parameters satisfy
\[
\alpha\in(\theta_*,\pi),\quad
\kappa>\frac{\|{\mathfrak V}\|_\infty}{J(\alpha)\sin\alpha},\quad
\mu_I>0.
\]
Then system \eqref{E-1} admits a unique equilibrium ${\mathcal U}^\infty\in O(\mathrm{U}(d))_{\alpha^\infty}^N$. 
Moreover, this is the only equilibrium in \(O(\mathrm{U}(d))_{\alpha}^N\).
\end{proposition}
\begin{proof} We split the proof into two steps. 

\vspace{0.1cm}

\noindent $\bullet$~Step A (Existence in the stable region):~We fix initial data in \(O(\mathrm{U}(d))_{\alpha^\infty}^N\), and let \({\mathcal U}(t)\) be the corresponding solution. Then, Lemma~\ref{L5.1} with $T=0$ shows that this region is positively invariant. For \(s\geq0\), the time translate \(\widetilde{\mathcal U}(t)\coloneqq {\mathcal U}(t+s)\) is another solution in the same stable region. Lemma~\ref{L5.2} therefore gives, with \(p=q=2\),
\begin{equation}\label{E-10}
\|{\mathcal U}(t+s)-{\mathcal U}(t)\|_{\ell_2^N(S_2^d)}
\leq e^{-\kappa\mu_I t}
\|{\mathcal U}(s)-{\mathcal U}(0)\|_{\ell_2^N(S_2^d)}.
\end{equation}
Because $\mathrm{U}(d)^N$ is compact, the last factor in \eqref{E-10} is bounded uniformly in $s$. For $t_2\geq t_1\geq0$, we take 
\[ t=t_1 \quad \mbox{and} \quad s=t_2-t_1 \quad \mbox{in \eqref{E-10}.}  \]
Then, we have
\[
\|\mathcal U(t_2)-\mathcal U(t_1)\|_{\ell_2^N(S_2^d)}
\leq Ce^{-\kappa\mu_I t_1},
\]
where $C$ is independent of $t_1$ and $t_2$. Thus $\mathcal U(t)$ is Cauchy and it converges to some $\mathcal U^\infty\in\mathrm{U}(d)^N$. Positive invariance and continuity yield
\[
\mathcal U^\infty\in\overline{O(\mathrm U(d))_{\alpha^\infty}}^{\,N}
\subset O(\mathrm U(d))_\alpha^N.
\]
Hence the flow is defined at the limiting state by Lemma~\ref{L5.1}.
Let $\Phi_s$ denote the flow of \eqref{E-1}. By continuity of this flow, we have
\[
\Phi_s({\mathcal U}^\infty)
=\lim_{t\to\infty}\Phi_s({\mathcal U}(t))
=\lim_{t\to\infty}{\mathcal U}(t+s)
={\mathcal U}^\infty,
\quad s\geq0,
\]
so \({\mathcal U}^\infty\) is an equilibrium. The strict inward estimate in Lemma~\ref{L5.1} excludes the boundary of \(O(\mathrm{U}(d))_{\alpha^\infty}^N\), and hence the equilibrium belongs to the stated open set. \newline

\noindent $\bullet$~Step B (Uniqueness of equilibria in \(O(\mathrm{U}(d))_{\alpha}^N\)):~If \({\mathcal V}^\infty\) is another equilibrium in \(O(\mathrm{U}(d))_{\alpha}^N\), Lemma~\ref{L5.1} places the constant trajectory \({\mathcal V}^\infty\) in \(O(\mathrm{U}(d))_{\alpha^\infty}^N\). Applying Lemma~\ref{L5.2} to two constant trajectories gives
\[
\|{\mathcal U}^\infty-{\mathcal V}^\infty\|_{\ell_2^N(S_2^d)}
\leq e^{-\kappa\mu_I t}
\|{\mathcal U}^\infty-{\mathcal V}^\infty\|_{\ell_2^N(S_2^d)},
\]
for every \(t\geq0\). Letting \(t\to\infty\), we obtain uniqueness.
\end{proof}
Next, we establish exponential convergence to the equilibrium.
\begin{theorem}\label{T5.1}
Suppose that $({\mathcal F}_{A}1)$--$({\mathcal F}_{A}3)$ hold, and parameters satisfy 
\[ 1\leq p\leq\infty,\quad 1\leq q\leq\infty,~~
\alpha\in(\theta_*,\pi),~~
\kappa>\frac{\|{\mathfrak V}\|_\infty}{J(\alpha)\sin\alpha},~~
\mu_I>0,\quad U_i^0\in\overline{O(\mathrm{U}(d))_{\alpha}},\quad i\in[N],
\]
and let \({\mathcal U}(t)\) be the global solution to \eqref{E-1}. Then the following assertions hold.
\begin{enumerate}
\item
System \eqref{E-1} admits a unique equilibrium $\mathcal U^\infty\in O(\mathrm{U}(d))_{\alpha^\infty}^N
$, and it is the only equilibrium in $O(\mathrm{U}(d))_\alpha^N$, as established in Proposition~\ref{P5.1}.
\vspace{0.1cm}
\item
There exists \(T<T_3^*\) such that
\begin{equation}\label{E-11}
\|{\mathcal U}(t)-{\mathcal U}^\infty\|_{\ell_q^N(S_p^d)}
\leq e^{-\kappa\mu_I(t-T)}
\|{\mathcal U}(T)-{\mathcal U}^\infty\|_{\ell_q^N(S_p^d)},
\quad t\geq T.
\end{equation}
If~\({\mathcal U}^0\in O(\mathrm{U}(d))_{\alpha^\infty}^N\), then one may take \(T=0\).
\end{enumerate}
\end{theorem}
\begin{proof}
Lemma~\ref{L5.1} provides the entrance time \(T\). The shifted solution \({\mathcal U}(T+\cdot)\) and the constant solution \({\mathcal U}^\infty\) both lie in \(O(\mathrm{U}(d))_{\alpha^\infty}^N\). Finally, we apply Lemma~\ref{L5.2} to these two solutions to find \eqref{E-11}.
\end{proof}
\begin{corollary}\label{C:general-unitary-basin}
	Suppose that $({\mathcal F}_{A}1)$--$({\mathcal F}_{A}3)$ hold.
	Fix $\mathfrak V$ and initial data with
	$-1\notin\operatorname{spec}(U_i^0)$ for each $i\in[N]$.
	There exists $\kappa_0\geq0$ such that, for every $\kappa>\kappa_0$,
	the solution of \eqref{E-1} exists globally and converges
	exponentially to an equilibrium.
\end{corollary}
\begin{proof}
	Write $F(\theta)=J(\theta)\sin\theta$. By $({\mathcal F}_{A}3)$,
	$F$ is positive on $(0,\pi)$, so $J(0)>0$ by monotonicity of $J$.
	As $\alpha\uparrow\pi$, \eqref{E-2} gives
	$F(\alpha^\infty)=F(\alpha)\to0$. Since $F$ is strictly increasing
	on $(0,\theta_*]$, we have $a=\alpha^\infty\to0$.
	Fix $\beta\in(0,\pi)$ and let $L_\beta$ be the Lipschitz constant
	of $I$ on $O(\mathrm U(d))_\beta$ from $({\mathcal F}_{A}1)$.
	For $a<\beta$, \eqref{E-5} gives
	\[
	\mu_I\geq\cos(a)J(a)-d\sin(a)L_\beta\longrightarrow J(0)>0.
	\]
	Choose $\alpha\in(\theta_*,\pi)$ close enough to $\pi$ that
	$\mu_I>0$ and $U_i^0\in O(\mathrm U(d))_\alpha$ for all $i$.
	Theorem~\ref{T5.1} now applies with
	$\kappa_0=\|\mathfrak V\|_\infty/[J(\alpha)\sin\alpha]$.
\end{proof}
\section{\texorpdfstring{Winfree--Lohe models as Kuramoto--Lohe subsystems}{Winfree-Lohe models as Kuramoto-Lohe subsystems}}\label{sec:6}
\setcounter{equation}{0}
There is another route to Winfree-type models on Lie groups: one can realize them as the invariant subsystems of a larger Kuramoto--Lohe ensemble. We begin with the circle construction of Moreno-Earle, Ryoo, and To \cite[Lemma 3.3]{morenoearle2026relaxation}, and then give a symmetry construction within the Lie-group framework of \cite{ryoo2025asymptotic}. The resulting examples include the unitary model of Section~\ref{sec:4} and a model on the special unitary group.

\subsection{The circle embedding}\label{subsec:lie-circle}
For a solution $\Theta=(\theta_1,\ldots,\theta_N)$ of \eqref{A-1}, we define an ensemble of $4N$ phases and intrinsic frequencies by
\begin{equation}\label{LE-circle-embedding}
\Phi=(\Theta,-\Theta,0^{2N}),\quad
\widetilde{\mathcal V}=(\mathcal V,-\mathcal V,0^{2N}),
\end{equation}
where $0^{2N}$ denotes $2N$ zero entries. Then $\Phi$ solves the Kuramoto model
\begin{equation*}\label{LE-circle-kuramoto}
\dot\varphi_a=\widetilde\nu_a+
\frac{2\kappa}{4N}\sum_{b=1}^{4N}\sin(\varphi_b-\varphi_a),
\quad a\in[4N].
\end{equation*}
Indeed, the interaction in the equation for the phase $\theta_i$ is given as 
\[ \frac{\kappa}{2N}\left[
\sum_{j=1}^N\bigl(\sin(\theta_j-\theta_i)+\sin(-\theta_j-\theta_i)\bigr)
-2N\sin\theta_i\right] =-\frac{\kappa}{N}
\sum_{j=1}^N(1+\cos\theta_j)\sin\theta_i.
\]
Oddness gives the reflected equations and makes the velocities of all $2N$ zero phases vanish. Thus, the state \eqref{LE-circle-embedding} describes an invariant subsystem, and the first $N$ coordinates recover the Winfree dynamics. The ambient Kuramoto coupling strength is $2\kappa$. The same construction in \cite{morenoearle2026relaxation} also applies to the inertial equations. We only employ its first-order version here.

\subsection{A symmetry construction on Lie groups}\label{subsec:lie-symmetry}
Let $G$ be a Lie group with identity $e$ and Lie algebra $\mathfrak g$. For $x\in G$, we write $R_x(y)=yx$. Then, the generalized Kuramoto--Lohe model of \cite{ryoo2025asymptotic} is given as follows.
\begin{equation}\label{LE-ambient}
\dot Z_a=(dR_{Z_a})_e\left(
B_a+\frac{K}{M}\sum_{b=1}^M\phi(Z_bZ_a^{-1})\right),
\quad a\in[M],
\end{equation}
where $B_a\in\mathfrak g$ and $\phi:G\to\mathfrak g$ is $C^1$. The Cauchy--Lipschitz theory only gives local existence and uniqueness; under the local attraction assumptions $\phi(e)=0$ and the positivity of the real parts of the eigenvalues of $(d\phi)_e$, global existence holds if the initial states are sufficiently close to one another and the coupling strength $K$ is sufficiently large relative to $\max_{a\in[M]}|B_a|$ \cite[Theorem~1.2(2)]{ryoo2025asymptotic}.

On a nonabelian group, the scalar reflection cannot simply be replaced by inversion. For example, the drift $\dot X=HX$ of a matrix group gives
\[
\frac{d}{dt}X^{-1}=-X^{-1}H,
\]
whereas a reflected oscillator with intrinsic frequency $-H$ has drift $-HX^{-1}$. These need not agree. We therefore replace the two scalar symmetries by an orbit under genuine group automorphisms.

Let a finite group $\Gamma$ of order $m$ act on $G$ by automorphisms; the action need not be faithful. We write $\sigma(x)$ for the image of $x\in G$ under $\sigma\in\Gamma$, and $A_\sigma=(d\sigma)_e$ for the differential of this automorphism. We assume
\begin{equation}\label{LE-symmetry-assumptions}
\phi(\sigma(x))=A_\sigma\phi(x),\quad
\frac1m\sum_{\sigma\in\Gamma}A_\sigma=0,
\quad \sigma\in\Gamma,\quad x\in G.
\end{equation}
The second condition implies in particular that the induced action $A_\sigma$ on $\mathfrak g$ has no nonzero common fixed vector. Define the averaged interaction:
\begin{equation*}\label{LE-averaged-kernel}
\mathcal K_\Gamma(x,y)\coloneqq
\frac1m\sum_{\sigma\in\Gamma}\phi(\sigma(y)x^{-1}).
\end{equation*}
We denote the identity element of $\Gamma$ by $\mathrm{id}$.

\begin{proposition}\label{P:lie-subsystem}
Suppose that $\phi(e)=0$ and \eqref{LE-symmetry-assumptions} hold. Fix $H_1,\ldots,H_N\in\mathfrak g$ and initial data $x_1^0,\ldots,x_N^0\in G$. We index $mN$ oscillators by $(i,\sigma)\in[N]\times\Gamma$ and the remaining $mN$ oscillators by $(\mathrm{anc},\ell)$, $\ell\in[mN]$. Their intrinsic frequencies and initial states are prescribed by
\[
\begin{aligned}
B_{i,\sigma}&=A_\sigma H_i,
& Z_{i,\sigma}(0)&=\sigma(x_i^0),
&& i\in[N],\quad \sigma\in\Gamma,\\
B_{\mathrm{anc},\ell}&=0,
& Z_{\mathrm{anc},\ell}(0)&=e,
&& \ell\in[mN].
\end{aligned}
\]
Each pair $(i,\sigma)$ labels a separate oscillator, even when the action of $\Gamma$ is not faithful. Let $Z_{i, \sigma}(t)$ and $Z_{\mathrm{anc}, \ell}(t)$ denote the solution to \eqref{LE-ambient} with $M=2mN$ and $K=2\kappa$.

Then, for $x_i(t) \coloneqq Z_{i, \mathrm{id}}(t)$, the following assertions hold:
\begin{enumerate}
\item The solution preserves the relations
\begin{equation}\label{LE-orbit-embedding}
\begin{aligned}
Z_{i,\sigma}(t)&=\sigma(x_i(t)), && i\in[N],\quad \sigma\in\Gamma,\\
Z_{\mathrm{anc},\ell}(t)&=e, && \ell\in[mN].
\end{aligned}
\end{equation}
\vspace{0.1cm}
\item The variables $x_1,\ldots,x_N$ solve the reduced Cauchy problem
\begin{equation}\label{LE-reduced}
\begin{cases}
\displaystyle \dot x_i=(dR_{x_i})_e\left[H_i +\kappa\phi(x_i^{-1})+\frac{\kappa}{N}\sum_{j=1}^N \mathcal K_\Gamma(x_i,x_j)\right],\\[2mm]
x_i(0)=x_i^0,\quad i\in[N].
\end{cases}
\end{equation}
\end{enumerate} 
\end{proposition}

\begin{proof}
Since the vector field in \eqref{LE-reduced} is $C^1$, the reduced Cauchy problem \eqref{LE-reduced} has a unique maximal solution $X=(x_1,\ldots,x_N)$ on an interval $[0,T)$. We define the curves $Z_{i,\sigma}$ and $Z_{\mathrm{anc},\ell}$ by \eqref{LE-orbit-embedding} and verify that they solve the ambient equation \eqref{LE-ambient} with the prescribed initial data.

\medskip
\noindent\textbf{$\bullet$ Step A}
(Ancillary oscillators):
For an ancillary oscillator, the intrinsic frequency is zero. At $Z_{\mathrm{anc},\ell}=e$, its interaction term is
\[
\frac{\kappa}{mN}\left[\sum_{j=1}^N\sum_{\sigma\in\Gamma}\phi(\sigma(x_j))+mN\phi(e)\right]=\frac{\kappa}{mN}\sum_{j=1}^N \left(\sum_{\sigma\in\Gamma}A_\sigma\right)\phi(x_j)=0.
\]
Here, we used $\phi(e)=0$ and both identities in \eqref{LE-symmetry-assumptions}. Thus the constant curves $Z_{\mathrm{anc},\ell}(t)=e$ satisfy \eqref{LE-ambient}.

\medskip
\noindent\textbf{$\bullet$ Step B}
(Dynamics of $Z_{i, \mathrm{id}}$):
For the oscillator indexed by $(i,\mathrm{id})$, we have $Z_{i,\mathrm{id}}=x_i$ and $B_{i,\mathrm{id}}=H_i$. The contributions from the orbit oscillators and the identity oscillators are, respectively,
\[
\frac{\kappa}{mN}\sum_{j=1}^N\sum_{\sigma\in\Gamma}\phi(\sigma(x_j)x_i^{-1})=\frac{\kappa}{N}\sum_{j=1}^N\mathcal K_\Gamma(x_i,x_j),
\quad \frac{\kappa}{mN}\sum_{\ell=1}^{mN}\phi(ex_i^{-1})=\kappa\phi(x_i^{-1}).
\]
Substituting these expressions into \eqref{LE-ambient} gives \eqref{LE-reduced}. Hence the curves $Z_{i,\mathrm{id}}$ satisfy their ambient equations.

\medskip
\noindent\textbf{$\bullet$ Step C}
(Dynamics of $Z_{i, \tau}$):
Fix $\tau\in\Gamma$.
Since $\tau$ is a group automorphism, we have
\[
\tau\circ R_x=R_{\tau(x)}\circ\tau,
\quad
\tau(x^{-1})=\tau(x)^{-1}.
\]
Differentiating the first identity at $e$ gives
\[
(d\tau)_x\circ(dR_x)_e
=
(dR_{\tau(x)})_e\circ A_\tau.
\]
We next compare the interaction terms at $x_i$ and $\tau(x_i)$.
For each $j\in[N]$, the change of index $\sigma=\tau\rho$
and \eqref{LE-symmetry-assumptions} yield
\[
\sum_{\sigma\in\Gamma}\phi\bigl(\sigma(x_j)\tau(x_i)^{-1}\bigr)=\sum_{\rho\in\Gamma}\phi\bigl(\tau(\rho(x_j)x_i^{-1})\bigr)=A_\tau\sum_{\rho\in\Gamma}\phi\bigl(\rho(x_j)x_i^{-1}\bigr).
\]
The contribution from the identity oscillators transforms in
the same way:
\[
\phi\bigl(\tau(x_i)^{-1}\bigr)
=
\phi\bigl(\tau(x_i^{-1})\bigr)
=
A_\tau\phi(x_i^{-1}).
\]
Combining these identities with \eqref{LE-reduced}, we obtain
\[
\dot Z_{i,\tau}=(d\tau)_{x_i}\dot x_i=(dR_{Z_{i,\tau}})_e \left(A_\tau H_i+\frac{K}{M}\sum_{b=1}^{M}\phi(Z_bZ_{i,\tau}^{-1})\right).
\]
Here, the sum over $b$ includes
both the orbit oscillators and the identity oscillators.
Thus every curve $Z_{i,\tau}$ satisfies
\eqref{LE-ambient}.

The constructed curves \eqref{LE-orbit-embedding} therefore solve \eqref{LE-ambient} with the prescribed initial data. By uniqueness, they coincide with the ambient solution
on $[0,T)$. It remains to compare the maximal intervals of existence. 

The construction already gives an ambient solution on the entire interval $[0,T)$. If $T<\infty$ and the ambient solution extended beyond $T$, then $x_i(t)=Z_{i,\mathrm{id}}(t)$ would have a limit in $G$ as $t\uparrow T$ for every $i$. Local existence for \eqref{LE-reduced} at this limiting configuration would extend $X$ beyond $T$, contradicting its maximality. Hence the maximal intervals agree, and \eqref{LE-orbit-embedding} holds throughout this interval.
\end{proof}

Using the embedding in Proposition~\ref{P:lie-subsystem}, we apply
\cite[Theorem~1.2]{ryoo2025asymptotic} to the ambient system
\eqref{LE-ambient} to obtain the following relaxation result for
the reduced system \eqref{LE-reduced}.\footnote{The theorem number
	refers to \href{https://arxiv.org/abs/2109.14761v3}{arXiv:2109.14761v3}.}

Fix an inner product on $\mathfrak g$, let $|\cdot|$ denote its norm, and let $d_G$ denote the induced right-invariant Riemannian distance on $G$. We denote
\[
\mathfrak H=(H_1,\ldots,H_N),\quad X=(x_1,\ldots,x_N),\quad \|\mathfrak H\|_\Gamma \coloneqq \max_{i\in[N],\,\sigma\in\Gamma}|A_\sigma H_i|,
\]
and set
\[
B_G(e,r)\coloneqq\{x\in G:d_G(x,e)<r\}.
\]

\begin{theorem}[Relaxation of the Winfree--Lohe model \eqref{LE-reduced}]\label{T:lie-relaxation}
Suppose that $\phi(e)=0$ and \eqref{LE-symmetry-assumptions} hold, and suppose that every eigenvalue of $(d\phi)_e$ has positive real part. Then, there exist constants
\[
0<r_1<r_2,\quad C>0,\quad c>0,
\]
depending only on the Lie-group structure, $\phi$, the metric, and the fixed action of $\Gamma$, but not on $N$, such that the following assertions hold.

If
\[
\kappa>C\|\mathfrak H\|_\Gamma,\quad \max_{i\in[N]}d_G(x_i^0,e)<r_1,
\]
then:
\begin{enumerate}
\item
The solution of \eqref{LE-reduced} exists uniquely for all $t\ge0$ and satisfies
\[
\max_{i\in[N]}d_G(x_i(t),e)<r_2, \quad t\ge0.
\]

\item
There is a unique equilibrium $X^\infty$ of \eqref{LE-reduced} in $B_G(e,r_1)^N$, and
\[
\max_{i\in[N]}d_G(x_i^\infty,e) \le C\frac{\|\mathfrak H\|_\Gamma}{\kappa}.
\]

\item
Every solution with initial data in $B_G(e,r_1)^N$ converges exponentially to the same equilibrium $X^\infty$, and its speeds converge exponentially to zero. More precisely,
\[
\max_{i\in[N]}d_G(x_i(t),x_i^\infty) +\kappa^{-1}\max_{i\in[N]} \left| (dR_{x_i(t)^{-1}})_{x_i(t)}\dot x_i(t) \right| \le Ce^{-c\kappa t}, \quad t\ge0.
\]
\end{enumerate}
\end{theorem}

\begin{proof}
It suffices to prove the theorem for an inner product adapted to
$A=(d\phi)_e$. Indeed, since every eigenvalue of $A$ has positive
real part, there is an inner product $\langle\cdot,\cdot\rangle_*$
such that $\langle Av,v\rangle_*\geq\lambda |v|_*^2$ for some
$\lambda>0$ and all $v\in\mathfrak g$
\cite[Section~2, (2.1)]{ryoo2025asymptotic}.
After rescaling this inner product, we may assume that
$|v|\leq |v|_*\leq b|v|$ for some $b\geq1$.
Comparison of curve lengths gives $d_G\leq d_*\leq b d_G$ for the
corresponding right-invariant distances.

Suppose the theorem holds for $d_*$ with constants
$r_{1,*},r_{2,*},C_*,c_*$. For $d_G$, we take
\[
r_1=r_{1,*}/b,\qquad r_2=r_{2,*},\qquad c=c_*,
\qquad C=bC_*\max\{1,r_1^{-1}\}.
\]
Then $B_G(e,r_1)\subset B_*(e,r_{1,*})$ and
$\|\mathfrak H\|_{\Gamma,*}\leq b\|\mathfrak H\|_\Gamma$.
Thus the adapted-metric theorem applies, and its equilibrium satisfies
\[
\max_i d_G(x_i^\infty,e)
\leq bC_*\|\mathfrak H\|_\Gamma/\kappa < r_1.
\]
Uniqueness in $B_G(e,r_1)^N$ follows from the ball inclusion above.
The trajectory and speed bounds follow from $d_G\leq d_*$ and
$|\cdot|\leq|\cdot|_*$. We may therefore assume for the rest of
the proof that the chosen inner product is adapted to $A$.
	
We first verify the hypotheses of \cite[Theorem 1.2]{ryoo2025asymptotic} for the ambient system. We then use the stationary identity oscillators to obtain the corresponding conclusions for the reduced system.

\medskip
\noindent $\bullet$ Step A (Initial concentration of the ambient system): We set
\[
L_\Gamma \coloneqq \max_{\sigma\in\Gamma}\|A_\sigma\|_{\mathrm{op}}.
\]
Since $\mathrm{id}\in\Gamma$, we have $L_\Gamma\ge1$. The identity
\[
(d\sigma)_x\circ(dR_x)_e = (dR_{\sigma(x)})_e\circ A_\sigma
\]
and right invariance of the metric imply, by comparison of curve lengths, that
\[
d_G(\sigma(x),\sigma(y)) \le L_\Gamma d_G(x,y).
\]
In particular, since $\sigma(e)=e$,
\[
d_G(\sigma(x),e)\le L_\Gamma d_G(x,e).
\]

We fix the constants $C_{\mathrm R},c_{1,\mathrm R},c_{\mathrm R},c_{2,\mathrm R}>0$ appearing in \cite[Theorem 1.2]{ryoo2025asymptotic}, applied with the ambient coupling $K$. Here $c_{1,\mathrm R}$ is the threshold for the initial and propagated pairwise distances, $c_{\mathrm R}$ is the coefficient in the coupling condition $K>c_{\mathrm R}\max_a|B_a|$, $C_{\mathrm R}$ is the coefficient in the asymptotic distance bound $C_{\mathrm R}\max_a|B_a|/K$, and $c_{2,\mathrm R}$ is the exponential rate in part (3). We also denote by $D_{\mathrm R}>0$ a constant implicit in the $O_G(e^{-c_{2,\mathrm R}Kt})$ estimate in that part. These constants depend only on the ambient Lie-group structure, $\phi$, and the chosen metric.

We choose
\[
r_2=c_{1,\mathrm R},\quad r_1=\frac{c_{1,\mathrm R}}{4L_\Gamma},
\]
and take the constant $C$ in the statement sufficiently large that
\[
C\ge\max\left\{c_{\mathrm R},\frac{C_{\mathrm R}}{r_1}\right\}.
\]
We enlarge $C$ further at the end of the proof.

In \eqref{LE-ambient}, take
\[
M=2mN,\quad K=2\kappa,
\]
and prescribe
\[
\begin{aligned}
Z_{i,\sigma}(0)&=\sigma(x_i^0), & B_{i,\sigma}&=A_\sigma H_i, && i\in[N],\quad \sigma\in\Gamma,\\
Z_{\mathrm{anc},\ell}(0)&=e, & B_{\mathrm{anc},\ell}&=0, && \ell\in[mN].
\end{aligned}
\]
This choice and the assumption on $\kappa$ imply
\[
\max_{a\in[M]}|B_a|=\|\mathfrak H\|_\Gamma,\quad K=2\kappa>c_{\mathrm R}\max_{a\in[M]}|B_a|.
\]
Moreover, the initial concentration condition gives
\[
d_G(Z_{i,\sigma}(0),Z_{j,\tau}(0)) \le d_G(\sigma(x_i^0),e)+d_G(\tau(x_j^0),e)\le L_\Gamma\bigl(d_G(x_i^0,e)+d_G(x_j^0,e)\bigr) <2L_\Gamma r_1,
\]
and
\[
d_G(Z_{i,\sigma}(0),Z_{\mathrm{anc},\ell}(0)) \le L_\Gamma d_G(x_i^0,e)<L_\Gamma r_1.
\]
The distances between ancillary oscillators are zero. Consequently,
\[
\max_{a,b\in[M]}d_G(Z_a(0),Z_b(0)) <2L_\Gamma r_1=\frac{c_{1,\mathrm R}}{2}<c_{1,\mathrm R}.
\]
Thus the ambient initial data and parameters satisfy all hypotheses of \cite[Theorem 1.2]{ryoo2025asymptotic}.

\medskip
\noindent $\bullet$ Step B (Global existence and boundedness of the reduced solution): By \cite[Theorem 1.2(2)]{ryoo2025asymptotic}, the ambient Cauchy problem has a unique global solution $Z$ satisfying
\[
\max_{a,b\in[M]}d_G(Z_a(t),Z_b(t))<c_{1,\mathrm R}, \quad t\ge0.
\]

Let $X$ be the unique maximal solution of \eqref{LE-reduced} with initial data $X^0$, defined on $[0,T_{\max})$. By Proposition~\ref{P:lie-subsystem}, its embedding \eqref{LE-orbit-embedding} solves \eqref{LE-ambient} with the initial data prescribed in Step A. Uniqueness of the ambient solution therefore gives
\[
Z_{i,\sigma}(t)=\sigma(x_i(t)),\quad Z_{\mathrm{anc},\ell}(t)=e, \quad 0\le t<T_{\max}.
\]
If $T_{\max}<\infty$, then the global ambient solution gives
\[
\lim_{t\uparrow T_{\max}}x_i(t) = Z_{i,\mathrm{id}}(T_{\max})\in G, \quad i\in[N].
\]
Since the reduced vector field is $C^1$, local existence at this limiting configuration extends $X$ beyond $T_{\max}$, contradicting maximality. Hence $T_{\max}=\infty$, and the embedding relations hold for all $t\ge0$. Fix the ancillary index $a_0=(\mathrm{anc},1)$. Using $Z_{a_0}(t)=e$ and $Z_{i,\mathrm{id}}(t)=x_i(t)$, we obtain
\[
\max_{i\in[N]}d_G(x_i(t),e) = \max_{i\in[N]} d_G(Z_{i,\mathrm{id}}(t),Z_{a_0}(t)) <c_{1,\mathrm R}=r_2, \quad t\ge0.
\]
This proves assertion (1).

\medskip
\noindent $\bullet$ Step C (Convergence to a unique equilibrium): By \cite[Theorem 1.2(3)]{ryoo2025asymptotic}, there exist phase representatives $W_a\in G$, $a\in[M]$, of a phase-locked solution of the ambient system such that
\[
\max_{a,b\in[M]} d_G\bigl( Z_a(t)Z_b(t)^{-1},W_aW_b^{-1} \bigr) \le D_{\mathrm R}e^{-c_{2,\mathrm R}Kt}, \quad t\ge0.
\]
We set
\[
x_i^\infty \coloneqq W_{i,\mathrm{id}}W_{a_0}^{-1}, \quad i\in[N].
\]
Since
\[
Z_{i,\mathrm{id}}(t)Z_{a_0}(t)^{-1}=x_i(t),
\]
the preceding estimate gives
\[
\max_{i\in[N]}d_G(x_i(t),x_i^\infty) \le D_{\mathrm R}e^{-2c_{2,\mathrm R}\kappa t}, \quad t\ge0.
\]
Thus the convergence of the ambient relative phases gives convergence of the reduced states themselves.

We next locate the limiting configuration. The distance estimate in \cite[Theorem 1.2(2)]{ryoo2025asymptotic} yields
\[
d_G(x_i^\infty,e) = \lim_{t\to\infty} d_G(Z_{i,\mathrm{id}}(t),Z_{a_0}(t))\le C_{\mathrm R}\frac{\max_{a\in[M]}|B_a|}{K} = \frac{C_{\mathrm R}\|\mathfrak H\|_\Gamma}{2\kappa} <\frac{r_1}{2}.
\]
The last inequality follows from the choice of $C$; when $\|\mathfrak H\|_\Gamma=0$, the preceding upper bound is zero. Hence
\[
X^\infty\in B_G(e,r_1)^N.
\]

To verify that $X^\infty$ is an equilibrium of \eqref{LE-reduced}, let $\Phi_s$ denote its local flow. For sufficiently small $s\ge0$, continuity of the flow gives
\[
\Phi_s(X^\infty) = \lim_{t\to\infty}\Phi_s(X(t)) = \lim_{t\to\infty}X(t+s) = X^\infty.
\]
Therefore $X^\infty$ is an equilibrium.

It remains to prove its uniqueness in $B_G(e,r_1)^N$. Let $Y=(y_1,\ldots,y_N)$ be another equilibrium in this set. By Proposition~\ref{P:lie-subsystem}, the arrays
\[
\begin{aligned}
Z_{i,\sigma}^\infty&=\sigma(x_i^\infty), & Z_{\mathrm{anc},\ell}^\infty&=e,\\
\widetilde Z_{i,\sigma}^\infty&=\sigma(y_i), & \widetilde Z_{\mathrm{anc},\ell}^\infty&=e
\end{aligned}
\]
are stationary solutions of the same ambient equation with the same intrinsic frequencies $B_a$. The estimates in Step A show that both arrays have diameter smaller than $2L_\Gamma r_1<c_{1,\mathrm R}$. Thus \cite[Theorem 1.2(1)]{ryoo2025asymptotic} provides $g\in G$ such that
\[
\widetilde Z_a^\infty=Z_a^\infty g, \quad a\in[M].
\]
At the ancillary coordinate $a_0$, this becomes $e=eg$, so $g=e$. Taking $a=(i,\mathrm{id})$ gives $y_i=x_i^\infty$ for every $i\in[N]$. This proves assertion (2). It also shows that, for fixed $\mathfrak H$ and $\kappa$, all initial data in $B_G(e,r_1)^N$ lead to the same limit.

\medskip
\noindent $\bullet$ Step D (Exponential decay of speeds): We define the right-trivialized reduced vector field by
\[
V_i(X) \coloneqq H_i+\kappa\phi(x_i^{-1}) +\frac{\kappa}{N}\sum_{j=1}^N \mathcal K_\Gamma(x_i,x_j).
\]
Then
\[
(dR_{x_i(t)^{-1}})_{x_i(t)}\dot x_i(t)=V_i(X(t)), \quad V_i(X^\infty)=0.
\]

Since $\phi$ is $C^1$ and $\Gamma$ is finite, the maps $x\mapsto\phi(x^{-1})$ and $(x,y)\mapsto\mathcal K_\Gamma(x,y)$ are $C^1$. Since a right-invariant Riemannian metric on a Lie group is complete, closed metric balls are compact. Any minimizing geodesic joining two points of $B_G(e,r_2)$ has length less than $2r_2$ and therefore lies in $B_G(e,3r_2)$. The derivatives of $x\mapsto\phi(x^{-1})$ and $(x,y)\mapsto\mathcal K_\Gamma(x,y)$ are bounded on the compact sets $\overline{B_G(e,3r_2)}$ and $\overline{B_G(e,3r_2)}^{\,2}$, respectively. Integrating along such geodesics, one variable at a time for $\mathcal K_\Gamma$, gives a constant $L_0>0$ such that
\[
|\phi(x^{-1})-\phi(y^{-1})| \le L_0d_G(x,y), \quad |\mathcal K_\Gamma(x,z)-\mathcal K_\Gamma(y,w)| \le L_0\bigl(d_G(x,y)+d_G(z,w)\bigr)
\]
for $x,y,z,w\in B_G(e,r_2)$. Consequently, for $X,Y\in B_G(e,r_2)^N$,
\[
|V_i(X)-V_i(Y)| \le \kappa L_0d_G(x_i,y_i) +\frac{\kappa L_0}{N}\sum_{j=1}^N \bigl(d_G(x_i,y_i)+d_G(x_j,y_j)\bigr)\le 3\kappa L_0\max_{j\in[N]}d_G(x_j,y_j).
\]

By Step B, $X(t)\in B_G(e,r_2)^N$ for all $t\ge0$, and Step C gives $X^\infty\in B_G(e,r_1)^N$. Applying the preceding estimate with $Y=X^\infty$, we obtain
\[
\kappa^{-1}\max_{i\in[N]} \left| (dR_{x_i(t)^{-1}})_{x_i(t)}\dot x_i(t) \right| \le 3L_0\max_{j\in[N]}d_G(x_j(t),x_j^\infty)\le 3L_0D_{\mathrm R}e^{-2c_{2,\mathrm R}\kappa t}.
\]
Taking $c=2c_{2,\mathrm R}$ and enlarging $C$ so that $C\ge(1+3L_0)D_{\mathrm R}$, we obtain assertion (3).
\end{proof}

\begin{corollary}[Identical intrinsic frequencies]\label{C:lie-identical}
Suppose that $\phi(e)=0$ and \eqref{LE-symmetry-assumptions} hold, and suppose that every eigenvalue of $(d\phi)_e$ has positive real part. Let $0<r_1<r_2$, $C>0$, and $c>0$ be as in Theorem~\ref{T:lie-relaxation}. Let $H_i=H$ for every $i\in[N]$ for a fixed $H\in\mathfrak{g}$. If
\[
\kappa>C\|\mathfrak{H}\|_\Gamma,\quad X^0\in B_G(e,r_1)^N,
\]
then the following assertions hold.
\begin{enumerate}
\item The solution of \eqref{LE-reduced} exists uniquely for all $t\geq0$ and remains in $B_G(e,r_2)^N$.
\item There exists $x_*^\infty\in B_G(e,r_1)$ such that $X^\infty=(x_*^\infty,\ldots,x_*^\infty)$ is the unique equilibrium of \eqref{LE-reduced} in $B_G(e,r_1)^N$, and
\[
d_G(x_*^\infty,e)\le C\frac{\|\mathfrak{H}\|_\Gamma}{\kappa}.
\]
\item If $X(t)$ is a solution to \eqref{LE-reduced} with initial data $X^0\in B_G(e,r_1)^N$, then
\[
\max_i d_G(x_i(t),x_*^\infty)+\kappa^{-1}\max_i\left|(dR_{x_i(t)^{-1}})_{x_i(t)}\dot x_i(t)\right|\leq Ce^{-c\kappa t},\quad t\ge0.
\]
\item If in addition $H=0$, then $x_*^\infty=e$.
\end{enumerate}
\end{corollary}
\begin{proof}
	When $H_i=H$ for all $i$, equation \eqref{LE-reduced} is invariant
	under permutations of the oscillator indices. Hence every permutation
	of $X^\infty$ is also an equilibrium in $B_G(e,r_1)^N$.
	By uniqueness in Theorem~\ref{T:lie-relaxation}, each transposition
	fixes $X^\infty$, so all its components coincide. The stated bounds
	now follow from the same theorem. If $H=0$, then $\phi(e)=0$ implies
	that $(e,\ldots,e)$ is an equilibrium, and uniqueness gives $x_*^\infty=e$.
\end{proof}

To obtain the separable Winfree interaction used in this paper, we impose the additional identity
\begin{equation}\label{LE-factorization}
\mathcal K_\Gamma(x,y)=a(y)\phi(x^{-1}),\quad
a\in C^1(G;\mathbb R).
\end{equation}
Then, system \eqref{LE-reduced} becomes
\begin{equation}\label{LE-winfree}
\dot x_i=F_i(x_i)+\frac{\kappa}{N}\sum_{j=1}^N I(x_j)S(x_i),
\end{equation}
with
\begin{equation}\label{LE-fields}
F_i(x)=(dR_x)_eH_i,\quad
I(x)=1+a(x),\quad
S(x)=(dR_x)_e\phi(x^{-1}).
\end{equation}

The factorization \eqref{LE-factorization} is an additional structural property: it does not follow merely from the local attraction assumptions on $\phi$. For $G=\mathbb S^1$, the two automorphisms are $z\mapsto z$ and $z\mapsto\overline z$, and the sine addition formula gives $a(e^{\mathrm i\theta})=\cos\theta$. Hence this construction recovers \eqref{LE-circle-embedding}.

\subsection{Unitary and special unitary examples}\label{subsec:lie-unitary}
The preceding identities can be realized without restrictions on the intrinsic frequency matrices. \newline

Let $G=\mathrm U(d)$, $d\geq1$, or $G=\mathrm{SU}(d)$, $d\geq2$, and set
\begin{equation}\label{LE-unitary-phi}
\phi_G(U)=P_{\mathfrak g}\left(\frac{U-U^\dagger}{2}\right),\quad P_{\mathfrak u(d)}(A)=A,\quad P_{\mathfrak{su}(d)}(A)=A-\frac{\operatorname{tr}A}{d}I_d.
\end{equation}
These are the standard unitary interaction and its traceless projection, respectively; both occur in the Lie-group framework of \cite{ryoo2025asymptotic}. They satisfy
\[ \phi_G(I_d)=0\quad\mbox{and}\quad(d\phi_G)_{I_d}=\operatorname{id}_{\mathfrak g}. \]
For completeness, we give an explicit finite family of symmetries that yields \eqref{LE-factorization}. Let $e_1,\ldots,e_d$ be the standard orthonormal basis of $\mathbb C^d$. We extend the indices periodically by $e_{k+d}=e_k$. Put $\omega=e^{2\pi\mathrm i/d}$ and define the diagonal matrix $D$ and the cyclic shift matrix $P$ via
\[
De_k=\omega^k e_k,\quad Pe_k=e_{k+1},\quad k\in[d],
\]
and define the mixed shift matrices
\[
W_{rs}=D^rP^s,\quad 0\leq r,s<d.
\]
Note that both $D$ and $P$ are unitary and satisfy $D^d=P^d=I_d$ and $DP=\omega PD$.

\begin{lemma}\label{L:unitary-symmetry-action}
Let
\[
\Gamma=\bigl((\mathbb Z/d\mathbb Z)^2\bigr)\rtimes\mathbb Z/2\mathbb Z,
\]
where the nontrivial element of $\mathbb Z/2\mathbb Z$ sends $(r,s)$ to $(-r,s)$. For $0\leq r,s<d$, define
\begin{equation}\label{LE-unitary-symmetries}
\sigma_{rs,0}(U)=W_{rs}UW_{rs}^\dagger,\quad \sigma_{rs,1}(U)=W_{rs}\overline U W_{rs}^\dagger,
\end{equation}
where $\overline U$ denotes the entrywise complex conjugate of $U$. Then $(r,s,\varepsilon)\cdot U=\sigma_{rs,\varepsilon}(U)$ defines an action of $\Gamma$ on $G$ by Lie group automorphisms.
\end{lemma}
\begin{proof}
Since $W_{rs}$ is unitary, conjugation by $W_{rs}$ and entrywise complex conjugation are smooth group automorphisms of $\mathrm U(d)$. They also preserve $\mathrm{SU}(d)$, since unitary conjugation preserves the determinant and $\det\overline U=\overline{\det U}$. Thus each map in \eqref{LE-unitary-symmetries} is a Lie group automorphism of $G$.

Writing the elements of $\Gamma$ as $(r,s,\varepsilon)$, we have the multiplication law
\[
(r,s,\varepsilon)(r',s',\varepsilon')=\bigl(r+(-1)^\varepsilon r',s+s',\varepsilon+\varepsilon'\bigr),
\]
where the first two entries are taken modulo $d$ and the last entry modulo $2$. The relations $DP=\omega PD$, $\overline D=D^{-1}$, and $\overline P=P$ give
\[
W_{rs}W_{r's'}=\omega^{-sr'}W_{r+r',s+s'},\quad \overline{W_{rs}}=W_{-r,s}.
\]
Since scalar factors disappear under conjugation, we obtain
\[
\sigma_{rs,\varepsilon}\circ\sigma_{r's',\varepsilon'}=\sigma_{r+(-1)^\varepsilon r',s+s',\varepsilon+\varepsilon'}.
\]
Together with $\sigma_{00,0}=\operatorname{id}_G$, this proves the action law.
\end{proof}

The action may be nonfaithful. All sums over $\Gamma$ are taken over its $2d^2$ group elements, including multiplicities when distinct elements induce the same automorphism. Differentiating \eqref{LE-unitary-symmetries} at the identity gives
\[
A_{\sigma_{rs,0}}H=W_{rs}HW_{rs}^\dagger,\quad A_{\sigma_{rs,1}}H=W_{rs}\overline H W_{rs}^\dagger,\quad H\in\mathfrak g.
\]

For the matrix computations below, we use the formulas in \eqref{LE-unitary-symmetries} also for arbitrary $Y\in\mathbb C^{d\times d}$. Let $Y=(Y_{k\ell})\in\mathbb C^{d\times d}$, with the row and column indices interpreted modulo $d$. Then we have
\[
(W_{rs}YW_{rs}^\dagger)_{k\ell}=(D^r P^s Y (P^s)^\dagger (D^r)^\dagger)_{k\ell} = \omega^{r(k-\ell)} (P^s Y (P^s)^\dagger)_{k\ell} = \omega^{r(k-\ell)}Y_{k-s,\ell-s},
\]
and
\[
\frac1d\sum_{r=0}^{d-1}\omega^{r(k-\ell)}=\delta_{k\ell} \coloneqq \begin{cases}
1, & k=\ell, \\
0, & k\ne \ell.
\end{cases}
\]
Consequently, since $s$ runs through all diagonal positions,
\[
\left(\frac1{d^2}\sum_{r,s=0}^{d-1}W_{rs}YW_{rs}^\dagger\right)_{k\ell}=\frac 1 {d^2} \sum_{r, s=0}^{d-1} \omega^{r(k-\ell)}Y_{k-s, \ell-s} = \frac{\delta_{k\ell}}d\sum_{s=0}^{d-1}Y_{k-s,\ell-s}=\frac{\operatorname{tr}Y}{d}\delta_{k\ell}.
\]
Applying the same identity to $\overline Y$ and averaging the two identities, we obtain
\begin{equation}\label{LE-matrix-average}
\frac1{d^2}\sum_{r,s=0}^{d-1}\sigma_{rs,0}(Y)=\frac{\operatorname{tr}Y}{d}I_d, \quad \frac1{2d^2}\sum_{\sigma\in\Gamma}\sigma(Y)=\frac12\left(\frac{\operatorname{tr}Y}{d}+\frac{\operatorname{tr}\overline Y}{d}\right)I_d=\frac{\operatorname{Re}\operatorname{tr}Y}{d}I_d.
\end{equation}
Consequently, for skew-Hermitian $H$, we have
\[
\frac1{2d^2}\sum_{\sigma\in\Gamma}A_\sigma H=\frac{\operatorname{Re}\operatorname{tr}H}{d}I_d=0.
\]
The maps $\phi_G$ are equivariant under all the symmetries in \eqref{LE-unitary-symmetries}, so \eqref{LE-symmetry-assumptions} holds. Indeed, both the skew-Hermitian part and its traceless projection commute with unitary conjugation and entrywise complex conjugation. The expression defining $\phi_G$ extends real-linearly to all complex matrices. Therefore, the relations \eqref{LE-matrix-average} imply
\begin{equation*}\label{LE-matrix-factorization}
\frac1{2d^2}\sum_{\sigma\in\Gamma}\phi_G(\sigma(V)U^{-1})=\phi_G\left(\frac{\operatorname{Re}\operatorname{tr}V}{d}U^{-1}\right)=\frac{\operatorname{Re}\operatorname{tr}V}{d}\phi_G(U^{-1}).
\end{equation*}
Thus, we have
\[ a(U)=d^{-1}\operatorname{Re}\operatorname{tr}U\in[-1,1], \]
and Proposition~\ref{P:lie-subsystem} realizes the corresponding Winfree model as a subsystem of a $4d^2N$-particle Kuramoto--Lohe system with coupling $2\kappa$. For $G=\mathrm U(d)$, the fields in \eqref{LE-fields} are
\[
F_i(U)=H_iU,\quad I(U)=1+\frac{\operatorname{Re}\operatorname{tr}U}{d},\quad S(U)=\frac12(I_d-U^2),
\]
so the reduced equation is exactly \eqref{D-3}, with $H_i=\Omega_i$. For $G=\mathrm{SU}(d)$ and $H_i\in\mathfrak{su}(d)$, it instead becomes
\begin{equation}\label{LE-SU-model}
\begin{cases}
\displaystyle \dot U_i=H_iU_i+\frac{\kappa}{N}\sum_{j=1}^N\left(1+\frac{\operatorname{Re}\operatorname{tr}U_j}{d}\right)S_{\mathrm{SU}}(U_i),\\
\displaystyle S_{\mathrm{SU}}(U)=\frac12(I_d-U^2)+\frac{\operatorname{tr}(U-U^\dagger)}{2d}U.
\end{cases}
\end{equation}
\begin{remark}
The trace correction ensures $S_{\mathrm{SU}}(U)U^{-1}\in\mathfrak{su}(d)$ and hence preserves the constraint $\det U=1$. The ordinary unitary Winfree vector field need not preserve this constraint. When $d=2$, the trace correction vanishes because every $U\in\mathrm{SU}(2)$ has real trace; for higher $d$ it is generally nonzero.
\end{remark}
For both matrix groups, the Frobenius inner product is invariant under the symmetries in \eqref{LE-unitary-symmetries}, so $\|\mathfrak{H}\|_\Gamma=\max_i\|H_i\|_{S_2^d}$, and $(d\phi_G)_{I_d}=\operatorname{id}_{\mathfrak g}$ verifies the local attraction hypothesis. Theorem~\ref{T:lie-relaxation} therefore gives, for models \eqref{D-3} and \eqref{LE-SU-model}, a unique equilibrium near $I_d$ and exponential convergence of the states and speeds for sufficiently concentrated initial data and $\kappa>C\max_i\|H_i\|_{S_2^d}$. The constants and neighborhood may depend on $d$, but not on $N$. For the canonical unitary model, Theorem~\ref{T4.1} additionally treats the larger class of initial data whose spectra avoid $-1$, with a coupling threshold depending on those data.

\section{Winfree oscillators on general manifolds}\label{sec:7}
\setcounter{equation}{0}
In this section, we formulate an abstract Winfree model motivated by the sphere, unitary, and Lie-group constructions of the preceding sections.

Let \((\mathcal M,g)\) be an \(m\)-dimensional Riemannian manifold, and, motivated by the sphere model \eqref{B-3}, $\mathrm{SO}(n)$-model \eqref{B-5}, $\mathrm{U}(d)$-models \eqref{D-3} and \eqref{E-1}, the abstract Lie group model \eqref{LE-reduced}, and $\mathrm{SU}(d)$-model \eqref{LE-SU-model}, we consider the Cauchy problem for the abstract  Winfree manifold model:
\begin{equation}\label{F-1}
\begin{cases}
\displaystyle \dot{x}_i=F_i(x_i)+\frac{\kappa}{N}\sum_{j=1}^N I(x_j)S(x_i),\quad t>0,\\
\displaystyle x_i\big|_{t=0}=x_i^0,\quad i\in[N],
\end{cases}
\end{equation}
where $F_1,\ldots,F_N$ and $S$ are Lipschitz vector fields on $\mathcal M$, and $I:\mathcal M\to[0,\infty)$ is a Lipschitz real-valued influence function. We distinguish a point $p\in\mathcal M$ toward which the interaction field $S$ is attractive. The sphere, unitary, and special unitary models above are special cases of \eqref{F-1} with specific choices of $F_i$, $I$, and $S$: 
\vspace{0.1cm}
\begin{itemize}
\item
Example A~(The sphere model): For \(\mathcal M=\mathbb{S}^{d-1}\) and a fixed \(p\in\mathbb{S}^{d-1}\), we set
\begin{equation}\notag
F_i(x)=A_ix,\quad I(x)=1+\langle x,p\rangle,\quad S(x)=p-\langle x,p\rangle x.
\end{equation}
Then system $\eqref{F-1}_1$ reduces to the Winfree sphere model \eqref{B-3}:
\[  \dot{x}_i = A_i x_i+\frac \kappa {N} \sum_{j=1}^N (1+\langle x_j,p\rangle)(p-\langle x_i,p\rangle x_i). \]
\item
Example B (The unitary model):
For \(\mathcal M=\mathrm{U}(d)\), we take \(p=I_d\) and set
\begin{equation}\notag
F_i(U)=\Omega_iU,\quad
I(U)=\frac{d+\operatorname{Re}\operatorname{tr}U}{d},\quad
S(U)=\frac{1}{2}(I_d-U^2).
\end{equation}
Indeed, we multiply \eqref{D-5}  on the right by \(U_i\) to find 
\[ \dot{U}_i=F_i(U_i)+\kappa RS(U_i), \]
so system $\eqref{F-1}_1$ recovers the mean-field form of the Winfree matrix model.
\vspace{0.2cm}
\item
Example C (The special unitary model):
For $\mathcal M=\mathrm{SU}(d)$, $p=I_d$, and $H_i\in\mathfrak{su}(d)$, take
\[
F_i(U)=H_iU,\quad
I(U)=1+\frac{\operatorname{Re}\operatorname{tr}U}{d},\quad
S(U)=S_{\mathrm{SU}}(U),
\]
where $S_{\mathrm{SU}}$ is given in \eqref{LE-SU-model}. Then \eqref{F-1} is the Winfree--Lohe subsystem constructed in Section~\ref{sec:6}.
\vspace{0.2cm}
\item
Example D (The Lie-group model):
Under the factorization assumption \eqref{LE-factorization}, the reduced Lie-group model \eqref{LE-reduced} has the form \eqref{F-1} with $p=e$, $I=1+a$, and the vector fields given in \eqref{LE-fields}. For $\phi\not\equiv0$, evaluating \eqref{LE-factorization} at $y=e$ gives $a(e)=1$. Hence, on a sufficiently small coordinate ball about $e$, the influence $I$ is positive, and the coordinate representatives of $F_i$, $S$, and $I$ are Lipschitz.
\end{itemize}

For $r>0$, we set 
\[ B_r^m\coloneqq\{x\in\mathbb R^m:\|x\|<r\}. \]
To obtain a quantitative local relaxation theorem in one fixed coordinate system, we fix a chart
\[
\varphi:\mathcal N\longrightarrow B_1^m\subset\mathbb{R}^m,
\quad \varphi(p)=0.
\]
From this point onward, we write $x_i$ for the coordinate variables and use $F_i$, $S$, and $I$ for the corresponding coordinate representatives. All inner products, norms, and Lipschitz constants below are Euclidean quantities in this fixed chart; consequently, the resulting theorem is local on $\mathcal M$ and chart-dependent. Henceforth,
\[
F_1,\ldots,F_N,S:B_1^m\to\mathbb{R}^m,
\quad
I:B_1^m\to[0,\infty) \quad \mbox{are Lipschitz functions}.
\]
We introduce the order parameter:
\begin{equation}\label{F-2}
R\coloneqq\frac{1}{N}\sum_{j=1}^N I(x_j),
\end{equation}
so that $\eqref{F-1}_1$ becomes
\begin{equation}\label{F-3}
\dot{x}_i=F_i(x_i)+\kappa RS(x_i),\quad i\in[N].
\end{equation}
We also set
\[
{\mathfrak F}\coloneqq(F_1,\ldots,F_N),\quad X\coloneqq(x_1,\ldots,x_N),
\]
and
\[
\|{\mathfrak F}\|_\infty\coloneqq\max_{i\in[N]}\|F_i\|_{L^\infty(B_1^m)},\quad
\|\operatorname{Lip}{\mathfrak F}\|_\infty\coloneqq\max_{i\in[N]}\operatorname{Lip}F_i,
\quad
\|S\|_\infty\coloneqq\|S\|_{L^\infty(B_1^m;\mathbb{R}^m)}.
\]
Next, we impose the following framework for relaxation dynamics in terms of $F_i$, $S$, and $I$:
\begin{itemize}
\item
$({\mathcal F}_{B}1)$: The vector fields \(F_1,\ldots,F_N\) and \(S\), and the influence function \(I\), are Lipschitz continuous on \(B_1^m\).
\vspace{0.1cm}
\item
$({\mathcal F}_{B}2)$: The interaction field \(S\) is strictly attractive toward the origin: there exists \(c>0\) such that
\begin{equation}\label{F-4}
S(0)=0,\quad
\langle S(x)-S(y),x-y\rangle\leq-c\|x-y\|^2,
\quad x,y\in B_1^m.
\end{equation}
\item
$({\mathcal F}_{B}3)$: The influence function is strictly positive: there exists \(I_0>0\) such that
\begin{equation}\label{F-5}
I(x)\geq I_0,
\quad x\in B_1^m.
\end{equation}
\item
$({\mathcal F}_{B}4)$: The attraction dominates the variation of the influence function, and the coupling strength is sufficiently large:
\begin{equation}\notag
cI_0>\|S\|_\infty\operatorname{Lip}I,
\quad
\kappa>\max\left\{
\frac{\|{\mathfrak F}\|_\infty}{cI_0},
\frac{\|\operatorname{Lip}{\mathfrak F}\|_\infty}{cI_0-\|S\|_\infty\operatorname{Lip}I}
\right\}.
\end{equation}
\end{itemize}

\subsection{Global well-posedness}\label{subsec:6-1}
In this subsection, we discuss the global well-posedness of the Cauchy problem \eqref{F-1}. 

\begin{lemma}\label{L6.1}
Suppose that the assumptions $({\mathcal F}_{B}1)$--$({\mathcal F}_{B}4)$ hold and that initial data satisfy 
\[
x_i^0\in B_1^m,\quad i\in[N].
\]
Then the Cauchy problem \eqref{F-1} admits a unique global solution that is Lipschitz continuous in time.
\end{lemma}

\begin{proof}
The proof consists of two steps. We first choose a smaller ball containing the initial data and then we show that no component can leave that ball. \newline

\noindent\textbf{$\bullet$ Step A} (Choice of an invariant ball): We choose \(r\in(0,1)\) sufficiently close to \(1\) such that
\begin{equation}\notag
x_i^0\in B_r^m,\quad i\in[N],
\quad
\kappa>\frac{\|{\mathfrak F}\|_\infty}{crI_0}.
\end{equation}
The Cauchy--Lipschitz theorem gives a unique local solution as long as \(X(t)\in(B_1^m)^N\).

\vspace{0.5em}

\noindent\textbf{$\bullet$ Step B} (First-exit contradiction):
Suppose that \(T<\infty\) is the first time at which some component reaches \(\partial B_r^m\). Then, we have
\[ \|x_i(T)\|=r \quad \mbox{for some \(i\in[N]\)}. \]
By \eqref{F-3}, \eqref{F-4}, and \eqref{F-5}, we have
\begin{align*}
\frac{1}{2}\left.\frac{d}{dt}\|x_i(t)\|^2\right|_{t=T}
&=\langle F_i(x_i(T)),x_i(T)\rangle
+\kappa R(T)\langle S(x_i(T)),x_i(T)\rangle\\
&\leq \|{\mathfrak F}\|_\infty r-\kappa cI_0r^2\\
&=r\left(\|{\mathfrak F}\|_\infty-\kappa crI_0\right)<0,
\end{align*}
where we used \(S(0)=0\) in the attraction estimate. At the first exit time, however, the left derivative of \(\|x_i(t)\|^2\) cannot be negative. This contradiction shows that every component remains in \(B_r^m\). Since \(\overline{B_r^m}\) is compactly contained in \(B_1^m\), the local solution extends for all \(t\geq0\). Moreover, the Lipschitz assumptions on the bounded ball imply that the vector fields and influence function are bounded, and
\[
\|\dot x_i(t)\|\leq\|\mathfrak F\|_\infty+\kappa\|I\|_{L^\infty(B_1^m)}\|S\|_\infty,
\quad t\geq0,\quad i\in[N].
\]
Thus the solution is Lipschitz continuous in time, as asserted.
\end{proof}
\subsection{Exponential stability}\label{subsec:6-2}
In this subsection, we discuss the exponential stability of \eqref{F-1}.  For an $N$-tuple $Y=(y_1,\ldots,y_N)\in(\mathbb{R}^m)^N$, we use the mixed norm $\|Y\|_{q,2}$ defined in \eqref{A-2}.
We also set
\begin{equation}\notag
\mu_{\mathcal M}\coloneqq
\kappa\left(cI_0-\|S\|_\infty\operatorname{Lip}I\right)
-\|\operatorname{Lip}{\mathfrak F}\|_\infty.
\end{equation}
Then, by $({\mathcal F}_{B}4)$,  we have 
\[ \mu_{\mathcal M}>0. \]
\begin{lemma}\label{L6.2}
Suppose that $({\mathcal F}_{B}1)$--$({\mathcal F}_{B}4)$ hold. Let \(X=(x_1,\ldots,x_N)\) and \(\widetilde{X}=(\widetilde{x}_1,\ldots,\widetilde{x}_N)\) be two global solutions to \eqref{F-1} with
\(X(0),\widetilde X(0)\in(B_1^m)^N\). Then, for every \(1\leq q\leq\infty\),
\begin{equation}\label{F-6}
\|X(t)-\widetilde{X}(t)\|_{q,2}
\leq e^{-\mu_{\mathcal M}t}
\|X(0)-\widetilde{X}(0)\|_{q,2},
\quad t\geq0.
\end{equation}
\end{lemma}
\begin{proof}
The proof consists of two steps. We first derive a pointwise differential inequality for $X - \widetilde{X}$, and we then pass to the mixed norm. \newline

\noindent\textbf{$\bullet$ Step A} (Pointwise estimate):
Let \(R\) and \(\widetilde{R}\) be the order parameters associated with \(X\) and \(\widetilde{X}\), respectively. Subtracting the two equations in \eqref{F-3}, we obtain
\begin{equation}\label{F-7}
\dot{x}_i-\dot{\widetilde{x}}_i
=F_i(x_i)-F_i(\widetilde{x}_i)
+\kappa R\big(S(x_i)-S(\widetilde{x}_i)\big)
+\kappa(R-\widetilde{R})S(\widetilde{x}_i).
\end{equation}
Moreover, we use \eqref{F-2} to find 
\begin{equation}\label{F-8}
|R-\widetilde{R}|
\leq\frac{\operatorname{Lip}I}{N}\sum_{j=1}^N\|x_j-\widetilde{x}_j\|
=(\operatorname{Lip}I)\|X-\widetilde{X}\|_{1,2}.
\end{equation}
Now, we take the inner product of \eqref{F-7} with \(x_i-\widetilde{x}_i\), and use \eqref{F-4}, \eqref{F-5}, and \eqref{F-8} to find
\[
\frac{1}{2}\frac{d}{dt}\|x_i-\widetilde{x}_i\|^2 \leq\left(\|\operatorname{Lip}{\mathfrak F}\|_\infty-\kappa cI_0\right)\|x_i-\widetilde{x}_i\|^2 +\kappa\|S\|_\infty(\operatorname{Lip}I)
\|X-\widetilde{X}\|_{1,2}\|x_i-\widetilde{x}_i\|.
\]
For $x_i\neq\widetilde x_i$, we divide the squared-distance estimate by $\|x_i-\widetilde x_i\|$. At $x_i=\widetilde x_i$, the relation \eqref{F-7} gives
\[
D^+\|x_i-\widetilde x_i\|=\|\dot x_i-\dot{\widetilde x}_i\|
\leq\kappa\|S\|_\infty(\operatorname{Lip}I)\|X-\widetilde X\|_{1,2}.
\]
Thus, at every time,
\begin{equation}\label{F-9}
D^+\|x_i-\widetilde{x}_i\|
\leq\left(\|\operatorname{Lip}{\mathfrak F}\|_\infty-\kappa cI_0\right)\|x_i-\widetilde{x}_i\|
+\kappa\|S\|_\infty(\operatorname{Lip}I)
\|X-\widetilde{X}\|_{1,2}.
\end{equation}

\vspace{0.5em}
\noindent\textbf{$\bullet$ Step B} (Mixed-norm estimate):
For $q=1$, we average \eqref{F-9}. For $1<q<\infty$, at times with $\|X-\widetilde X\|_{q,2}>0$, we multiply \eqref{F-9} by $\|x_i-\widetilde x_i\|^{q-1}$, average the resulting relation over $i$, and multiply by $\|X-\widetilde X\|_{q,2}^{1-q}$. Using H\"older's inequality and
\[
\|X-\widetilde X\|_{1,2}\leq\|X-\widetilde X\|_{q,2},
\]
we obtain
\[
D^+\|X-\widetilde{X}\|_{q,2} \leq\left(\|\operatorname{Lip}{\mathfrak F}\|_\infty-\kappa cI_0
+\kappa\|S\|_\infty\operatorname{Lip}I\right)
\|X-\widetilde{X}\|_{q,2} =-\mu_{\mathcal M}\|X-\widetilde{X}\|_{q,2}.
\]
If $\|X-\widetilde X\|_{q,2}=0$, uniqueness of solutions gives $X=\widetilde X$ from that time onward, so the same inequality holds. The case \(q=\infty\) follows from the upper right Dini derivative rule for the finite maximum of the component distances. Gronwall's lemma yields \eqref{F-6}.
\end{proof}

\subsection{Relaxation dynamics}\label{subsec:6-3}
In this subsection, we study the relaxation dynamics of \eqref{F-1} using Lemma \ref{L6.1} and Lemma \ref{L6.2}.

\begin{proposition}\label{P6.1}
Suppose that $({\mathcal F}_{B}1)$--$({\mathcal F}_{B}4)$ hold. Then system \eqref{F-1} admits a unique equilibrium
\[
X^\infty\in(B_1^m)^N.
\]
\end{proposition}

\begin{proof}
Choose arbitrary initial data $X^0=(x_1^0,\ldots,x_N^0)\in(B_1^m)^N$, and let $X(t)$ be the corresponding unique global solution given by Lemma~\ref{L6.1}. By its proof, there exists $r\in(0,1)$ such that 
\[ X(t)\in(B_r^m)^N \quad \mbox{for all $t\geq0$}. \]
For $s\geq0$, the time translate $\widetilde{X}(t)\coloneqq X(t+s)$ is another global solution. Lemma~\ref{L6.2} with $q=2$ implies
\begin{equation}\label{F-10}
\|X(t+s)-X(t)\|_{2,2}
\leq e^{-\mu_{\mathcal M}t}\|X(s)-X(0)\|_{2,2}.
\end{equation}
Because $X(t)\in(B_r^m)^N$ and $\overline{B_r^m}$ is compact, the last factor in \eqref{F-10} is bounded uniformly in $s$. For $t_2\geq t_1\geq0$, we take $t=t_1$ and $s=t_2-t_1$ in \eqref{F-10}. Then, we have
\[
\|X(t_2)-X(t_1)\|_{2,2}
\leq Ce^{-\mu_{\mathcal M}t_1},
\]
where $C$ is independent of $t_1$ and $t_2$. Thus $X(t)$ is Cauchy and it converges to some $X^\infty\in(\overline{B_r^m})^N$.

Let \(\Phi_s\) denote the flow generated by \eqref{F-1}. By the continuity of the flow,
\[
\Phi_s(X^\infty)
=\lim_{t\to\infty}\Phi_s(X(t))
=\lim_{t\to\infty}X(t+s)
=X^\infty,
\quad s\geq0.
\]
Thus \(X^\infty\) is an equilibrium. The strict inequality in the boundary estimate from the proof of Lemma~\ref{L6.1} rules out \(\|x_i^\infty\|=r\), so
\[ X^\infty\in(B_r^m)^N\subset(B_1^m)^N. \]

If \(Y^\infty\) is another equilibrium, Lemma~\ref{L6.2} applied to the two constant solutions gives
\[
\|X^\infty-Y^\infty\|_{2,2}
\leq e^{-\mu_{\mathcal M}t}\|X^\infty-Y^\infty\|_{2,2},
\quad t\geq0.
\]
Letting \(t\to\infty\) yields 
\[ X^\infty=Y^\infty. \]
\end{proof}

\begin{theorem}\label{T6.1}
Suppose that $({\mathcal F}_{B}1)$--$({\mathcal F}_{B}4)$ hold, and let \(X(t)\) be the global solution to \eqref{F-1} with initial data
\[
x_i^0\in B_1^m,\quad i\in[N].
\]
Then \(X(t)\) converges exponentially to the unique equilibrium \(X^\infty\) from Proposition~\ref{P6.1}. More precisely, for every \(1\leq q\leq\infty\),
\begin{equation}\label{F-11}
\|X(t)-X^\infty\|_{q,2}
\leq e^{-\mu_{\mathcal M}t}\|X(0)-X^\infty\|_{q,2},
\quad t\geq0.
\end{equation}
\end{theorem}
\begin{proof}
The equilibrium \(X^\infty\) generates the constant solution \(\widetilde{X}(t)\equiv X^\infty\), and we apply Lemma~\ref{L6.2} to \(X\) and \(\widetilde{X}\) to derive the desired estimate \eqref{F-11}.
\end{proof}

\section{Concluding remarks}\label{sec:8}
\setcounter{equation}{0}
In this paper, we have proposed higher-dimensional counterparts of the phase-locking results of \cite{ha2015emergence} for the Winfree model on the unit sphere $\mathbb{S}^{d-1}$, the unitary and special unitary groups, Lie groups, and general Riemannian manifolds. In Section~\ref{sec:3}--Section \ref{sec:5}, trajectories enter a positively	invariant region in finite time. Contraction in that region allows us to	compare each trajectory with its time translates and prove exponential convergence to an equilibrium. In Section~\ref{sec:7}, the same contraction argument applies from time zero because the initial data already lie in the invariant coordinate neighborhood. In Section~\ref{sec:6}, we embed the reduced system in a Kuramoto--Lohe ensemble with an oscillator fixed at the identity. The relative phase of each original oscillator with respect to this fixed oscillator is its state itself, so convergence of the ambient relative phases yields convergence of the reduced states to equilibrium.

Some natural questions remain open. First, how sharp are the sufficient coupling bounds for convergence in Theorems~\ref{T3.1} and~\ref{T4.1}? For prescribed intrinsic frequencies and initial data, define the pathwise critical coupling as the infimum of all $k\geq0$ such that, for every $\kappa>k$, the corresponding solution converges to an equilibrium, as in \cite[Definition~8]{ryoo2026oscillator}. For the canonical unitary model, Remark~\ref{R:unitary-critical-coupling} gives the exact threshold for the existence of an equilibrium. When does the pathwise critical coupling coincide with this threshold, in particular when some $U_i^0$ fails to commute with $\Omega_i$? For the canonical sphere model \eqref{C-1}, one may likewise ask for the equilibrium-existence threshold and its relation to the pathwise critical coupling. Second, it is natural to ask whether the results of this paper admit mean-field limits as $N\to\infty$, yielding convergence  theorems for the corresponding Vlasov-type kinetic equations on $\mathbb{S}^{d-1}$, $\mathrm{U}(d)$, $\mathrm{SU}(d)$, Lie groups, and general manifolds.

\appendix
\section{Proof of Lemma \ref{L4.1}}\label{app:1}
In this appendix, we verify Lemma~\ref{L4.1}. We thank Otte Hein\"avaara for the following proof. Note that if $A$ is a $d\times d$ matrix with eigenvalues $\lambda_1,\cdots,\lambda_d$, we have 
\[ \sum_{i=1}^d |\lambda_i|\le \|A\|_{S_1^d}. \]
The proof goes as follows.\footnote{This argument is from \url{https://math.stackexchange.com/questions/381808/sum-of-eigenvalues-and-singular-values}.} We assume that $A$ is upper triangular by Schur decomposition, with $\lambda_1,\cdots,\lambda_d$ on the diagonal. Then, we write
\[ \lambda_i= |\lambda_i|\alpha_i \quad \mbox{with $|\alpha_i|=1$} \]
and define 
\[ D=\operatorname{diag}(\alpha_1,\cdots,\alpha_d)\in \mathrm{U}(d) \quad \mbox{and} \quad  B=D^{-1}A. \]
Then, we have 
\[ \operatorname{tr}B=\sum_{i=1}^d|\lambda_i|. \] On the other hand, $B$ has the same singular values $\sigma_1,\cdots,\sigma_d$ as $A$, and 
\[ B=U\operatorname{diag}(\sigma_1,\cdots,\sigma_d)V \quad \mbox{for some $U,V\in \mathrm{U}(d)$}. \]
Now
    \begin{align*}
    \begin{aligned}
    \sum_{i=1}^d|\lambda_i| &=\operatorname{tr}B=\operatorname{tr}(U\operatorname{diag}(\sigma_1,\cdots,\sigma_d)V) =\operatorname{tr}(\operatorname{diag}(\sigma_1,\cdots,\sigma_d)VU) \\
    &=|\sum_{i=1}^d \sigma_i[VU]_{ii}|\le \sum_{i=1}^d \sigma_i|[VU]_{ii}| \le \sum_{i=1}^d\sigma_i=\|B\|_{S_1^d}=\|A\|_{S_1^d},
    \end{aligned}
    \end{align*}
    where $[\cdot]_{ij}$ denotes the $(i,j)$-component. \newline

    Because of this fact, it is enough to show that
\[ \operatorname{Re}\lambda\le \sin\alpha^\infty|\lambda|, \]
for each eigenvalue $\lambda$ of $U-V$.
Equivalently, we show that if a complex number $\lambda$ satisfies 
\[ \operatorname{Re}\lambda> \sin\alpha^\infty|\lambda|, \]
then $\lambda$ is not an eigenvalue of $U-V$. By \cite[Lemma 2]{wielandt1955eigenvalues}, if we denote the eigenvalues of $U$ by $\alpha_1,\cdots,\alpha_d$ and those of $V$ by $\beta_1,\cdots,\beta_d$, it is enough to show that there is a circle separating $\alpha_1,\cdots,\alpha_d$ from $\beta_1+\lambda,\cdots,\beta_d+\lambda$. We set
\[ \Gamma\coloneqq\{e^{i\theta}:-\alpha^\infty<\theta<\alpha^\infty\}. \]
Since $\alpha_i,\beta_i\in \Gamma$, it is enough to show that there is a circle that separates $\Gamma$ from $\Gamma+\lambda$. But for any $e^{i\theta}\in \Gamma$, it is clear that $\operatorname{Re}(\lambda e^{-i\theta})> 0$, so that 
\[ |e^{i\theta}+\lambda|^2 > 1+|\lambda|^2. \]
So the circle $|z|=\sqrt{1+|\lambda|^2}$ is a circle separating $\Gamma$ from $\Gamma+\lambda$, completing the proof.

\section{Proof of Lemma \ref{L4.2}}\label{app:2}
In this appendix, we verify Lemma~\ref{L4.2}. Since $U$ is unitary, $K$ and $H$ commute:
\[H K = \frac 1 4 (U + U^\dagger)(U - U^\dagger) = \frac 1 4 (U^2 - (U^\dagger)^2) = \frac 1 4 (U - U^\dagger)(U + U^\dagger) =  K H.\]
Then, for each $x\in \bbc^d$, we have
\[\|x\|^2 = \|Ux\|^2 = \|(H + K)x\|^2 = \|Hx\|^2 + \|Kx\|^2 + x^\dagger (K^\dagger H + H^\dagger K) x.\]
Since $H$ and $K$ are Hermitian and skew-Hermitian, respectively, we have
\[K^\dagger H + H^\dagger K = -KH + HK = 0.\]
This implies the desired result:
\[\|x\|^2 = \|Hx\|^2 + \|Kx\|^2.\]

\bibliographystyle{myalpha}
\bibliography{bib}


\end{document}